\documentclass[twoside]{article}

\usepackage[accepted]{aistats2020}

\usepackage{epsf}
\usepackage{graphics}
\usepackage{graphicx}
\usepackage{psfrag}
\usepackage{microtype}
\usepackage{subfigure}
\usepackage{algorithmic}
\usepackage{color,xcolor}
\usepackage[linesnumbered,ruled]{algorithm2e}% http://ctan.org/pkg/algorithm2e
\usepackage{color}

\usepackage{amsthm}
\usepackage{amsfonts}
\usepackage{amsmath}
\usepackage{amssymb,bbm}
\usepackage{multirow}

\usepackage{etoc}

\usepackage{url}% for url's in bib
\usepackage[colorlinks=True,citecolor=magenta,linkcolor=blue,urlcolor=blue,pagebackref=true,backref=true]
{hyperref}
\renewcommand*{\backrefalt}[4]{%
    \ifcase #1 \footnotesize{(Not cited.)}%
    \or        \footnotesize{(Cited on page~#2.)}%
    \else      \footnotesize{(Cited on pages~#2.)}%
    \fi}

\usepackage{nicefrac}
\usepackage{comment}

 \usepackage{tabularx}%

\usepackage{enumitem}
\usepackage{booktabs}
\usepackage{caption}

\usepackage{bm,bbm}
\usepackage{mathtools}

\newcommand{\mixfit}[1]{\mathcal{G}_{\mathrm{#1}}}

\newcommand{\mix}{\mathcal{G}}

\newcommand{\lstar}{\ind_\star}
\newcommand{\scalartwo}{b}
\newcommand{\remainder}{\varepsilon}
\newcommand{\myfrac}[2]{\ensuremath{\frac{#1}{#2}}}
\newcommand{\Mtil}{\widetilde{M}}
\newcommand{\PseudoND}{\Mtil_{n,d}}

\newcommand{\qtext}[1]{\quad\text{#1}\quad}
\newcommand{\radem}{\varepsilon}

\newcommand{\event}{\mathcal{E}}
\newcommand{\staterr}{\omega}
\newcommand{\newm}{\widehat{M}}
\newcommand{\newtheta}{\widetilde\theta}
\newcommand{\newr}{\widetilde{r}}

\newcommand{\dims}{\ensuremath{d}}
\newcommand{\samples}{\ensuremath{n}}
\newcommand{\obs}{\samples}
\newcommand{\real}{\ensuremath{\mathbb{R}}}
\newcommand{\realdim}{\ensuremath{\real^\dims}}

\newcommand{\smallthreshold}{\beta}

\newcommand{\order}[1]{\ensuremath{\mathcal{O}\parenth{#1}}}

\newcommand{\thetastar}{\ensuremath{\theta^*}}
\newcommand{\sdstar}{\ensuremath{\sigma^*}}

\newcommand{\NORMAL}{\ensuremath{\mathcal{N}}}

\newcommand{\brackets}[1]{\left[ #1 \right]}
\newcommand{\parenth}[1]{\left( #1 \right)}

\newcommand{\biggparenth}[1]{\bigg( #1 \bigg)}
\newcommand{\braces}[1]{\left\{ #1 \right \}}
\newcommand{\abss}[1]{\left| #1 \right |}

\newcommand{\tp}{^\top}
\newcommand{\ceil}[1]{\left\lceil #1 \right\rceil}

\newcommand{\interval}{{I}}

\newcommand{\term}{U}

\newcommand{\basis}{e}
\newcommand{\locacontramutilo}{\ensuremath{\underline\rho}}
\newcommand{\locacontramutiup}{\ensuremath{\overline\rho}}

\newcommand{\unicontraup}{\ensuremath{\overline\gamma}}
\newcommand{\unicontralo}{\ensuremath{\underline\gamma}}
\newcommand{\truncate}{\ensuremath{\tau}}
\newcommand{\corop}{\ensuremath{\overline{M}_{1}}}
\newcommand{\gencorop}{\ensuremath{\overline{M}_{d}}}

\newcommand{\ind}{\ensuremath{\ell}}

\newcommand{\samm}{M_{n,1}}
\newcommand{\popmean}[1]{\mu_{#1}}
\newcommand{\sammean}[1]{\widehat\mu_{#1}}
\newcommand{\idx}{k}
\newcommand{\constant}{c}
\newcommand{\simconst}{\beta}
\newcommand{\simconsthat}{\widehat\simconst}
\newcommand{\jset}{\mathcal{J}}

\newcommand{\Rspace}{\ensuremath{\mathbb{R}}}

\newcommand{\Ncal}{\ensuremath{\mathcal{N}}}

\newcommand{\ball}{\ensuremath{\mathbb{B}}}

\newcommand{\estloca}{\ensuremath{\widehat{\theta}_{n}}}
\newcommand{\estsca}{\ensuremath{\widehat{\sigma}_{n}}}

\newcommand{\normDensity}{\ensuremath{\phi}}
\newcommand{\mean}{\ensuremath{\theta}}
\newcommand{\sd}{\ensuremath{\sigma}}

\newcommand{\weightFun}{\ensuremath{w}}
\newcommand{\mydefn}{\ensuremath{:=}}

\newcommand{\samsurro}{\ensuremath{Q_{n}}}

\newcommand{\unicon}{c}
\newcommand{\uniconnew}{c'}

\newcommand{\unicontwo}{c_2}

\newcommand{\defn}{:=}
\newcommand{\rdefn}{=:}
\newcommand{\etal}{{et al.}}

\newcommand{\matsnorm}[2]{|\!|\!| #1 | \! | \!|_{{#2}}}
\newcommand{\vecnorm}[2]{\| #1 \|_{#2}}

\newcommand{\enorm}[1]{\vecnorm{#1}{2}} % euclidean norm

\newcommand{\fronorm}[1]{\ensuremath{\matsnorm{#1}{\footnotesize{\mathbb{F}}}}}

\newcommand{\Exs}{\ensuremath{{\mathbb{E}}}}
\newcommand{\Prob}{\ensuremath{{\mathbb{P}}}}

\newtheoremstyle{named}{}{}{\itshape}{}{\bfseries}{.}{.5em}{\thmnote{#3's }#1}
\theoremstyle{named}

\theoremstyle{plain}

\newtheorem{theorem}{Theorem}
\newtheorem{proposition}{Proposition}
\newtheorem{lemma}{Lemma}

\newtheorem{corollary}{Corollary}

\newlength{\widebarargwidth}
\newlength{\widebarargheight}
\newlength{\widebarargdepth}

\makeatletter
\long\def\@makecaption#1#2{
        \vskip 0.8ex
        \setbox\@tempboxa\hbox{\small {\bf #1:} #2}
        \parindent 1.5em  %% How can we use the global value of this???
        \dimen0=\hsize
        \advance\dimen0 by -3em
        \ifdim \wd\@tempboxa >\dimen0
                \hbox to \hsize{
                        \parindent 0em
                        \hfil
                        \parbox{\dimen0}{\def\baselinestretch{0.96}\small
                                {\bf #1.} #2
                                }
                        \hfil}
        \else \hbox to \hsize{\hfil \box\@tempboxa \hfil}
        \fi
        }
\makeatother

\long\def\comment#1{}
\definecolor{battleshipgrey}{rgb}{0.52, 0.52, 0.51}
\definecolor{darkgray}{rgb}{0.66, 0.66, 0.66}
\definecolor{darkgreen}{rgb}{0.0, 0.2, 0.13}
\definecolor{darkspringgreen}{rgb}{0.09, 0.45, 0.27}
\definecolor{dukeblue}{rgb}{0.0, 0.0, 0.61}
\definecolor{olivedrab7}{rgb}{0.24, 0.2, 0.12}
\definecolor{darkblue}{rgb}{0.0, 0.0, 0.55}
\definecolor{darkscarlet}{rgb}{0.34, 0.01, 0.1}
\definecolor{candyapplered}{rgb}{1.0, 0.03, 0.0}
\definecolor{ao(english)}{rgb}{0.0, 0.5, 0.0}
\definecolor{applegreen}{rgb}{0.55, 0.71, 0.0}

\newcommand{\widgraph}[2]{\includegraphics[keepaspectratio,width=#1]{#2}}

\begin{document}

\etocdepthtag.toc{mtchapter}
\etocsettagdepth{mtchapter}{subsection}
\etocsettagdepth{mtappendix}{none}

\twocolumn[

\aistatstitle{Sharp Analysis of Expectation-Maximization\\
for Weakly Identifiable Models}
        
\aistatsauthor{ Raaz Dwivedi$^{\dagger,\star}$ \And Koulik Khamaru$^{\Diamond,\star}$
\And Nhat Ho$^{\dagger,\star}$}
\aistatsauthor{Martin J. Wainwright$^{\dagger,\Diamond}$ \And 
Michael I. Jordan$^{\dagger,\Diamond}$ \And Bin Yu$^{\dagger,\Diamond}$}

\aistatsaddress{ Department of $^\dagger$EECS and $^\Diamond$Statistics,
UC Berkeley} ]
% \aistatsauthor{ Author 1 \And Author 2 \And  Author 3 }

% \aistatsaddress{ Institution 1 \And  Institution 2 \And Institution 3 } 

%%\vspace*{.2in}
% {\large{
% \begin{tabular}{ccc}
%  Raaz Dwivedi$^{\star, \diamond}$ & Nhat Ho$^{\star, \diamond}$ &  Koulik
%  Khamaru$^{\star, \dagger}$ \\
% \end{tabular}
% \begin{tabular}
% {ccc}
%  Martin J. Wainwright$^{\diamond, \dagger,\ddagger}$ & Michael I. Jordan$^{\diamond, \dagger}$ & 
%  Bin Yu$^{\diamond, \dagger}$
% \end{tabular}
%}}

% %\vspace*{.2in}

% \begin{tabular}{c}
% Department of Electrical Engineering and Computer Sciences$^\diamond$\\
% Department of Statistics$^\dagger$ \\
% UC Berkeley\\
% \end{tabular}

% %\vspace*{.1in}
% \begin{tabular}{c}
% The Voleon Group$^\ddagger$
% \end{tabular}

%%\vspace*{.2in}

%\today

%%\vspace*{.2in}

\begin{abstract}
 We study a class of weakly identifiable location-scale mixture models for which the maximum likelihood estimates based on $n$ i.i.d. samples are known to have lower accuracy than the classical  $n^{- \frac{1}{2}}$ error. We investigate whether the  Expectation-Maximization (EM) algorithm also converges slowly for these models.  We provide a rigorous characterization of EM for fitting a weakly identifiable Gaussian mixture in a univariate setting where we prove that the EM algorithm converges in order $n^{\frac{3}{4}}$ steps and returns estimates that are at a Euclidean distance of order ${ n^{- \frac{1}{8}}}$ and ${ n^{-\frac{1} {4}}}$ from the true location and scale parameter respectively. Establishing the slow rates in the univariate setting requires a novel localization argument with two stages, with each stage involving an epoch-based argument applied to a different surrogate EM operator at the population level.  We demonstrate several multivariate ($d \geq 2$) examples that exhibit the same slow rates as the univariate case. We also prove slow statistical rates in higher dimensions in a special case, when the fitted covariance is constrained to be a multiple of the identity.
\end{abstract}
%\let\thefootnote\relax\footnotetext{$^\star$Raaz Dwivedi, Nhat Ho,
%  and Koulik Khamaru contributed equally to this work.}
%\end{center}

%\vspace{-2mm}
\section{Introduction} 
\label{sec:introduction}
%\vspace{-2mm}
Gaussian mixture models~\cite{Pearson-1894} have been used widely to
model heterogeneous data in many applications arising from physical
and the biological sciences. In several scenarios, the data has a large
number of sub-populations and the mixture components in the data may
not be well-separated. In such settings, estimating the true number of
components may be difficult, so that one may end up fitting a mixture model with a number of components larger than that present in the data.  Such mixture fits, referred to as \emph{over-specified
  mixture distributions}, are commonly used by practitioners in order
to deal with uncertainty in the number of components in the
data~\cite{Rousseau-2011, Rousseau-2015}.  However, a deficiency of
such models is that they are \emph{singular}, meaning that their
Fisher information matrices are degenerate.  Given the popularity of
over-specified models in practice, it is important to understand how
methods for parameter estimation, including maximum likelihood and the
EM algorithm, behave when applied to such models.
%\vspace{-3mm}

\subsection{Background and past work} % (fold)
\label{sub:background_and_past_work}
%\vspace{-2mm}
In the context of singular mixture models, an important distinction is
between those that are \emph{strongly} versus \emph{weakly}
identifiable.  Chen~\cite{Chen1992} studied the class of strongly
identifiable models in which, while the Fisher information matrix may
be degenerate at a point, and it is not degenerate over a larger set.
Studying over-specified Gaussian mixtures with known scale parameters,
he showed that the accuracy of the MLE for the unknown location parameter is of the order $n^{-\frac{1}{4}}$, which should be
contrasted with the classical $n^{-\frac{1}{2}}$ rate achieved in regular settings.  A line of follow-up work has extended this type of analysis to other types of strongly identifiable mixture models; see
the papers~\cite{Ishwaran-2001, Rousseau-2011, Nguyen-13, Jonas-2016}
as well as the references therein for more details.

A more challenging class of mixture models are those that are only
\emph{weakly identifiable}, meaning that the Fisher information is
degenerate over some larger set.  This stronger form of singularity
arises, for instance, when the scale parameter in an over-specified
Gaussian mixture is also unknown~\cite{Chen-jrssb_2001, Chen_2009}. Ho et al.~\cite{Ho-Nguyen-AOS-17} characterized the
behavior of MLE for a class of weakly identifiable models. They showed
that the convergence rates of MLE in these models could be very slow,
with the precise rates determined by algebraic relations among the
partial derivatives.  However, this past work has not addressed the
computational complexity of computing the MLE in a weakly identifiable
model.

The focus of this paper is the intersection of statistical and
computational issues associated with fitting the parameters of weakly
identifiable mixture models.  In particular, we study the
expectation-maximization (EM)
algorithm~\cite{Rubin-1977,Jeff_Wu-1983,redner1984mixture}, which is
the most popular algorithm for computing (approximate) MLEs in the mixture models.  It is an instance of a minorization-maximization algorithm,
in which at each step, a suitably chosen lower bound of the
log-likelihood is maximized. There is now a lengthy line of work on
the behavior of EM when applied to regular models.  The classical
papers~\cite{Jeff_Wu-1983,tseng2004analysis,chretien2008algorithms}
establish the asymptotic convergence of EM to a local maximum of the
log-likelihood function for a general class of incomplete data models.
Other papers~\cite{Xu_Jordan-1995,Jordan-1996,Jordan-2000}
characterized the rate of convergence of EM for regular Gaussian
mixtures.  More recent years have witnessed a flurry of work on the behavior of EM for various kinds of regular mixture
models~\cite{Siva_2017, HanLiu_nips2015, Caramanis-nips2015,
  Hsu-nips2016, Daskalakis_colt2017, Sarkar_nips2017, Cheng_2018,
  Cai_2018}; as a consequence, our understanding of EM in such cases
is now relatively mature.  More precisely, it is known that for
Gaussian mixtures, EM converges in $\mathcal{O}(\log(n/d))$-steps to
parameter estimates that lie within Euclidean distance
$\mathcal{O}((d/n)^ {1/2})$ of the true location parameters, assuming minimal separation between the mixture components.

In our recent work \cite{Raaz_Ho_Koulik_2018}, we studied the behavior of EM for fitting a class of \emph{non-regular} mixture
models, namely those in which the Fisher information is degenerate at
a point, but the model remains strongly identifiable.  One such class
of models are Gaussian location mixtures with known scale parameters
that are \emph{over-specified}, meaning that the number of components
in the mixture-fit exceeds the number of components in the data
generating distribution. For such non-regular but strongly
identifiable mixture models, they~\cite{Raaz_Ho_Koulik_2018} showed that
the EM algorithm takes
$\mathcal{O}((n/d)^{\myfrac{1}{2}})$ steps to converge to a Euclidean
ball of radius $\mathcal{O}((d/n)^{\myfrac{1}{4}})$ around the true
location parameter.  Recall that for such models, the MLE is known to
lie at a distance $\mathcal{O}(n^{-\myfrac{1}{4}})$ from the true
parameter~\cite{Chen1992}, so that even though its convergence rate as
an optimization algorithm is slow; the EM algorithm nonetheless produces a solution with a statistical error of the same order as the
MLE.  This past work does not consider the more realistic setting in
which both the location and scale parameters are unknown, and the EM algorithm is used to fit both simultaneously.  Indeed, as mentioned
earlier, such models may become weakly identifiable due to algebraic
relations among the partial
derivatives~\cite{Chen_2009}.  Thus, analyzing EM in the
case of weakly identifiable mixtures is challenging for two reasons:
the weak separation between the mixture components, and the algebraic interdependence of the partial derivatives of the log-likelihood.
The main contributions of this work are (a) to highlight the
dramatic differences in the convergence behavior of the EM algorithm, depending
on the structure of the fitted model relative to the data-generating
distribution; and (b) to analyze the EM algorithm under a few specific
yet representative settings of weakly identifiable models, giving a
precise analytical characterization of its convergence behavior.
%\vspace{-2mm}

\subsection{Some illustrative examples} % (fold)
\begin{figure*}[t]
  \begin{center}
    \begin{tabular}{ccc}
      \widgraph{0.30\textwidth}{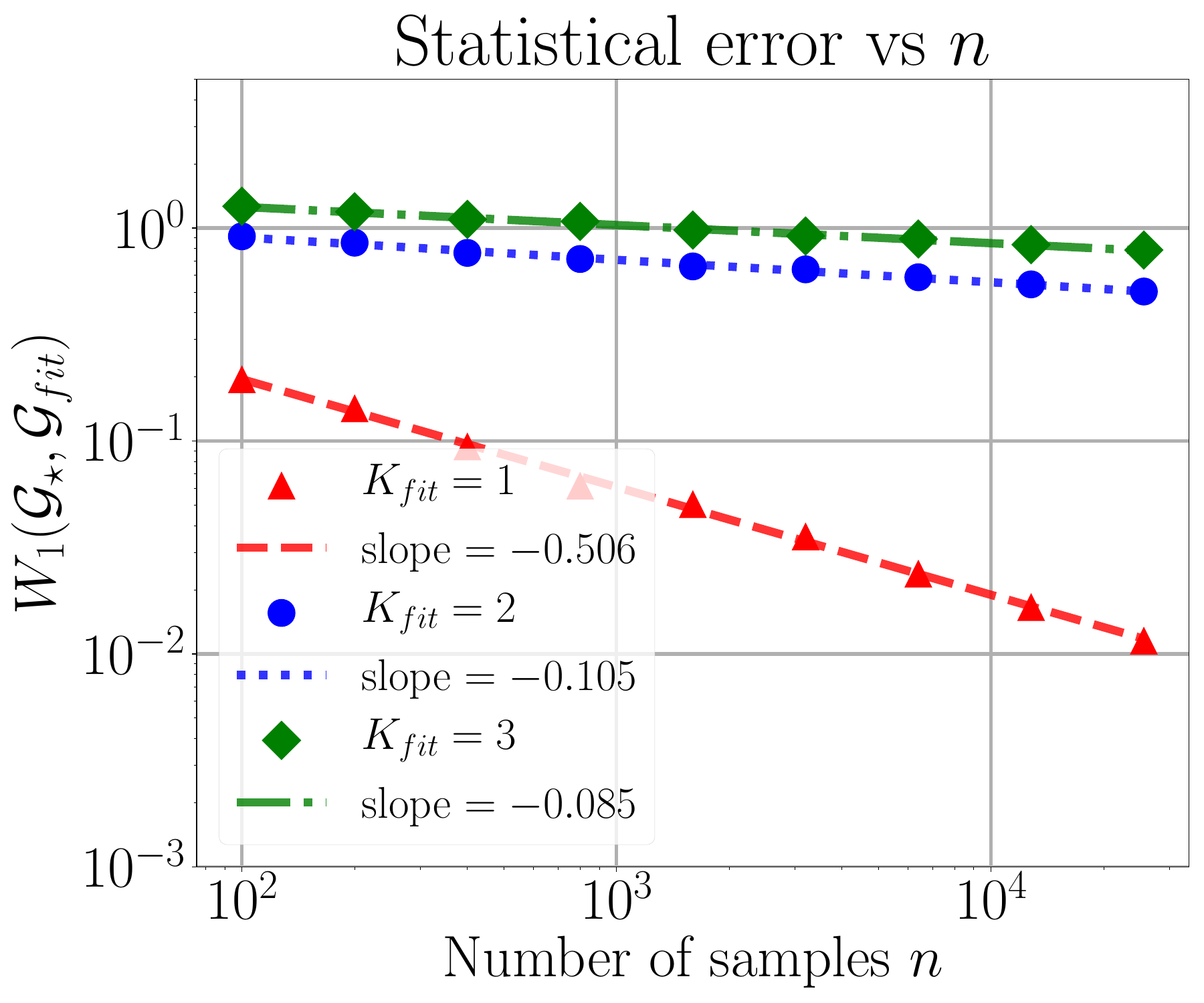} &
      \widgraph{0.30\textwidth}{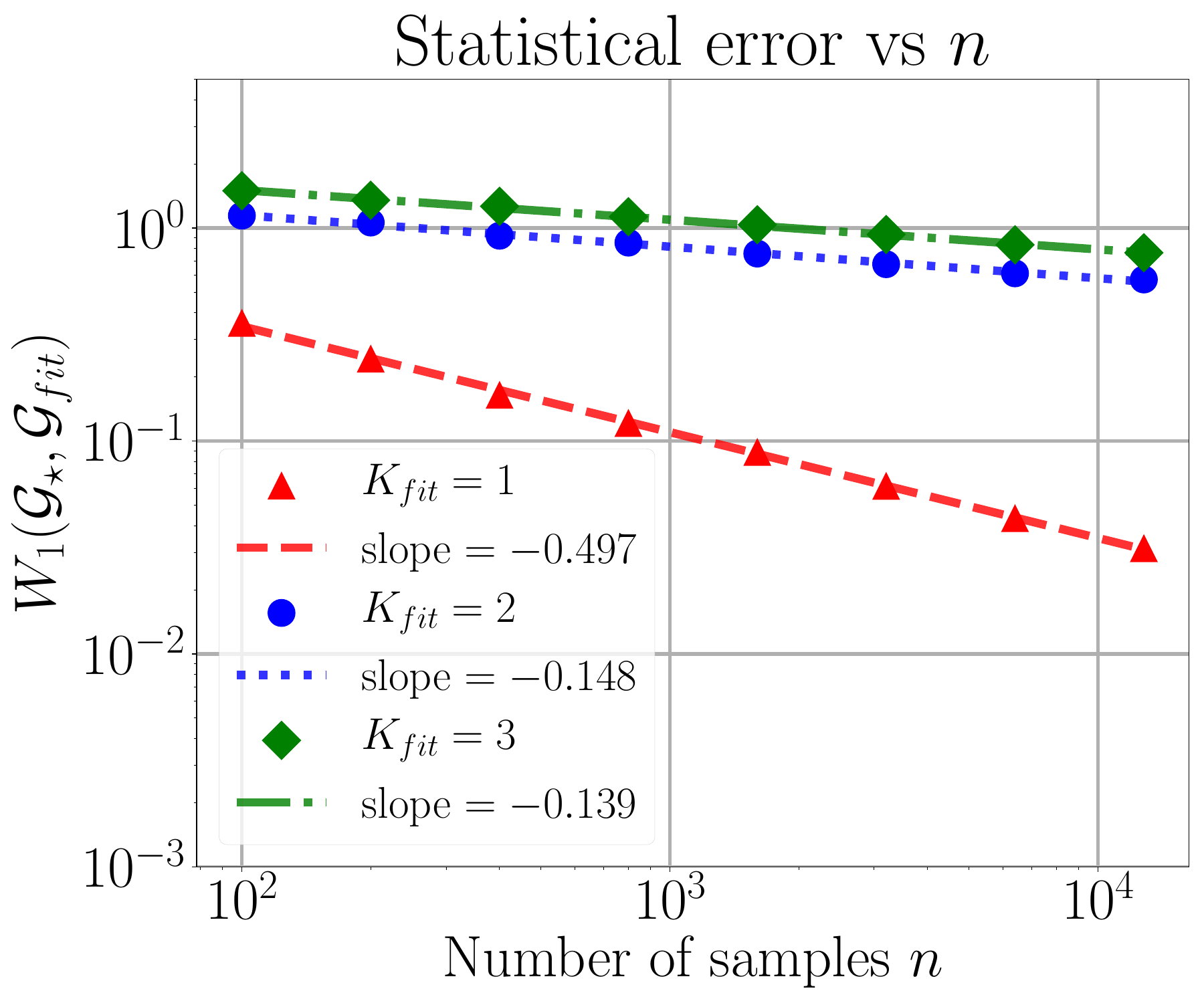} & 
      \widgraph{0.30\textwidth}{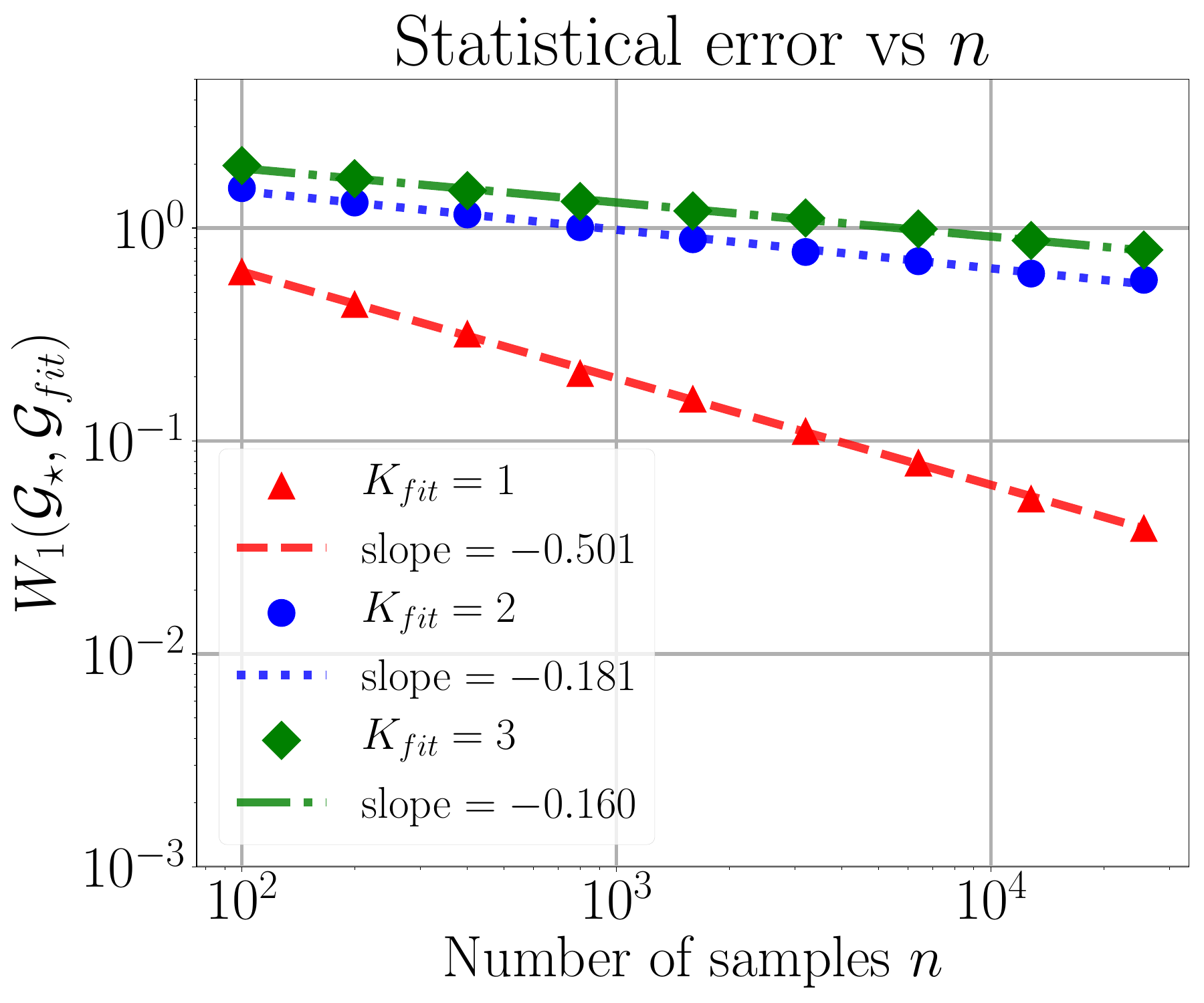}\\
      (a) $d=1$ & (b) $d=2$ & (c) $d=4$
    \end{tabular}
    \caption{Scaling of the Wasserstein error between the
    true parameters and the EM estimates, when EM is
    used to fit a Gaussian mixture model with $K_{\textrm{fit}} \in \braces{1,
    2, 3}$ components, i.e., $\mathcal{G}_{\textrm{fit}} =\sum_{i=1}^{K_{\textrm{fit}}}
    w_i \Ncal(\mu_i, \Sigma_i)$, on an $n$ sample-dataset generated from
    standard Gaussian distribution $\mathcal{G}_* = \Ncal(0, I_ {\dims})$.
    In all three examples, when the fitted model is over-specified, meaning that the fitted model has more components
    than the true model ($K_{\textrm{fit}} \in  \{2, 3\}$ in these examples),
    we
    observe a significant increase in the
    Wasserstein error. Stated differently, the simulations 
    suggest that the estimation accuracy of the EM algorithm 
    degrades dramatically when the fitted model is over-specified.}
    % %\vspace{-7mm}
    \label{fig:gaussian_vs_mixture}
  \end{center}
\end{figure*}

\label{sub:some_illustrative_examples}
%\vspace{-2mm}
Before proceeding further, we summarize a few common aspects of the numerical experiments and the associated
figures presented in the paper. Computations at the population-level were
done via numerical integration on a sufficiently fine grid.  For EM with
finite sample size $n$, we track its performance for several values of $n
\in \braces{100, 200, 400, \ldots, }$ and report the quantity $\widehat
m_e + 2\widehat s_e$ on the y-axis, where $\widehat m_e \text{ and } \widehat s_e$, 
respectively,
denote the mean and standard deviation across the experiments for the metric
under consideration (as a function of $n$ on the x-axis, e.g., Wasserstein
error for parameter estimation in Figure~\ref{fig:gaussian_vs_mixture}.
The stopping criteria for sample EM were: (a)
the change in the iterates was small enough $(<.001/n)$, or (b) the number
of
iterations was too large (greater than $100,000$); criteria~(a) led to convergence
in most experiments.
Furthermore, whenever we provide a slope, it is the slope for
the least-squares fit on the log-log scale for the quantity on
$y$-axis when fitted with the quantity reported on the $x$-axis.  
For
instance, in Figure~\ref{fig:gaussian_vs_mixture}(a), we plot the
Wasserstein error between the estimated mixture and the true mixture
on the $y$-axis value versus
the sample size $n$ on the $x$-axis and also provide the slopes for the least-squares fit.
In particular, in panel (a) the green dot-dashed line with the legend
`slope$=-0.09$' denotes the least-squares fit and the respective slope for
the logarithmic error $\log W_1(\mathcal{G}_*, G_\text{fit})$ (green diamonds) with
respect to the logarithmic sample size $\log n$ when the number of components
in the fitted mixture is $3$. Such a result implies that the error
$ W_1 (\mathcal{G}_*, G_\text{fit})$ scales as $n^{-0.09}$ with the sample size $n$ in
our experiments.

To begin with, we consider the simplest case of over-specification with
Gaussian mixture
models---when the true data is generated from a zero-mean standard Gaussian distribution in $d$ dimensions and EM is used to fit a general  multi-component mixture model with different number of mixtures.
% \footnote{
(We note that fitting by one mixture model is simply a Gaussian
fit.)
% }. 
Given the estimates for the mixture weights, location and scale parameters
returned by EM, we compute the first order Wasserstein distance\footnote{
First-order Wasserstein distance has been used in prior works to characterize
the error between the estimated and true parameters. See section 1.1~\cite{Ho-Nguyen-Ann-16}.}
between the true
and estimated parameters.
Results for $d\in\braces{1, 2, 4}$ and for various amount of over-specification
are plotted in Figure~\ref{fig:gaussian_vs_mixture}.
From these results, we notice that the decay in statistical error is $\obs^
{-1/2}$ when the fitted number of components is well-specified and equal
to the true number of components but has a much slower rate whenever the number of fitted components is two or more.
Moreover, in Section~\ref{sec:discussion} (see Figure~\ref{fig:mixture_vs_mixture})
we show that such a phenomenon occurs more generally in mixture
models.

 While a rigorous
theoretical analysis of EM under over-specification in general mixture models
is desirable, it remains beyond the scope of this paper. Instead, here we
provide a full characterization of EM when it is used  to fit the following
class of models to the data drawn from standard Gaussian $\NORMAL(0, I_\dims)$:
\begin{align}
        \mixfit{symm}((\theta, \sd^2)) 
    \!=\!\frac{1}{2} \NORMAL(\theta, \sigma^2 I_d) + \frac{1}{2}\NORMAL (-\theta,
      \sigma^2I_d).\! \label{eq:symmetric_fit}
\end{align}
In particular, in this symmetric fit, we fix the mixture
weights to be equal to $\frac{1}{2}$ and require that the two components
have same scale parameter. Given the estimates $\widehat{\theta}, \widehat{\sigma}$,
the Wasserstein error (see equation~\eqref{eq:wass_dist} in Appendix~\ref{sec:wass_dist})
in this case can be simplified
as
$\enorm{\widehat{\theta}} + \sqrt{d}\sqrt{\abss{\widehat{\sigma}^2-1}}$.
In our results to be stated later, we show that the two
terms are of the same order (equations~\eqref{eqn:sigma_rate}, \eqref{eq:sigma_multivariate}) 
and hence we primarily focus
on the error $\enorm{\widehat{\theta}- \theta_\star}$ going forward to simplify
the exposition.
We consider our set-up as a simple yet first step towards understanding
the behavior of EM in over-specified mixtures when \emph{both} location and
scale parameter are unknown.
% \footnote{
In our prior work~\cite{Raaz_Ho_Koulik_2018},
we studied the slow down of EM with over-specified mixtures for estimating
only the location parameter, but they assumed that the scale parameter was
known and fixed. Here a more general setting is considered.
% }

We now elaborate the choice of our class of models~\eqref{eq:symmetric_fit}
that may appear a bit restrictive at first glance. This model turns
out to be the simplest example of a weakly identifiable model in $d=1$.
Let $\normDensity$ denote the density of a Gaussian distribution with mean
$\theta$ and variance $\sd^2$, then we have%\vspace{-2mm}
\begin{align} 
\label{EqnAlgebra}
  \dfrac{ \partial^{2}{ \normDensity}}{ \partial{ \theta^{2}}} (x;
  \mean, \sd^2) = 2 \dfrac{ \partial { \normDensity}} { \partial{
      \sigma^{2}}} (x; \mean, \sd^2), 
\end{align}
valid for all $x \in \real$, $\mean \in \real$ and $\sd > 0$.
As alluded to earlier, models with algebraic dependence between partial
derivatives lead to weak identifiability and slow statistical estimation
with MLE. However, in the multivariate setting when the same 
parameter $\sigma$ is shared across multiple dimensions, this
algebraic relation does not hold and the model is strongly identifiable
(since the Fisher information matrix is singular at $(\theta^*,\sd^* ) : =
(0, 1)$). For this reason, we believe that analysis of EM for the special
fit~\eqref{eq:symmetric_fit} may provide important insight for more general 
over-specified weakly identifiable models.

\textbf{Population EM:} Given $n$ samples from a $d$-dimensional standard
Gaussian distribution, the
sample EM algorithm for location and scale parameters generates a sequence of 
the form $\theta^{t+1} = M_{n,d}(\theta^t)$ and $\sigma^{t + 1}$, which is some function of $\enorm{ \theta^{t + 1}}^2$; see
equation~\eqref{eq:sample_em_operators} for a precise definition.  An
abstract counterpart of the sample EM algorithm---not useful in
practice but rather for theoretical understanding---is the population
EM algorithm $\gencorop$, obtained in the limit of an infinite sample
size (cf. equation~\eqref{eq:stage_2_pop_em_operator}).

\begin{figure*}[t]
  \begin{center}
      \begin{tabular}{cc}
        \widgraph{0.35\textwidth}{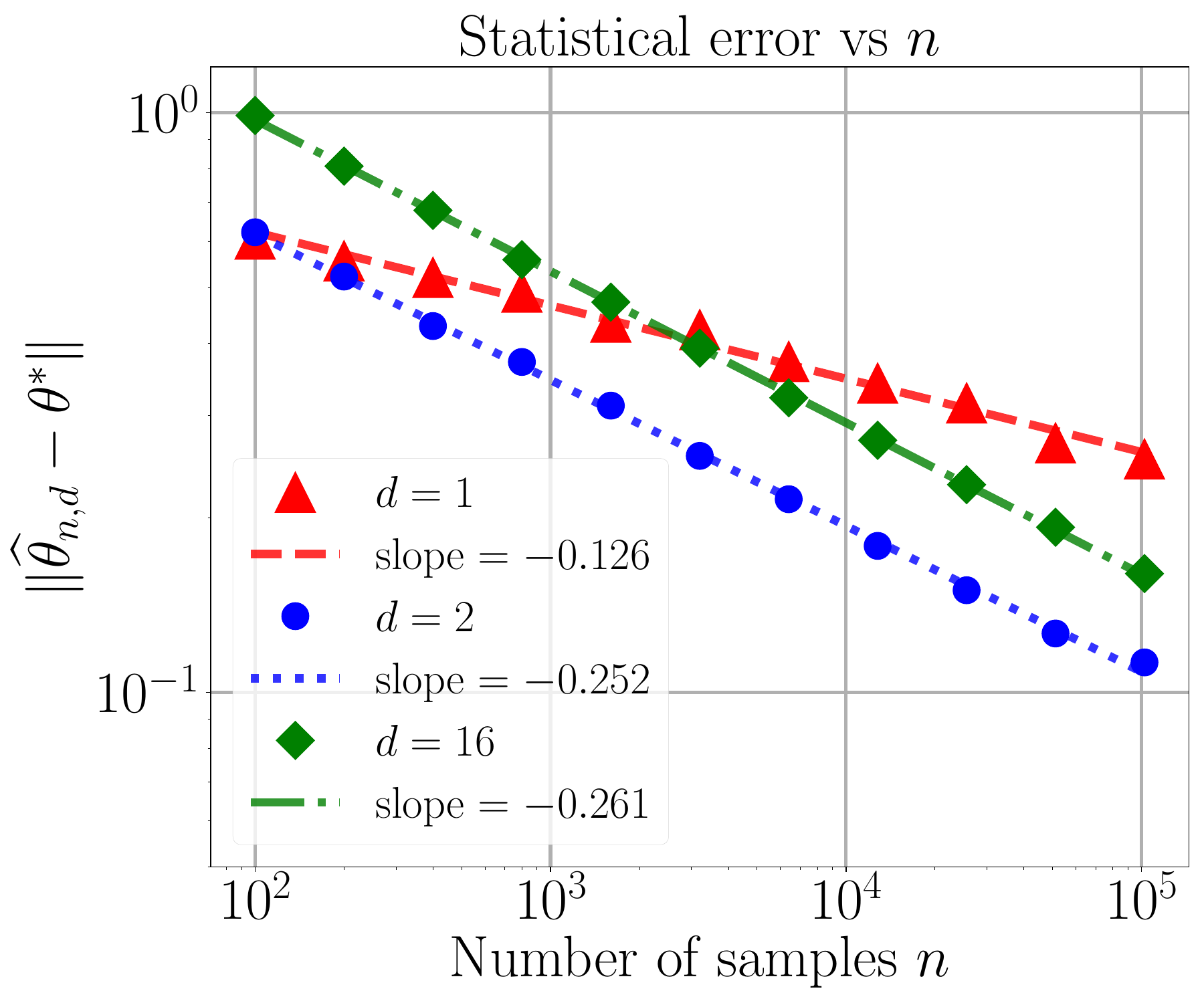}
        & 
        \widgraph{0.35\textwidth}{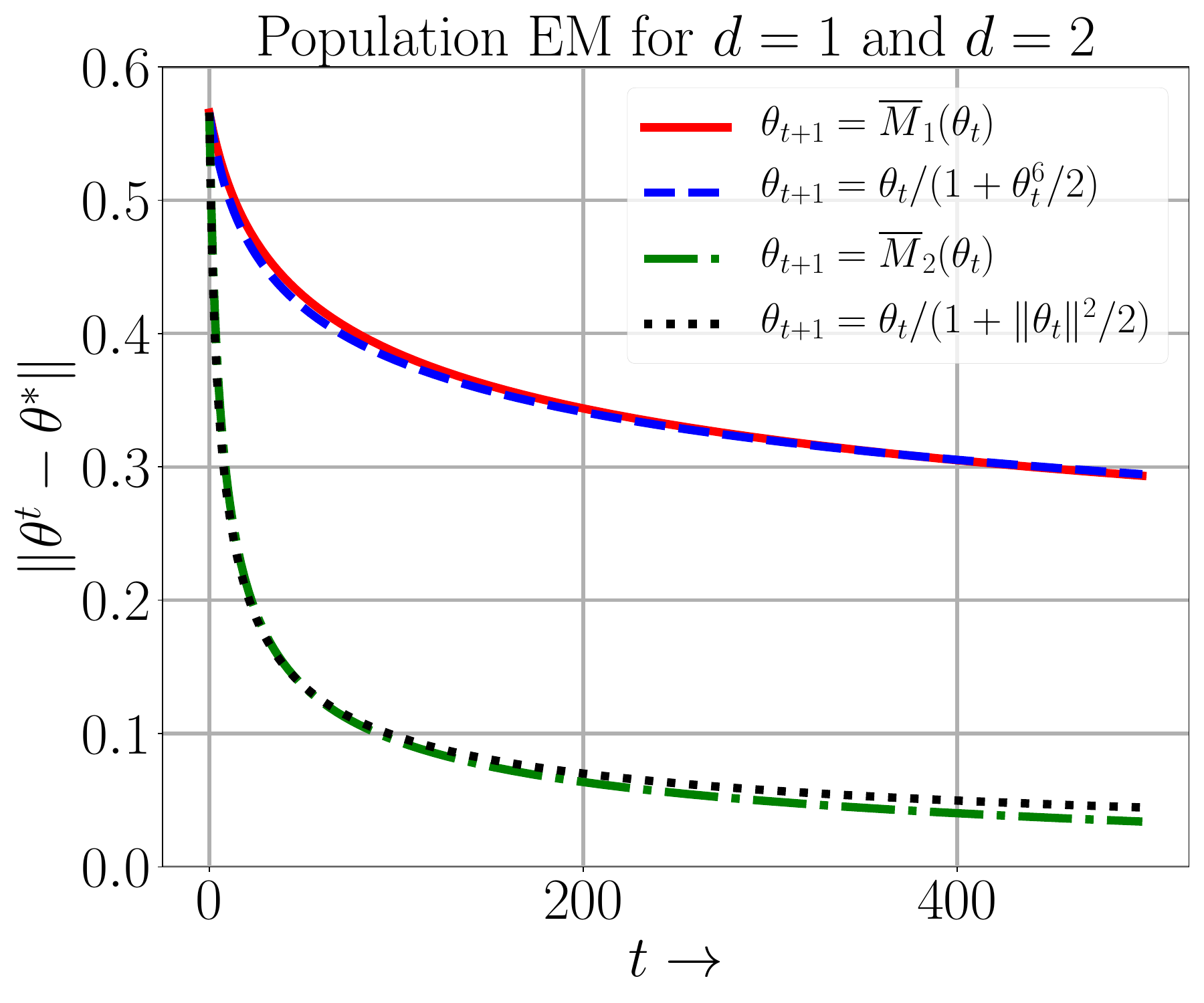}\\
        (a) & (b)
      \end{tabular}
      \caption{
       Behavior of the EM algorithm for the fitted model~\eqref{eq:symmetric_fit},
      where the data is being generated from $\Ncal(0, I_ {\dims})$. 
       (a) Scaling of the Euclidean error 
        $\|\widehat{\theta}_{\obs, \dims} - \thetastar\|_2$ 
        with respect to the sample size $\obs$ for dimension
        $\dims \in \{1, 2, 16\}$. Here, $\widehat{\theta}_{\obs, \dims}$
        denotes the EM algorithm estimate of the mean parameter $\theta$
        based on $\obs$ samples.
        Note that the simulations
        indicate two distinct error scaling for $d = 1$ and $d > 1$. 
      (b) Convergence behavior of the population-like EM
        sequence $\theta^{t+1} =
        \gencorop(\theta^t)$~\eqref{eq:stage_2_pop_em_operator}
        in dimensions $d = 1$ and $2$. The rate of convergence
        in dimension $d = 1$ is significantly slower compared to 
        the rate in dimension $d = 2$. Overall, both the plots provide strong
        empirical evidence towards two distinct behaviors of the 
        EM algorithm for dimension $d = 1$ and dimensions $d > 1$.
        See the Theorems~\ref{thm:em_univariate}-\ref{thm:em_multivariate}, 
        and Lemmas~\ref{lemma:correct_pop_contraction_and_radem_three} 
        and \ref{lemma:d_dim_operator} for a theoretical justification
        of trends in panels (a) and (b) respectively.
        }
        % %\vspace{-5mm}
      \label{fig:symm_fit_all_results}
    \end{center}
\end{figure*}

In practice, running the sample EM
algorithm yields an estimate $\widehat{\theta}_{\obs, \dims}$ of the
unknown location parameter $\thetastar$.  Panel (a) in Figure~\ref{fig:symm_fit_all_results} shows the scaling
of the statistical estimation error $\|\widehat{\theta}_{\obs, \dims}
- \thetastar\|_2$ of this sample EM estimate versus the sample size
$n$ on a log-log scale.  The three curves correspond to dimensions $d
\in \{1, 2, 16 \}$, along with least-squares fits (on the log-log scale)
to the data.
 In panel (b), we plot the Euclidean norm
$\|\theta^t\|_2$ of the population EM iterate\footnote{In fact, our
  analysis makes use of two slightly different population-level
  operators $\PseudoND$ and $\gencorop$ defined in
  equations~\eqref{eq:pop_like_em_operator} and
  \eqref{eq:stage_2_pop_em_operator} respectively.
  Figure~\ref{fig:symm_fit_all_results}(b) shows plots for the operator
  $\gencorop$, but the results are qualitatively similar for the operator
  $\PseudoND$.} versus the iteration number $t$, with solid red line
corresponding to $\dims = 1$ and the dash-dotted green line
corresponding to $\dims = 2$.  Observe that the algorithm converges
far more slowly in the univariate case than the multivariate case.
The theory to follow in this paper (see
Theorems~\ref{thm:em_univariate}, \ref{thm:em_multivariate} and Lemmas~
\ref{lemma:correct_pop_contraction_and_radem_three}
and \ref{lemma:d_dim_operator})
provides explicit predictions for the rate at which different quantities
plotted in Figure~\ref{fig:symm_fit_all_results} should decay.
We now summarize our theoretical results that are also consistent with
the trends observed in Figure~\ref{fig:symm_fit_all_results}.
%\vspace{-2mm}

\subsection{Our contributions} % (fold)
\label{sub:our_contributions}
%\vspace{-2mm}
The main contribution of this paper is to provide a precise analytical
characterization of the behavior of the EM algorithm for certain
special cases of over-specified mixture models~\eqref{eq:symmetric_fit}.

\textbf{Univariate over-specified Gaussian mixtures:} In the
  univariate setting ($d = 1$) of $\mixfit{symm}$ in~\eqref{eq:symmetric_fit},
  we prove that the EM estimate has statistical estimation error of
  the order $n^{-\myfrac{1}{8}}$ and $n^{-\myfrac{1}{4}}$ after order
  ${n^{ \frac{3}{4}}}$ steps for the location and scale parameters
  respectively. In particular, Theorem~\ref{thm:em_univariate} provides
  a theoretical justification for the slow rate observed in Figure~\ref{fig:symm_fit_all_results}
  (a) for $d=1$ (red dotted line with star marks).
  Proving these rates requires
  a novel analysis, and
  herein lies the main technical contribution of our paper. Indeed, we
  show that all the analysis techniques introduced in past work on EM,
  including work on both the regular~\cite{Siva_2017} and strongly
  identifiable cases~\cite{Raaz_Ho_Koulik_2018}, lead to sub-optimal
  rates.  Our novel method is a \emph{two-stage approach} that makes use
  of two different population level EM operators.
  Moreover, we also prove a matching lower bound (see Appendix~\ref{sec:minimax_bound})
  which ensures that the upper bound of order $\obs^{-\myfrac{1}{8}}$ for the
  statistical error of sample EM from Theorem~\ref{thm:em_univariate} is
  tight up to constant factors.

\textbf{Multivariate setting with shared covariance:} Given the technical
  challenges even in the simple univariate case, the symmetric spherical
  fit $\mixfit{symm}$ in~\eqref{eq:symmetric_fit} serves as a special case
  for the multivariate setting $d\geq 2$.
  In this case, we establish that the sharing of scale parameter proves
  beneficial in the convergence of EM. Theorem~\ref{thm:em_multivariate}
  shows that sample EM algorithm takes $\mathcal{O}((n/ d)^{1/2})$
  steps in order to converge to estimates, of the location and scale parameters
  respectively, that lie within distances $\mathcal{O}(d/ n)^{ 1/ 4}$
  and $\mathcal{O}(n d)^{-\myfrac{1}{2}}$ of the true location and
  scale parameters, respectively.

  %   This result is a consequence of a
  % localization argument introduced in our past work~\cite{Raaz_Ho_Koulik_2018}
  % and the faster convergence of population EM algorithm \emph{for
  %   the symmetric fit} (Lemma~\ref{lemma:d_dim_operator} in the sequel)
  %   which requires $\order{1/ \epsilon^2}$ and $\order{1/ \epsilon}$ steps
  %   in order to achieve $\epsilon$-accuracy for the location and scale parameters,
  %   respectively.  %\vspace{-3mm}

\textbf{General multivariate setting:} % (fold)
 We want to remind the readers that we expect the Wasserstein error to scale
much slowly than $n^{-\frac{1}{4}}$ (the rate mentioned in the previous
  paragraph) while estimating over-specified mixtures with no shared covariance.
  When the fitted variance parameters are not shared across dimensions
  our simulations under general multi-component fits
  in Figure~\ref{fig:gaussian_vs_mixture} demonstrate a much slower convergence
  of EM (for which a rigorous justification is beyond the scope of this
  paper).

%%%%%%%%%%%%%%%%%%%%%%%%%%%%%%%%%%%%%%%%%%%%%%%%%%%%%%%%%%%%%%%%%%%%%%%%

\textbf{Notation:} In the paper, the expressions $a_{n} 
\precsim b_{n}$ or $a_{n} \leq O(b_{n})$ will be used to denote 
$a_{n} \leq c b_{n}$ for some positive universal constant $c$ that does
not change with $n$. Additionally, we write $a_{n} \asymp b_{n}$ if both
$a_{n} \precsim b_{n}$ and $b_{n} \precsim a_{n}$ hold. Furthermore, we
denote $[n]$ as the set $\{1, \ldots, n\}$ for any $n \geq 1$. We define $\ceil{x}$ as the smallest integer greater than or equal to $x$
for any $x \in \Rspace$. The notation $\|x\|_{2}$ stands for the $\ell_
{2}$ norm of vector $x \in \Rspace^{d}$. We use $c, c', c_1$ etc. to denote
some universal constants independent of problem parameters (which might
change in value each time they appear).

% Before proceeding further we briefly state the EM updates that we use
% later in our results.

\subsection{EM updates for symmetric fit $\mixfit{symm}$} % (fold)
\label{sub:em_updates_for_symmetric_fit_symm_}

The EM updates for Gaussian mixture models are standard, so we simply
state them here. In terms of the shorthand notation $\eta \defn
(\theta, \sd)$, the E-step in the EM algorithm involves computing the
function
\begin{subequations}
  \begin{align*}
& \samsurro( \eta'; \eta) \mydefn \frac{1}{n} \sum_{i = 1}^{n}
\big[\weightFun_{ \theta, \sd}( X_{i}) \log \parenth{
    \normDensity(X_{i}; \mean', (\sd')^2 I_{d})} \nonumber \\
    & + \parenth{ 1 -
    \weightFun_{ \theta, \sd}( X_{i}) } \log \parenth{
    \normDensity(X_{i}; - \mean', (\sd')^2 I_{d})}\big],
\end{align*}
where the weight function is given by 
$\weightFun_{ \theta, \sigma}( x)~=~({1+ e^{\frac{-2\theta\tp x}{\sd^2}}})^{-1}$.
% \begin{align*}
% \weightFun_{ \theta, \sigma}( x) \mydefn \frac{ \exp \parenth{ \theta^{\top} x/ \sigma^2 }}{ \exp \parenth{\theta^{\top} x/ \sigma^2} + \exp \parenth{- \theta^{\top} x/ \sigma^2}}.
%  \end{align*}
The M-step involves maximizing the $\samsurro$-function
over the pair $(\theta', \sd')$ with $\eta$ fixed, which yields
\begin{align}
    % \label{eq:sample_em_theta}
    \mean' = \frac{1}{ \obs} \sum_{i =
    1}^\obs
    (2 \weightFun_{ \mean, \sd} (X_i) - 1) X_i, \quad \text{ and}
    \\
    \quad
    (\sd')^2 = \frac{ 1}{ d} \parenth{ \frac{ \sum_{i = 1}^\obs \enorm{
         X_i}^2}{ \obs} - \enorm{ \mean'}^2},
      \label{eq:sample_em_sigma}    
\end{align}
Doing some straightforward algebra, the 
EM updates $(\theta^{t}_n, \sd_n^{t})$ can be succinctly defined as
\begin{align}
\label{eq:sample_em_operators}
  \theta_{n}^{t+1} & =
  \frac{1}{n}\sum_{i=1}^n \tanh\parenth{\frac{X_i^\top\theta_{n}^t}{\sum_
  {i=1}^n \enorm{X_i}^2/(nd) - \enorm{\theta_n^t}^2/d}} \nonumber \\
  & \rdefn M_{n, d}( \theta_n^{t}),
\end{align}
and $\sd_n^{t+1} = \sum_{i=1}^n \enorm{X_i}^2/(nd) - \enorm{\theta_n^{t + 1}}^2/d$.
For simplicity in presentation, we refer to the operator $M_{n, d}$ as the 
\emph{sample EM operator}. 
\end{subequations}

\textbf{Organization:}
We present our main results in Section~\ref{sec:our_results}, with
Section~\ref{sub:em_for_univariate_singular_mixtures} devoted
to the univariate case, Section~\ref{sub:results_for_multivariate_case}
to the multivariate case and Section~\ref{sub:general_cases} to the simulations
with more general mixtures.
Our proof ideas are summarized in Section~\ref{sec:EM_analysis} and we conclude with a discussion in Section~\ref{sec:discussion}.
The detailed proofs of all our results are deferred to the Appendices.
%\vspace{-2mm}

%%%%%%%%%%%%%%%%%%%%%%%%%%%%%%%%%%%%%%%%%%%%%%%%%%%%%%%%%%%%%%%%%%%%%%%%%%%%%

\section{Main results} % (fold)
\label{sec:our_results}
%\vspace{-2mm}
In this section, we provide our main results for the behavior of EM
with the singular (symmetric) mixtures fit $\mixfit{symm}$~\eqref{eq:symmetric_fit}.
Theorem~\ref{thm:em_univariate} discusses the result for the 
univariate case, Theorem~\ref{thm:em_multivariate} discusses the result
for multivariate case. In Section~\ref{sub:general_cases} we discuss some simulated experiments 
for general multivariate location-scale Gaussian mixtures.   
%%%%%%%%%%%%%%%%%%%%%%%%%%%%%%%%%%%%%%%%%%%%%%%%%%%%%%%%%%%%%%%%%%%%%%%%%%%%%
%\vspace{-2mm}

%%%%%%%%%%%%%%%%%%%%%%%%%%%%%%%%%%%%%%%%%%%%%%%%%%%%%%%%%%%%%%%%%%%%%%%%%%%%%

\subsection{Results for the univariate case} 
\label{sub:em_for_univariate_singular_mixtures}
%\vspace{-2mm}
As discussed before, due to the relationship between the location and
scale parameter, namely the updates~\eqref{eq:sample_em_operators}, 
it suffices to analyze the sample EM operator for the location parameter. For the univariate Gaussian mixtures, given $n$ samples
$\braces{X_i, i\in [n]}$, the sample EM operator
is given by
\begin{align}
\label{eq:sample_em_simple_operator}
M_{\obs, 1}( \theta) & \mydefn \frac{ 1}{ n} \sum_{i = 1}^\obs X_i \tanh
\brackets{ \frac{ X_i \theta}{\sum_{j = 1}^\obs X_j^2/ \obs -
    \theta^2}}.
\end{align}  
We now state our first main result that characterizes the guarantees for
EM under the univariate setting.
Let $\interval_{\smallthreshold}'$ denote the interval $[ \unicon \obs^
{- \frac{1}{12} + \smallthreshold}, 1/ 10]$ where $c$ is a positive universal
constant.
\begin{theorem} 
\label{thm:em_univariate}
Fix $\delta \in (0, 1)$, $\smallthreshold \in (0, 1/ 8]$, and let 
$X_i\stackrel{i.i.d.}{\sim}\mathcal{N}(0, 1)$ for $i=1, \ldots, \obs$ such
that \mbox{$\samples \succsim \log \frac{ \log(1/ \smallthreshold)}
{ \delta}$}. Then for any initialization $\mean_{n}^{0}$ that satisfies
$\vert{ \mean_{n}^{0}}\vert \in \interval_{\smallthreshold}'$, 
the sample EM sequence $ \mean_{n}^{t+ 1} = \samm( \mean_{n}^t)$,
satisfies
\begin{align}
\label{eq:sample_em_univariate_rate}
\vert{ \mean_{n}^t - \thetastar} \vert 
\leq 
c_1\frac{1}{n^{1/ 8 -\smallthreshold}} \log^{5/4}\parenth{\frac{ 10 n \log
(8/ \smallthreshold)}{ \delta}},
  % \quad\text{for all }
  % t \succsim \obs^{\frac{3}{ 4} - 6 \smallthreshold}\cdot
% \log \obs \log\frac{1}{\smallthreshold},
  \end{align}
  for all $t \geq c_2 \obs^{\frac{3}{ 4} - 6 \smallthreshold}\cdot
\log \obs \log\frac{1}{\smallthreshold}$
with probability at least $1 -\delta$. 
\end{theorem}
\noindent See
Appendix~\ref{sub:proof_of_theorem_thm:em_univariate} for
the proof. 

The bound~\eqref{eq:sample_em_univariate_rate} shows that with high
probability after $\mathcal{O}(n^{ 3/ 4})$ steps the sample EM
iterates converge to a ball around $\thetastar$ whose radius is
arbitrarily close to $\obs^{-1/8}$. Moreover, as a direct consequence of
the relation~\eqref{eq:sample_em_sigma}, we conclude
that the EM estimate for the scale parameter is of order $\obs^{-\frac{1}
{4}}$ with high probability:
\begin{align}
 \label{eqn:sigma_rate}    
\abss{ ( \sd_n^{t})^2 - ( \sd^{*})^2} &= \abss{ \frac{ \sum_{i =
      1}^\obs { X_i}^2}{ \obs}  - {\parenth{
          \theta^{t}_n - \thetastar}^2} - ( \sd^{*})^2} \notag\\
    &\precsim \obs^{-\frac{1}{2}} + \obs^{-\frac{1}{4}} 
    = O(\obs^{-\frac{1}{4}})
\end{align}
where we have used the standard chi-squared concentration
for the sum ${ \sum_{i =1}^\obs { X_i}^2}/{ \obs} $.

\begin{figure*}[t]
  \begin{center}
    \begin{tabular}{cc}
      \widgraph{0.4\textwidth}{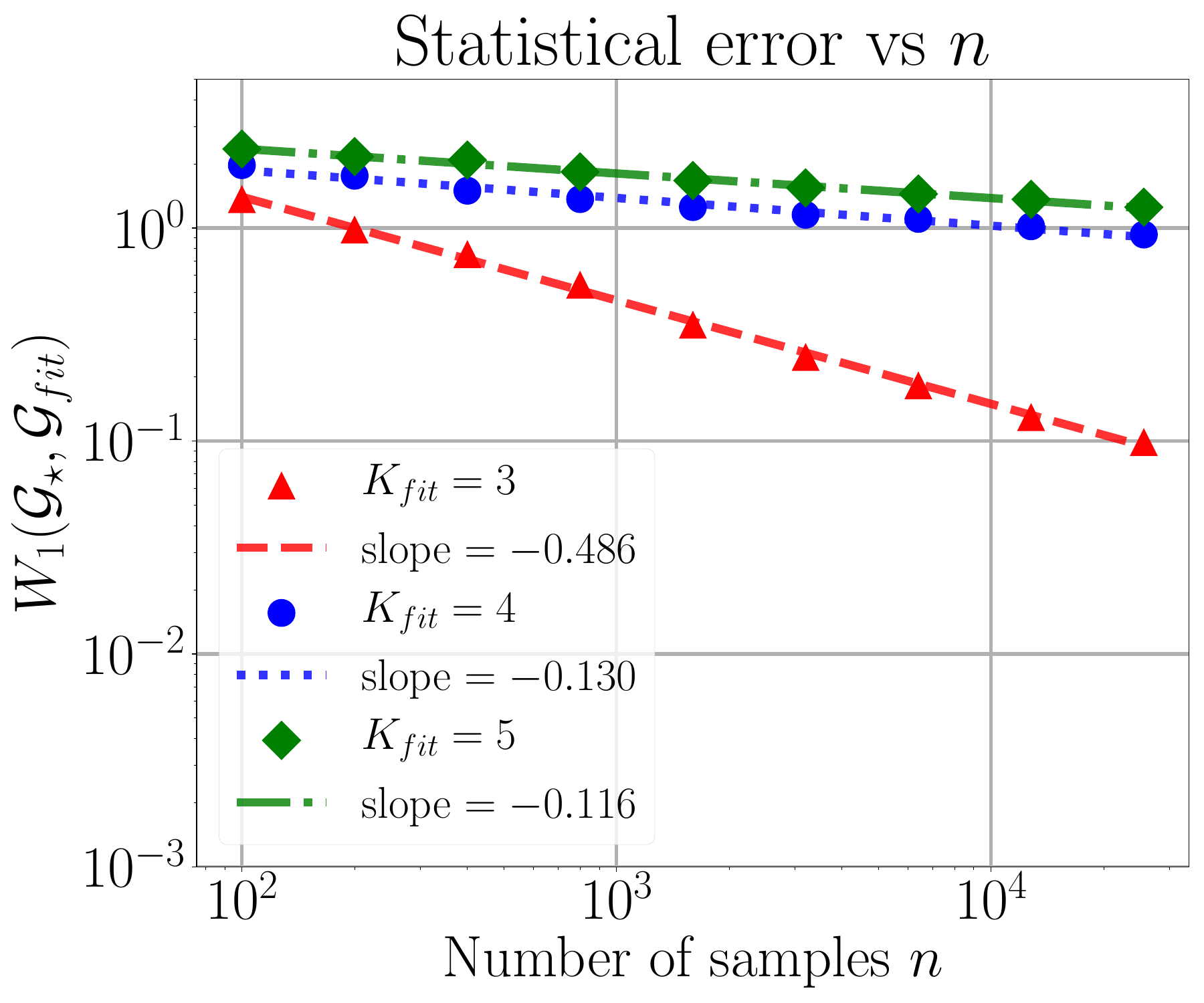} &
      \widgraph{0.4\textwidth}{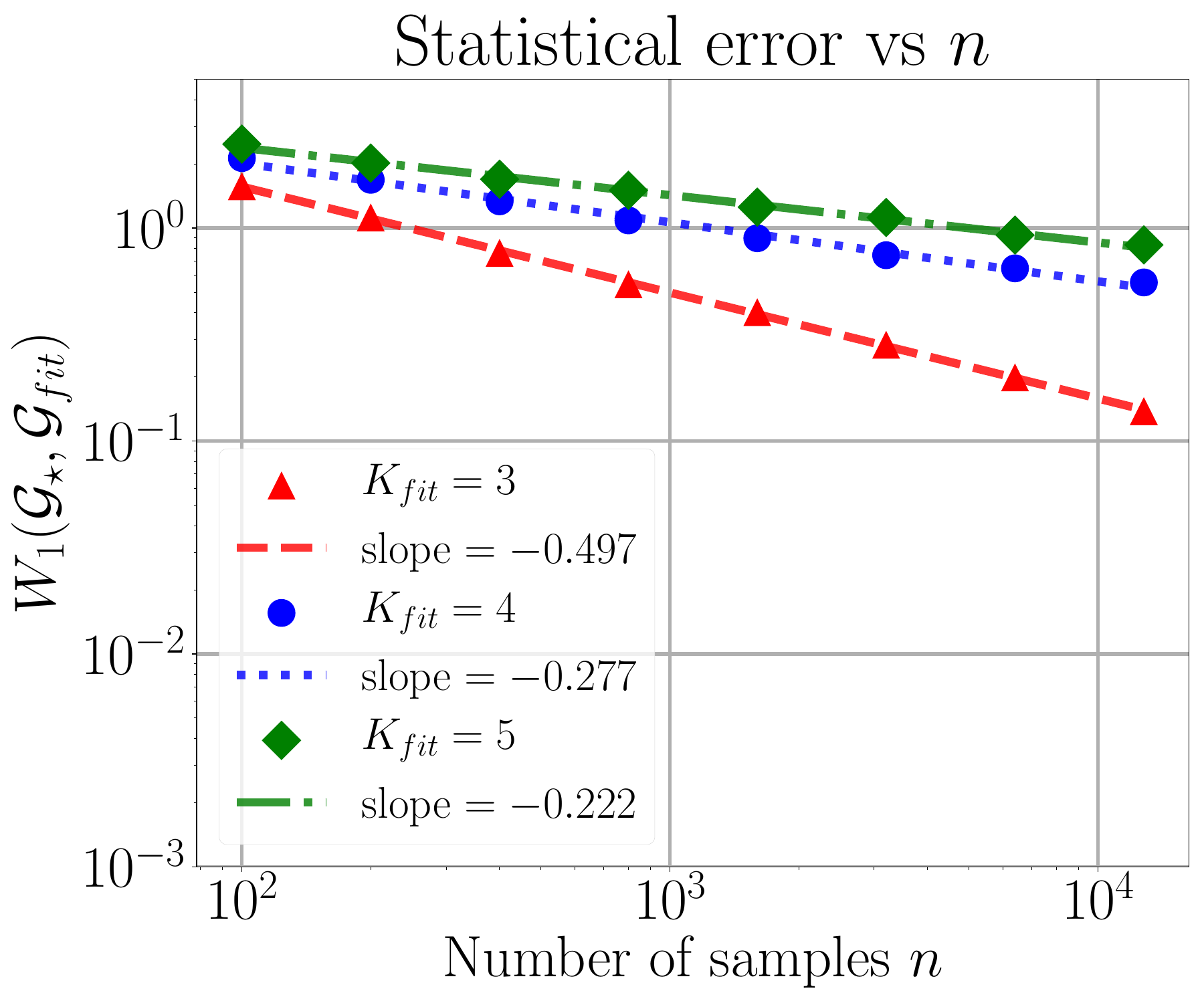} \\
      (a) $d=1$ with $\mathcal{G}_*$ given by equation~\eqref{eq:three_mixture_one_d}
      & (a) $d=2$ with $\mathcal{G}_*$ given by equation~\eqref{eq:three_mixture_two_d}
    \end{tabular}
    \caption{Scaling of the 
    first-order Wasserstein error for EM estimates when fitting a Gaussian
    mixture with $K_{\mathrm{fit}} \in \braces{3, 4, 5}$,
    i.e., $\mathcal{G}_{\textrm{fit}} =\sum_{i=1}^{K_{\textrm{fit}}}
    w_i \Ncal(\mu_i, \Sigma_i)$, on $n$ i.i.d. samples 
    from a $3$-Gaussian mixture model (equations \eqref{eq:three_mixture_one_d}
    and \eqref{eq:three_mixture_two_d}).  
    In the case of no over-specification, i.e.,  $K_{\textrm{fit}}
    = K_{\textrm{true}} = 3$, the error scales as
    $n^{-1/2}$, but when the fitted
    model is over-specified ($K_{\textrm{fit}} \in \braces{4, 5}$), the
    scaling is much worse (and degrades further for any given $n$ as 
    $K_{\textrm{fit}}$ gets large). See Section~\ref{sub:general_cases} for further details.}
    % \vspace{-5mm}
    \label{fig:mixture_vs_mixture}
  \end{center}
\end{figure*}

\textbf{Matching lower bound:} 
In Appendix~\ref{sec:minimax_bound}, we prove a matching lower
bound and thereby
conclude that the upper bound of order $\obs^{-\myfrac{1}{8}}$ for the
statistical error of sample EM from Theorem~\ref{thm:em_univariate} is
tight up to constant factors. In Section~\ref{sub:general_cases}, we provide
further evidence (cf. Figure~\ref{fig:mixture_vs_mixture}) that the slow 
statistical rates of EM with location parameter that we derived in Theorem~\ref{thm:em_univariate} 
might appear in more general settings of location-scale Gaussian mixtures
as well.
%\vspace{-2mm}

\subsection{Results for the multivariate case} % (fold)
\label{sub:results_for_multivariate_case}
%\vspace{-2mm}
Analyzing the EM updates for higher dimensions turns out to be
challenging. However, for the symmetric fit in higher dimensions
given by
\begin{align}
\label{eq:multivariate_fit}
  \mixfit{symm}((\theta, \sd^2)) 
    \!=\! \frac{1}{2} \NORMAL(\theta, \sigma^2 I_d)\!  +\!  \frac{1}{2}\NORMAL
    (-\theta, \sigma^2I_d),\! 
\end{align}
the sample EM operator $M_{n, d}(\theta)$ has a closed form as already
noted in the updates~\eqref{eq:sample_em_sigma} and \eqref{eq:sample_em_operators}.
Note that for the fit~\eqref{eq:multivariate_fit}, we have assumed
the same scale parameter for all dimensions.
Such a fit is over-specified for data drawn from Gaussian
distribution $\NORMAL(0, I_\dims)$. 
We now show that the sharing of scale parameter in the model fit across
dimensions~\eqref{eq:multivariate_fit},
leads to a faster convergence of EM in $\dims \geq 2$---both in terms of
number of steps and the final statistical accuracy.
In the following result, we denote $\interval_{\smallthreshold}:=[{5 \parenth{\frac{d}
{ n}}^{\frac{1}{4} + \smallthreshold},\frac{1}
{8}}]$.

\begin{theorem} 
\label{thm:em_multivariate}
Fix $\delta \in (0, 1)$, $\smallthreshold \in (0, 1/ 4]$, and let 
$X_i\stackrel{i.i.d.}{\sim}\mathcal{N}(0, I_\dims)$ for $i=1, \ldots, \obs$
such that $\dims\geq 2$ and  \mbox{$\obs \succsim \dims
\log^{\frac{1}{4 \smallthreshold}}(\log\frac{1/ \smallthreshold}{\delta})$}. Then with any starting
point $\theta_{n}^{0}$ such that $ \enorm{\theta_{n}^{0}} \in
\interval_{ \smallthreshold}$, the sample EM sequence $\theta_{n}^{t +
  1} = M_{n, d}( \theta_{n}^t)$ satisfies
\begin{align}
\label{eq:sample_em_multivariate_rate}
\enorm{ \theta_{n}^{t} - \thetastar} \leq c_1\parenth{ \frac{ \dims}{
    \obs} \log \frac{ \log( 1/ \smallthreshold) } {\delta}}^{\frac{1}{ 4} -
  \smallthreshold},
  \end{align}
  for all $t \geq c_2\parenth{ \frac{\obs}{
    \dims}}^{ \frac{1}{2} - 2 \smallthreshold} \log\frac{\obs}{\dims} \log
    \frac{1}{
    \smallthreshold}$
with probability at least $1 - \delta$.
\end{theorem}
See Appendix~\ref{sub:proof_of_theorem_thm:em_multivariate} for the proof.

The results in Theorem~\ref{thm:em_multivariate} show that the 
that the sample EM updates converge to a ball around $\thetastar = 0$ with
radius arbitrarily close to $(d/ n)^{\frac{1}{4}}$ when $d\geq 2$.
At first sight, the initialization condition $\enorm{ \theta_{n}^{0}} \leq
1/8$, assumed in Theorem~\ref{thm:em_multivariate}, might seem pretty 
restrictive but Lemma~\ref{lemma:one_step_update} (in
Appendix~\ref{sub:proof_of_lemma_lemma:one_step_update}) shows that for any
$\theta^0_n$ satisfying $\enorm{\theta_{n}^{0}} \leq \sqrt{d}$, we have
\mbox{$\PseudoND(\theta_{n}^{0}) \leq \sqrt{2/\pi}$}, with high probability.
In light of this result, we may conclude that the initialization condition
is Theorem~\ref{thm:em_multivariate} is
not overly restrictive. 

\textbf{Comparison with Theorem~\ref{thm:em_univariate}:}
The scaling of order $\obs^{-\frac{1}{4}}$ with $\obs$ is significantly better 
than the univariate case ($\obs^ {-\frac{1} {8}}$) stated in 
Theorem~\ref{thm:em_univariate}. 
We note that this faster statistical rate  is a consequence of the sharing
of the scale parameter across dimensions, and does not hold when the 
fit~\eqref{eq:multivariate_fit} has different variance parameters. Indeed,
as we demonstrated in Figure~\ref{fig:gaussian_vs_mixture}, when the fitted
components have freely varying scale parameter, the
statistical rate slows down (and can be of the order $\obs^{-\frac{1}
{8}}$ in higher dimensions).

\subsection{Simulations with general cases} % (fold)
\label{sub:general_cases}
%\vspace{-2mm}
We now present preliminary evidence that the slow statistical rates of EM with location parameter that we derived
in Theorem~\ref{thm:em_univariate} might appear in more general settings. In Figure~\ref{fig:mixture_vs_mixture}, we plot the statistical
error of estimates returned by sample EM when estimating \emph{all} the 
parameters (namely weights, location and scale) simultaneously, as a function
of sample size $\obs$,
for the following two cases:
\begin{align}
% \textbf{($d=1$): }
  \mix_{\star}^{d=1} &\!=\! \frac{1}{6}\NORMAL(-5, 1) \!+ \!
  \frac{1}{2}\NORMAL(1, 3) \!+\! 
  \frac{1}{3}\NORMAL(7, 2);\!
  %  \notag\\
  % % \qquad\qquad\qquad 
  % \mix_{\text{fit}} 
  % &= \sum_{k=1}^{K_\text{fit}}w_i \NORMAL(\theta_i,
  % \sigma_i^2), 
  % \notag
   \label{eq:three_mixture_one_d}
   \\
  % \textbf{($d=2$): }
  \mix_{\star}^{d=2} &=\! \frac{1}{2}\NORMAL\big({
  \begin{bmatrix} 0 \\ 0 \end{bmatrix} \!,\! \; I}\big)
  +
  \frac{1}{6}\NORMAL\big(\begin{bmatrix}
    7 \\ 5
  \end{bmatrix} \!,\! \; 2I\big) \!  
  % \notag
  % \\ &\qquad+ \!
  \frac{1}{3}\NORMAL\big(\begin{bmatrix}
    -4 \\ -7
  \end{bmatrix}\!,\! \; 3I\big). \!
  % \notag\\
  % \mix_{\text{fit}} 
  % =& \sum_{k=1}^{K_\text{fit}}w_i \NORMAL(\theta_i,  \Sigma_i) 
  \label{eq:three_mixture_two_d}
\end{align}
We plot the results for a $K_\text{fit} \in \braces{3, 4, 5}$-mixture
Gaussian model fit.
When $K_{\text{fit}}$ is equal to the number of components ($=3$) in the
true mixture the statistical rate is $\obs^{ -1/2}$. When it is larger,
i.e.,  $K_{\text{fit}} \in \braces{4, 5}$, the statistical rate of EM is much larger, $\obs^{- 0.12}$
in panel (a) (for $K_{\text{fit}} = 5$) and $\obs^{-0.20}$ in panel (b) (for $K_{\text{fit}} = 5$) of Figure~\ref{fig:mixture_vs_mixture}.
These simulations suggest that the statistical rates slower than $n^{-\frac14}$
and of order $n^{-\frac18}$ may arise in more general settings, and moreover
that the rates get slower as the over-specification of the number of mixtures
increases. See Section~\ref{sec:discussion} for possible future work in this direction.
%\vspace{-2mm}
% subsection general_cases (end)
%%%%%%%%%%%%%%%%%%%%%%%%%%%%%%%%%%%%%%%%%%%%%%%%%%%%%%%%%%%%%%%%%%%%%%%%%%%%%

\section{Analysis of EM} % (fold)
\label{sec:EM_analysis}
%\vspace{-2mm}
Deriving a sharp rate for univariate case (Theorem~\ref{thm:em_univariate})
turns out be pretty challenging
and requires a thorough discussion. On the other hand, the multivariate-case
considered in the paper (Theorem~\ref{thm:em_multivariate}) 
is \emph{relatively} easy due to the shared scale
parameter given the techniques developed in prior works~\cite{Siva_2017,Raaz-misspecified}.
See Appendix~\ref{sub:proof_of_theorem_thm:em_multivariate} for details.
We now outline the analysis of the univariate case
in which we make use of several novel techniques.

\subsection{Proof outline} % (fold)
\label{sub:proof_outline}

Our proof makes use of the population-to-sample analysis framework
of Balakrishnan~\etal~\cite{Siva_2017} albeit with several new ideas.
Let $Y \sim \NORMAL(0, 1)$, then the population-level analog of the 
operator~\eqref{eq:sample_em_operators}
can be defined in two ways:
\begin{subequations}
\begin{align}
\label{eq:pop_one_d}
\Mtil_{\obs,1}( \theta)\! &:=\! \Exs_{Y} \brackets{Y \tanh
  \parenth{ \frac{Y \theta}{ \sum_{j = 1}^n X_j^2/ n\! -\! \theta^2}}},\! \\
  \label{eq:stage_2_pop_em_operator}
\corop( \theta) & \defn \Exs_{Y}\brackets{Y \tanh
    \parenth{ \frac{ Y \theta}{1 - \theta}}}.
\end{align}  
The particular choice of the population-like operator $\Mtil_{n, 1}$ in
equation~\eqref{eq:pop_one_d} was motivated by the previous
works~\cite{Cai_2018} with the location-scale Gaussian mixtures.
We refer to this operator as the \emph{pseudo-population operator} since
it depends on the samples $\braces{X_i, i=1, \ldots n}$ and involves
an expectation.
Nonetheless, as we show in the sequel, analyzing $\Mtil_{n, 1}$ is not enough
to derive sharp rates for sample EM in the over-specified setting considered in
Theorem~\ref{thm:em_univariate}.
A careful inspection reveals that a ``better'' choice of the
population operator is required, which leads us to define
the operator $\corop$ in equation~\eqref{eq:stage_2_pop_em_operator}. 
Unlike the pseudo-population operator $\Mtil_{n, 1}$,
the operator $\corop$ is indeed a population operator as it
does not depend on samples $X_{1}, \ldots, X_{n}$. 
Note that, this operator
is obtained when we replace the sum $\sum_{j = 1}^{n} {X_{j}}^2/n$ in the
definition~\eqref{eq:pop_one_d} of the operator $\Mtil_{n, 1}$
by its corresponding expectation $\Exs[\enorm{X}^2] = 1$.  For this reason,
we also refer to this operator $\corop$ as the \emph{corrected population
  operator}. 
In the next lemma, we state the properties of the operators defined
above (here $\interval_\smallthreshold'$ denotes the interval $[ \unicon
\obs^{- \frac{1}{12}+\smallthreshold}, 1/ 10]$).
%\vspace{-2mm}
\end{subequations}
\begin{lemma}
\label{lemma:correct_pop_contraction_and_radem_three}
  The operators $\Mtil_{n,1}$ and $\corop$ satisfy
  %\vspace{-3.5mm}
\begin{subequations}
  \begin{align}
    \label{eq:lemma_pop_event}    
    \parenth{1 - \frac{3 \theta^6}{2}}\abss{ \theta} & \leq  \abss{
        \Mtil_{\obs,1}( \theta)} \leq \parenth{1 -
      \frac{ \theta^6}{5}}\abss{ \theta}, \\
      \parenth{1 - \frac{ \theta^6}{2}} \abss{\theta} &\leq \abss{
      \corop( \theta)} \leq \parenth{1 - \frac{ \theta^6}{5}}
    \abss{\theta},
    \label{eq:lemma_cor_pop}
\end{align}
where bound~\eqref{eq:lemma_pop_event} holds for all $\abss{ \theta} \in \interval_{\smallthreshold}'$ with high probability\footnote{
  Since the operator $\Mtil_{n,1}$ depends on the samples $\braces{X_j, j
    \in [\obs]}$, only a high probability bound (and not a
  deterministic one) is possible.} and the
bound~\eqref{eq:lemma_cor_pop} is deterministic and holds for all $\abss{ \theta} 
    \in \brackets{0, \frac{3}{20}}$. Furthermore, for any fixed $\delta \in (0, 1)$
    and any fixed $r\geq  O( \obs^{-\frac{1}{12}})$, we have that
  \begin{align}
    \Prob \brackets{ {\sup \limits_{ \theta \in \ball(0, r)} \abss{
          \samm( \theta)- \Mtil_{\obs,1}} \leq \unicon r \sqrt{
          \frac{ \log (1/ \delta)}{n}}}\,} & \nonumber \\
          & \hspace{- 7 em} \geq 1 - \delta.    \label{eq:define_e_r_old}
        \end{align}  
    On the other hand, for any fixed $r \leq O(\obs^{- \myfrac{1}{16}})$,
    we have
    \begin{align}
    \Prob \brackets{ {\sup \limits_{ \theta \in \ball(0, r)} \abss{
          \samm( \theta)- \corop( \theta)} \leq \unicontwo r^3 \sqrt{
          \frac{ \log^{ 10} (5 \obs/ \delta)}{n}}}\,} & \nonumber \\
          & \hspace{- 10 em} \geq 1 - \delta.  \label{eq:define_e_r_new}
  \end{align}
  \end{subequations}
\end{lemma}
See Appendix~\ref{sub:proof_of_lemma_lemma:correct_pop_contraction_and_radem_three}
for its proof where we also numerically verify the sharpness of the results
above (see Figure~\ref{fig:radem}).
Lemma~\ref{lemma:correct_pop_contraction_and_radem_three} establishes that,
as $\theta \to 0$, both the operators have similar contraction coefficient
 $\gamma(\theta)\asymp 1 -
c\theta^6$; thereby justifying the rates observed for $d=1$ in 
Figure~\ref{fig:symm_fit_all_results}(b).
However, their perturbation
bounds are significantly different: while the error
$\sup_{ \theta \in \ball(0, r)} \abss{ \samm( \theta) - \Mtil_{n, 1}(
  \theta)}$ scales linearly with the radius 
  $r$, the deviation error $\sup_{ \theta \in \ball (0, r)}
\abss{\samm( \theta) - \corop( \theta)}$ has a cubic scaling $r^3$.

%%%%%%%%%%%%%%%%%%%%%%%%%%%%%%%%%%%%%%%%%%%%%%%%%%%%%%%%%%%%%%%%%%%%%%%%%%%%%

\textbf{Remark:} 
A notable difference between the two bounds~\eqref{eq:define_e_r_old} and
\eqref{eq:define_e_r_new} is the
range of radius $r$ over which we \emph{prove} the validity of the 
bounds~\eqref{eq:define_e_r_old} and \eqref{eq:define_e_r_new}.
With our tools, we establish that the perturbation 
bound~\eqref{eq:define_e_r_old} for the
operator $\Mtil_{\obs,1}$ is valid for any $r \succsim n^{-\myfrac{1}{12}}$.
On the other hand, the corresponding bound~\eqref{eq:define_e_r_new}
for the operator $\corop$ is valid for any $r \precsim n^{-\myfrac{1}{16}}$.
We now elaborate why these different ranges of radii are helpful and make both the operators crucial to
in the analysis to follow. 
%\vspace{-2mm}
%%%%%%%%%%%%%%%%%%%%%%%%%%%%%%%%%%%%%%%%%%%%%%%%%%%%%%%%%%%%%%%%%%%%%%%%%%%%%

\subsection{ A sub-optimal analysis} 
 \label{subsection:slow_pop}
 %\vspace{-2mm}
Using the properties of the operator $\Mtil_{n, 1}$ from Lemma~\ref{lemma:correct_pop_contraction_and_radem_three}, we
now sketch the statistical rates for
the sample EM sequence, $\theta_n^{t+1}= \samm(\theta_n^t)$, that can be
obtained using (a) the generic procedure
outlined by Balakrishnan et al.~\cite{Siva_2017} and (b) the
localization argument introduced in our previous work~\cite{Raaz_Ho_Koulik_2018}.
As we show, both these arguments end up being \emph{sub-optimal} as they
do not provide us the rate of order $n^{-\myfrac{1}{8}}$ stated in
Theorem~\ref{thm:em_univariate}.  
We use the notation:
\begin{align*} 
\sup_{ \abss{ \theta} \geq \epsilon}  \abss{ \Mtil_{\obs,1}(\theta)}/
{\abss{ \theta}} \precsim 1 - \epsilon^6 \rdefn \gamma(\epsilon).
\end{align*}
%\vspace{-5mm}
%%%%%%%%%%%%%%%%%%%%%%%%%%%%%%%%%%%%%%%%%%%%%%%%%%%%%%%%%%%%%%%%%%%%%%%%%%%%%
  
\noindent\textbf{Sub-optimal rate I:} 
The eventual radius of convergence obtained using Theorem~5(a) from
the paper~\cite{Siva_2017} can be determined by
\begin{subequations}
%\vspace{-4mm}
\begin{align}
\label{eq:siva_simple}
 \frac{r/\sqrt{\obs}}{1 - \gamma(\epsilon)} = \epsilon
  \quad \Longrightarrow \quad
\epsilon \sim \obs^{- 1/ 14},
% %\vspace{-5mm}
\end{align}

where $r$ denotes the bound on the initialization radius $\vert \theta^0\vert$
but we have tracked dependency only on $\obs$.
This informal computation suggests that the the sample EM iterates for location
parameter are bounded by a term of order $n^{- 1/ 14}$. This rate is
clearly sub-optimal when compared to the EM rate of order $n^{-\myfrac{1}{8}}$
from Theorem~\ref{thm:em_univariate}.
%\vspace{-4mm}
%%%%%%%%%%%%%%%%%%%%%%%%%%%%%%%%%%%%%%%%%%%%%%%%%%%%%%%%%%%%%%%%%%%%%%%%%%%%%%%

\noindent\textbf{Sub-optimal rate II:} 
Next we apply the more sophisticated localization argument from the
paper~\cite{Raaz_Ho_Koulik_2018} in order to obtain a sharper rate.
In contrast to the computation~\eqref{eq:siva_simple}, this argument leads to
solving the equation
% \begin{subequations}
\begin{align}
\label{eq:local_simple}
  \frac{\epsilon \cdot r/\sqrt{\obs}}{1 - \gamma(\epsilon)} =
  \epsilon 
  \  \Longrightarrow  \
  \frac{ \epsilon r/ \sqrt{ n}}{ \epsilon^6} = \epsilon \
  \Longrightarrow  \epsilon \sim \obs^{-\myfrac{1}{12}}, \!
\end{align}
where, as before, we have only tracked dependency on $\obs$.
\end{subequations}
This calculation allows us to conclude that the EM algorithm converges
to an estimate which is at a distance of order $n^{-\frac{1}{12}}$ from the
true parameter, which is again
sub-optimal compared to the $\obs^{-\myfrac{1}{8}}$ rate of EM from Theorem~\ref{thm:em_univariate}.

Indeed both the conclusions above can be made rigorous (See Corollary~\ref{cor:convergence_rate_sample_EM_univariate_balanced}
for a formal statement) to conclude that, with high probability for any $\beta \in (0, \frac1{12}]$
\begin{align}
\label{suboptbnd}
  \abss{ \mean_{n}^t - \thetastar}  \leq O(n^{-\frac{1}{12}+\beta})
  \ \text{for} \ t \geq O(n^{\frac12-6\beta}).
\end{align}
%\vspace{-8mm}

%%%%%%%%%%%%%%%%%%%%%%%%%%%%%%%%%%%%%%%%%%%%%%%%%%%%%%%%%%%%%%%%%%%

\subsection{A two-staged analysis for sharp rates} % (fold)
\label{ssub:proof_sketch_for_theorem_thm:em_univariate_a_two_stage_analysis}
%\vspace{-2mm}
In lieu of the above observations, the proof of the sharp upper
bound~\eqref{eq:sample_em_univariate_rate} in
Theorem~\ref{thm:em_univariate} proceeds in two stages.  In the first
stage, invoking
Corollary~\ref{cor:convergence_rate_sample_EM_univariate_balanced} with
$\smallthreshold = \frac{1}{48}$, we
conclude that with high probability the sample EM iterates converge to
a ball of radius at most $r$ after $\sqrt{ \obs}$ steps, where $r
\ll\obs^{- 1/ 16}$.  Consequently, the sample EM iterates after
$\sqrt{\obs}$ steps satisfy the assumptions required to invoke the perturbation
bounds for the operator $\corop$ from 
Lemma~\ref{lemma:correct_pop_contraction_and_radem_three}.  Thereby,
in the second stage of the proof, we apply
the \mbox{$1 - c \theta^6$} contraction bound~\eqref{eq:lemma_cor_pop}
of the operator $\corop$ in conjunction with the cubic perturbation
bound~\eqref{eq:define_e_r_new}. Using localization argument for this
stage, we establish that the EM iterates obtain a statistical error of
order $\obs^ {- 1/ 8}$ in $\order{n^{3/4}}$ steps as stated in
Theorem~\ref{thm:em_univariate}. See Appendix~\ref{sub:proof_of_theorem_thm:em_univariate}
for a detailed proof.
\vspace{-4mm}

%%%%%%%%%%%%%%%%%%%%%%%%%%%%%%%%%%%%%%%%%%%%%%%%%%%%%%%%%%%%%%%%%%%%%%%

%%%%%%%%%%%%%%%%%%%%%%%%%%%%%%%%%%%%%%%%%%%%%%%%%%%%%%%%%%%%%%%%%%%%%%%
%%%%%%%%%%%%%%%%%%%%%%%%%%%%%%%%%%%%%%%%%%%%%%%%%%%%%%%%%%%%%%%%%%%%%%%%%%%%%%%%%%%%% References

\section{Discussion} % (fold)
\label{sec:discussion}
\vspace{-3mm}
In this paper, we established several results characterizing the
convergence behavior of EM algorithm for over-specified location-scale
Gaussian mixtures. We view our analysis of
EM for the symmetric singular Gaussian mixtures
as the first step toward a rigorous
understanding of EM for a broader class of weakly identifiable mixture
models. Such a study would provide a better understanding of the
singular models with weak identifiability which do
arise in practice since: (a) over-specification is a common phenomenon in fitting mixture models due to weak separation between
mixture components, and, (b) the parameters being estimated are often
inherently dependent due to the algebraic structures of the class of
kernel densities being fitted and the associated partial
derivatives. We now discuss a few other directions that can serve as a
natural follow-up of our work.

The slow rate of order $n^{-\myfrac{1}{8}}$ for EM updates with location parameter is in a sense
a worst-case guarantee.  In the univariate case, for the entire class of
two mixture Gaussian fits, MLE exhibits the slowest known statistical rate $n^{-\myfrac{1}{8}}$
for the settings that we analyzed. More precisely, for certain
asymmetric Gaussian mixture fits, the MLE convergence rate for the location parameter is faster than that of the symmetric equal-weighted
mixture considered in this paper
% \footnote{}.
E.g., for the fit $1/3\NORMAL
( - 2 \theta,
\sigma^2) + 2/3\NORMAL( \theta, \sigma^2)$ on $\NORMAL(0, 1)$ data, the
MLE converges at the rate $n^{-{1}/{6}}$ and $n^{-{1}/{3}}$ respectively~\cite{Ho-Nguyen-Ann-16}.
 It is interesting to understand the effect of such a geometric structure
of the global maxima on the convergence of the EM algorithm.

Our work analyzed over-specified mixtures with a specific structure and only
one extra component. As demonstrated above, the statistical rates for EM
appear to be slow for general covariance fits and further appear to slow
down as the number of over-specified components increases. The convergence
rate of the MLE for such over-specified models is known to further deteriorate
as a function of the number of extra components.  It remains to understand
how the EM algorithm responds to these more severe---and practically 
relevant---instances of over-specification.

\section*{Acknowledgments} 
\label{sec:acknowledgments}
This work was partially supported by Office of Naval Research grant
DOD ONR-N00014-18-1-2640 and National Science Foundation grant
NSF-DMS-1612948 to MJW, and National Science Foundation grant  
NSF-DMS-1613002 to BY, and by Army Research Office grant
W911NF-17-1-0304 to MIJ.

% \newpage
\bibliographystyle{plain}
\bibliography{Nhat}

%%%%%%%%%%%%%%%%%%%%%%%%%%%%%%%%%%%%%%%%%%%%%%%%%%%%%%%%%%%%%%%%%%%%%%%%%%%%%%%%%%%%%%%%%%%%%%%%%%%%%%%%%%%%%%%%%%%%%%%%%%%%%%%%%%%%%%%%%%%%%%%
\newpage
\onecolumn

\begin{center}
\textbf{\Large{Supplement for ``Sharp Analysis of Expectation-Maximization
for Weakly Identifiable Models"}}
\end{center}

\appendix
\etocdepthtag.toc{mtappendix}
\etocsettagdepth{mtchapter}{none}
\etocsettagdepth{mtappendix}{subsection}
\etocsettagdepth{mtappendix}{subsubsection}
\tableofcontents

%%%%%%%%%%%%%%%%%%%%%%%%%%%%%%%%%%%%%%%%%%%%%%%%%%%%%%%%%%%%%%%%%%%%%

\section{Proofs of main results} % (fold)
\label{sec:proofs}

In this section, we present the proofs for our main results while deferring
some technical results to the appendices.

\subsection{Proof of Theorem~\ref{thm:em_univariate}} % (fold)
\label{sub:proof_of_theorem_thm:em_univariate}
Our result makes use of the following corollary (proven in Appendix~\ref{sub:proof_of_prop_EM_univariate_balanced}):
\begin{corollary} 
\label{cor:convergence_rate_sample_EM_univariate_balanced}
Given constants $\delta \in (0, 1)$ and $\smallthreshold \in (0, 1/
12]$, suppose that we generate the the sample-level EM sequence $
  \mean_{n}^{t + 1} = \samm( \mean_{n}^t)$ starting from an
  initialization \mbox{$\vert{ \mean_{n}^{0}}\vert
  \in \interval'_\smallthreshold$}, and
  using a sample size $\samples$ lower bounded as \mbox{$\samples
    \succsim \log^{1/ (12 \smallthreshold)}(\log(1/ \smallthreshold)/
    \delta)$}.  Then for all iterations 
    \mbox{$t \geq \obs^{ 1/ 2 -6\smallthreshold} \log(n ) \log(1/ \smallthreshold)$},
    we have
\begin{align}
\label{eq:prop_rate}
    \vert{ \mean_{n}^t - \thetastar} \vert \leq c_1\parenth{
      \frac{1}{n} \log \frac{ \log( 1/ \smallthreshold)}{
        \delta}}^{\myfrac{1}{12} - \smallthreshold},
\end{align}
  with probability at least $1- \delta$.
\end{corollary}
\paragraph{Remark:}
We note that the sub-optimal bound~\eqref{eq:prop_rate}
obtained from
Corollary~\ref{cor:convergence_rate_sample_EM_univariate_balanced}
is not an artifact
of the localization argument and arises due to the definition
of the operator operator $\Mtil_{\obs, 1}$~\eqref{eq:pop_one_d}.  
As we have alluded to earlier, indeed a finer analysis with
the population EM operator $\corop$ is required to prove the rate of $\obs^
{-1/8}$ stated in Theorem~\ref{thm:em_univariate}.
However, a key assumption in the further derivation is that
the sample EM iterates $ \theta_{n}^{t}$ can converge to a ball of
radius $r \precsim \obs^{- 1/ 16}$ around $\thetastar$ in a finite
number of steps, for which Corollary~\ref{cor:convergence_rate_sample_EM_univariate_balanced}
comes in handy.

We now begin with a sketch the two stage-argument, and then provide
a rigorous proof
for Theorem~\ref{thm:em_univariate}.

\subsubsection{Proof sketch} % (fold)
\label{ssub:proof_sketch_of_theorem}

As mentioned earlier, the pseudo-population operator
$\Mtil_{\obs,1}$ is not sufficient to achieve the sharp rate of EM
iterates under the univariate symmetric Gaussian mixture fit. Therefore,
we make use of corrected-population operator
$\corop$ to get a sharp statistical rate of EM. Our proof for the tight 
convergence rate of sample EM updates relies on a novel two-stage localization
argument that we are going to sketch.

\paragraph{First stage argument:} 

Plugging in $\smallthreshold = 1/ 84$ in
Corollary~\ref{cor:convergence_rate_sample_EM_univariate_balanced}, we
obtain that for $t \succsim \sqrt{ \obs} \log ( \obs)$, with
probability at least $1 - \delta$ we have that
\begin{align}
\label{eq:cor_1_thm_3_bound}
  \abss{\theta_n^t - \thetastar} \leq \unicon
  \obs^{-\frac{1}{14}}\log^{\frac{1}{14}} \frac{ \log( 1/
    \smallthreshold)}{ \delta} \leq \obs^{- \frac{ 1}{ 16}},
\end{align}
where the second inequality follows from the large sample condition
$\obs \geq \uniconnew\log^8\frac{\log84}{\delta}$. All the following
claims are made conditional on the event~\eqref{eq:cor_1_thm_3_bound}.

%%%%%%%%%%%%%%%%%%%%%%%%%%%%%%%%%%%%%%%%%%%%%%%%%%%%%%%%%%%%%%%%%%%%%%%%%%%%%

\paragraph{Second stage argument:} 

In order to keep the presentation of the proof sketch simple, we do
not track constant and logarithmic factors in the arguments to follow.
In epoch $\ind$, for any iteration $t$ the EM iterates satisfy
$\theta_n^t\in [n^{- a_{\ind+1}},n^{- a_\ind} ]$ where $a_{\ind+1} >
a_{\ind}$ and $a_{\ind} \leq 1/ 16$. Applying
Lemma~\ref{lemma:correct_pop_contraction_and_radem_three} for such
iterations, we find that with high probability
\begin{align*}
\abss{ \corop( \theta_n^t )} \precsim \underbrace{ (1 - n^{- 6{
      a_{\ind+1}}})}_{ \rdefn \gamma_{\ind}} \abss{ \theta_n^t } \quad
\text{ and} \quad \abss{ \samm( \theta_n^t ) - \corop( \theta_n^t)}
\precsim \frac{n^{- 3 a_\ind}}{\sqrt{n}},
\end{align*}
where the first bound follows from the $1-c\theta^6$ contraction
bound~\eqref{eq:lemma_cor_pop} and the second bound follows from the
cubic-type Rademacher bound~\eqref{eq:define_e_r_new}.  Invoking the
basic triangle inequality $T$ times, we
obtain that
\begin{align*}
\abss{ \theta_n^{t + T}} \stackrel{(i)}{ \precsim} e^{- T n^{- 6
    {a_{\ind+1}}} } n^{- a_\ind} + \frac{1} {1 - \gamma_\ind} \cdot
\frac{n^{- 3 a_\ind}}{\sqrt{n}} \stackrel{(ii)}{\precsim} \frac{ 1}{1
  - \gamma_\ind} \cdot \frac{ n^{- 3 a_\ind}}{\sqrt{ n}} = n^{6
  a_{\ind+1} - 3 a_\ind - 1/2},
\end{align*}
where in step~(ii) we have used the fact that for large enough $T$,
the first term is dominated by the second term in the RHS of step~(i). To obtain a recursion for the sequence
$a_\ell$, we set the RHS equal to $n^{-a_{\ell+1}}$.  Doing so yields
the recursion 

\begin{subequations}
\begin{align}
\label{eq:a_sequence_define}
  a_{\ind+1} = \frac{3 a_\ind}{ 7} + \frac{ 1}{ 14}, \quad\text{where}\quad
  a_0 = 1/16.
\end{align}
Solving for the limit $a_{\ind+1} = a_{\ind} = a_\star$, we find that
$a_\star = 1/ 8$. Thus, we can conclude that sample EM iterates in the 
univariate setting converge to a ball of radius $n^{-1/8}$ as claimed
in the theorem statement.
\subsubsection{Formal proof of sample EM convergence rate}
\label{ssub:formal_proof_of_theorem}
We now turn to providing a formal proof for the
preceding arguments.

%%%%%%%%%%%%%%%%%%%%%%%%%%%%%%%%%%%%%%%%%%%%%%%%%%%%%%%%%%%%%%%%%%%%%%%%%%%%%%%
\paragraph{Notations:} To make the proof comprehensible, some additional notations
are necessary which we collect here.
Let 
$\lstar = \lceil \log( 8/ \smallthreshold)/ \log(7/ 3) \rceil$
so that $a_{ \lstar} \leq 1/ 8 - \smallthreshold$.
We define the following shorthand:
\begin{align}
\label{eq:defn_stat_err}
  \staterr \mydefn \frac{n}{c_{n, \delta}}, 
  \quad \text{ where} \quad 
  c_{n, \delta} \mydefn \log^{10}(10 \obs(\lstar + 1)/ \delta).
\end{align}
For $\ind=0, \ldots, \lstar$, we define the time sequences $t_\ind$ and $T_\ind$
as follows:
\begin{align}
  \label{eq:t_sequence_define}
  t_0 = \sqrt{n}, 
  \quad t_\ind = \ceil{10 \staterr^{6a_\ind}
  \log \staterr}, \quad \text{ and} \quad 
  T_\ind = \sum_{j = 0}^\ind t_j.
\end{align}
Direct computation leads to
\begin{align}
\label{eq:t_i_star}
T_{\lstar} \leq \sqrt{ \obs} + \lstar t_{\lstar} \precsim \log
\parenth{ \frac{ \obs \log \frac{1}{ \smallthreshold} }{ c_{n,\delta}
    \delta} } \parenth{ \frac{n}{c_{ n, \delta}}} ^{3/ 4 - 6
  \smallthreshold} \precsim \obs^{3/4}.
\end{align}
In order to facilitate the proof argument later, we define the
following set
\begin{align}
\label{eq:radii}
  \mathcal{R} \mydefn \braces{ \staterr^{- a_1}, \ldots, \staterr^{-
      a_{\lstar}}, \uniconnew \staterr^{- a_1}, \ldots, \uniconnew
    \staterr^{- a_{\lstar}}},
\end{align}
\end{subequations}
where $\uniconnew \mydefn (5 c_2 + 1)$. Here, $c_2$ is the universal
constant from Lemma~\ref{lemma:correct_pop_contraction_and_radem_three}.

%%%%%%%%%%%%%%%%%%%%%%%%%%%%%%%%%%%%%%%%%%%%%%%%%%%%%%%%%%%%%%%%%%%%%%%%%%%%

\paragraph{Formal argument:} We show that with probability at least $1 - \delta$
the following
holds:
\begin{align}
\label{eq:proof_general_claim}
  \abss{ \theta_n^t} \leq \parenth{ \frac{c_{ \obs, \delta}}{n}}^{
    a_\ind} = \staterr^{- a_\ind}, \quad \text{ for all } t \geq
  T_\ind, \text{ and } \ind \leq \lstar.
\end{align}
As a consequence of this claim and the
definitions~\eqref{eq:a_sequence_define}-\eqref{eq:t_i_star} of
$a_{\lstar}$ and $T_{ \lstar}$, we immediately obtain that
\begin{align*}
  \vert{ \mean_{n}^t - \thetastar}\vert \precsim  
  \parenth{ \frac{c_{n , \delta}}{ n}}^{1/ 8 - \smallthreshold} 
  \precsim \parenth{ \frac{1}{n} \log^{ 10} \frac{ 10 n \log( 8/ 
  \smallthreshold)}{ \delta}}^{1/ 8 - \smallthreshold},
\end{align*}
for all number of iterates $t \succsim \obs^{3/ 4 - 6 \smallthreshold}
\log( n) \log( 1/ \smallthreshold)$ with probability at least $1 -
\delta$ as claimed in Theorem~\ref{thm:em_univariate}.

We now define the high probability event that is crucial for our
proof.  For any $r \in \mathcal{ R}$, define the event $E_r$ as
follows
\begin{align*}
E_r & \mydefn \braces{ \sup \limits_{ \theta \in \ball(0, r)} \abss{
    \samm( \theta) - \corop( \theta)} \leq \unicon_{2} r^3 \sqrt{
    \frac{ \log^{10}(5 \obs \abss{ \mathcal{ R}}/ \delta)}{n}}}.
\end{align*}
Then, for the event
\begin{align}
\label{eq:intersection_events}
  \event \defn \bigcap_{r \in \mathcal{ R}} E_r \cap
  \braces{\text{Event } \eqref{eq:cor_1_thm_3_bound} \text{ holds }},
\end{align}
applying the union bound with
Lemma~\ref{lemma:correct_pop_contraction_and_radem_three} yields that
$\Prob[ \event] \geq 1 - \delta$.  All the arguments that follow are
conditional on the event $\event$ and hence hold with the claimed high
probability.

In order to prove the claim~\eqref{eq:proof_general_claim}, we make
use of the following intermediate claim:
\begin{lemma}
  \label{lemma:non_expansive}
Conditional on the event $\event$, if $\abss{ \theta} \leq \staterr^{-
  a_\ind}$, then $\abss{ \samm( \theta)} \leq \staterr^{- a_\ind}$ for
any $\ind \leq \lstar$.
\end{lemma}
\noindent Deferring the proof of Appendix~\ref{sub:proof_of_lemma_lemma:non_expansive}, we now establish the claim~\eqref{eq:proof_general_claim}
conditional on the event $\event$ only for $t = T_\ind$ and when $\vert{
\theta_n^t}\vert \in [\staterr^{- a_{\ind + 1}}, \staterr^{- a_\ind}]$ 
in which we now prove using induction.

%%%%%%%%%%%%%%%%%%%%%%%%%%%%%%%%%%%%%%%%%%%%%%%%%%%%%%%%%%%%%%%%%%%%%%%%%%%%%

\paragraph{Proof of base case $\ind = 0$:} 

Note that we have $a_0 = 1/ 16$ and that $ \obs^{- 1/ 16} \leq
\staterr^{ 1/ 16}$.  Also, by the
definition~\eqref{eq:intersection_events} we have that $\text{the
  event}~ \eqref{eq:cor_1_thm_3_bound} \subseteq \event$.  Hence,
under the event $\event$ we have that $\abss{ \theta_n^{t}} \leq
\obs^{- 1/ 16}$, for $t \succsim \sqrt{ \obs} \log ( \obs)$.  Putting
all the pieces together, we find that under the event $\event$, we
have $\abss{ \theta_n^{t}} \leq \obs^{- 1/ 16} \leq \staterr^{ 1/ 16}$
and the base case follows.

%%%%%%%%%%%%%%%%%%%%%%%%%%%%%%%%%%%%%%%%%%%%%%%%%%%%%%%%%%%%%%%%%%%%%%%%

\paragraph{Proof of inductive step:} 

We now establish the inductive step.  Note that
Lemma~\ref{lemma:non_expansive} implies that we need to show the
following: if $\abss{ \theta_n^{ t}} \leq \staterr^{- a_{ \ind}}$ for
all $t \in \braces{ T_\ind, T_\ind + 1, \ldots, T_{\ind + 1} - 1}$ for
any given $\ind \leq \lstar$, then $\vert{ \theta_n^{T_{\ind + 1}}}
\vert \leq \staterr^ {- a_{\ind + 1}}$.  We establish this claim in
two steps:
\begin{subequations}
\begin{align}
  \theta_n^{T_{ \ind} + t_{ \ind}/ 2} & \leq c' \staterr^{- a_{\ind
      + 1}}, \quad\text{and},
  \label{eq:first_step}\\
  \theta_n^{T_{\ind + 1}} & \leq \staterr^{- a_{\ind +
      1}}, \label{eq:second_step}
\end{align}
where $c' = (5c_2 + 1) \geq 1$ is a universal constant.  Note that the
inductive claim follows from the bound~\eqref{eq:second_step}.  It
remains to establish the two claims~\eqref{eq:first_step} and
\eqref{eq:second_step} which we now do one by one.
\end{subequations}

%%%%%%%%%%%%%%%%%%%%%%%%%%%%%%%%%%%%%%%%%%%%%%%%%%%%%%%%%%%%%%%%%%%%%%%%%%%%

\paragraph{Proof of claim~\eqref{eq:first_step}:} 

Let $\Theta_\ind = \braces{ \theta: \abss{ \theta} \in [\staterr^{-
      a_{\ind + 1}}, \staterr^{- a_{\ind}}]}$.  Now, conditional on
the event $\event$,
Lemma~\ref{lemma:correct_pop_contraction_and_radem_three} implies that
\begin{align*}
  \sup_{\theta \in \Theta_\ind} \abss{ \samm( \theta) -
    \corop(\theta)} \leq c_2 \staterr^{- 3 a_{\ind} - 1/ 2}, \quad
  \text{and} \quad \sup_{\theta \in \Theta_\ind} \abss{ \corop(
    \theta)/ \theta} \leq (1 - \staterr^{- 6 a_{\ind +1}}/ 5) \rdefn
  \gamma_\ind.
\end{align*}
We can check that $\gamma_\ind \leq e^{- \staterr^{ 6 a_{\ind + 1}}/
  5}$.  Unfolding the basic triangle inequality $t_\ind/ 2$ times and noting that $\theta_n^t \in \Theta_\ind$ for all
$t \in \braces{T_\ind, \ldots, T_\ind + t_\ind/ 2}$, we obtain that
\begin{align*}
  \abss{ \theta_n^{ T_\ind + t_\ind/ 2}} & \stackrel{}{ \leq}
  \gamma_\ind^{ t_\ind/ 2} \abss{ \theta_n^{ T_\ind}} + (1+\gamma_\ind
  + \ldots + \gamma_\ind^{t_\ind/ 2 - 1}) c_2 \staterr^{- 3 a_{\ind} -
    1 / 2} \\ & \stackrel{}{ \leq} e^{- t_\ind \staterr^{- 6 a_{\ind +
        1}}/ 10} \staterr^{- a_\ind} + \frac{1}{1 - \gamma_\ind} c_2
  \staterr^{- 3 a_{\ind} - 1 / 2} \\
& \stackrel{(i)}{ \leq} (1 + 5c_2) \staterr^{6 a_{\ind +1} - 3
    a_{\ind} - 1/ 2} \\
& \stackrel{(ii)}{ =} (5 c_2 + 1) \staterr^{- a_{\ind + 1}}
\end{align*}
where step~(i) follows from plugging in the value of $\gamma_\ind$ and
invoking the definition~\eqref{eq:t_sequence_define} of $t_\ind$,
which leads to
\begin{align*}
  e^{- t_\ind \staterr^{6 a_{\ind + 1}}/ 10} \staterr^{- a_\ind}
  \leq \staterr^{6 a_{\ind +1} - 3 a_{\ind} - 1/ 2}.
\end{align*}
Moreover, step~(ii) is a direct consequence of the
definition~\eqref{eq:a_sequence_define} of the sequence
$a_\ind$. Therefore, we achieve the conclusion of
claim~\eqref{eq:first_step}.

%%%%%%%%%%%%%%%%%%%%%%%%%%%%%%%%%%%%%%%%%%%%%%%%%%%%%%%%%%%%%%%%%%%%%%%%%%%%%%%

\paragraph{Proof of claim~\eqref{eq:second_step}:} 

The proof of this step is very similar to the previous step, except
that we now use the set $\Theta_\ind' = \braces{ \theta: \abss{
    \theta} \in [\staterr^{- a_{\ind + 1}}, c'\staterr^{- a_{\ind +
        1}}]}$ for our arguments.  Applying
Lemma~\ref{lemma:correct_pop_contraction_and_radem_three}, we have
\begin{align*}
  \sup_{\theta \in \Theta'_\ind} \abss{ \samm( \theta) - \corop(
    \theta)} \leq c_2 (c')^3 \staterr^{- 3 a_{\ind + 1} - 1/ 2}, \quad
  \text{and} \quad \sup_{ \theta \in \Theta'_\ind} \abss{ \corop(
    \theta)/ \theta} \leq \gamma_\ind.
\end{align*}
Using the similar argument as that from the previous case, we find that
\begin{align*}
  \abss{ \theta_n^{ T_\ind + t_\ind/ 2 + t_\ind/ s2 } } & \leq e^{-
    t_\ind \staterr^{6 a_{\ind +1}}/ 10} c' \staterr^{ - a_{\ind + 1}
  } + \frac{1}{1 - \gamma_\ind} c_2 (c')^3 \staterr^{- 3 a_{\ind + 1}
    - 1/ 2} \\ & \leq (5 c_2 + 1) (c')^3 \staterr^{ 4 a_{\ind + 1} -
    1/ 2 } \cdot \staterr^{- a_{\ind + 1}} \\ & \stackrel{(i)}{ \leq}
  \staterr^{- a_{\ind + 1}}
\end{align*}
where step~(i) follows from the inequality $e^{- t_\ind \staterr^{6 a_{\ind +1}}/ 10} \leq \staterr^{ 4 a_{\ind + 1} - 1/ 2 }$ and the inequality
\begin{align*}
\staterr^{4 a_{\ind + 1} - 1/ 2} \leq \staterr^{4 a_{\lstar} - 1/ 2}
\leq \staterr^{- 4 \smallthreshold} \leq 1/ (c')^4,
\end{align*}
since $n \geq (c')^{1 / \smallthreshold} c_{n, \delta}$.  The claim
now follows.
% subsection proof_of_theorem_thm:em_univariate (end)

%%%%%%%%%%%%%%%%%%%%%%%%%%%%%%%%%%%%%%%%%%%%%%%%%%%%%%%%%%%%%%%%%%%%%%%

\subsection{Proof of Theorem~\ref{thm:em_multivariate}} % (fold)
\label{sub:proof_of_theorem_thm:em_multivariate}
Before proceeding further, we first derive the  convergence rates for the
scale parameter $\sd^t_n$ using
Theorem~\ref{thm:em_multivariate}.
Noting that $(\thetastar, \sdstar)=(0, 1)$,
we obtain the following relation
\begin{align*}
 % \label{eqn:theta_sigma_upadte_relation}    
\abss{ ( \sd_n^{t})^2 - ( \sd^{*})^2} = \abss{ \frac{ \sum_{i =
      1}^\obs \enorm{ X_i}^2}{ d \obs} - ( \sd^{*})^2 - \frac{\enorm{
          \theta^{t}_n - \thetastar}^2}{d}}.
\end{align*}
Using standard chi-squared bounds, we obtain that
\begin{align*}
  \abss{\frac{\sum_{i =1}^\obs\enorm{ X_i}^2}{ d \obs} - (\sd^{*})^2} \precsim
  ({nd})^{-\myfrac{1}{2}},
\end{align*}
with high probability. 
From the bound~\eqref{eq:sample_em_multivariate_rate}, we also have $\enorm{
          \theta^{t}_n - \thetastar}^2/d \precsim ({nd})^{-\myfrac{1}{2}}
          $.
Putting the pieces together, we conclude that the statistical
error for the scale parameter satisfies 
\begin{align}
\label{eq:sigma_multivariate}
 \vert { ( \sd_n^
{t})^2 -
( \sd^{*})^2}\vert \precsim ({nd})^{-\myfrac{1}{2}}\quad 
\mbox{for all $t \succsim\parenth{\frac{n}{d}}^{\myfrac{1}{2}}$,}
\end{align}
with high probability.
Consequently, in the sequel, we focus primarily on the convergence rate
for the EM estimates $ \theta^{t}_n$ of the location parameter, as the 
corresponding guarantee for the scale parameter $\sd^t_n$ is readily implied
by it.

The proof of Theorem~\ref{thm:em_multivariate} is based on the 
population-to-sample analysis and follows a similar road-map
as of the proofs in the paper~\cite{Raaz_Ho_Koulik_2018}.
We first analyze the population-level EM operator and
then using epoch-based-localization argument derive the statistical
rates~\eqref{eq:sample_em_multivariate_rate}.
We make use of the following $d$-dimensional
analog of the pseudo-population operator (cf. equation~\eqref{eq:pop_one_d}):
\begin{align}
\label{eq:pop_like_em_operator}
  \PseudoND( \theta) \defn \Exs_{Y \sim \Ncal(0, I_\dims)} \brackets{Y
    \tanh \parenth{ \frac{ Y^\top \theta}{\sum_ {j = 1}^\obs
        \enorm{X_j}^2/( \obs \dims) - \Vert \theta \Vert^2/ \dims}}}.
\end{align}
In the next lemma, we establish the contraction properties
and the perturbation bounds for $\PseudoND$:
\begin{subequations}
\begin{lemma}
  \label{lemma:d_dim_operator}
  The operator $\PseudoND$ satisfies
  \begin{align}
    \label{eq:pop_em_multivariate_contraction_bound}
    \parenth{1 - \frac{3 \enorm{ \theta}^{2}}{4}} \leq \frac{ \enorm{
        \PseudoND(\theta)} } { \enorm{ \theta}} & \leq \parenth{1 -
      \frac{ (1 - 1/ d) \enorm{ \theta}^{2}}{4}}, \quad \text{ for all
    } \enorm{ \theta} \in \interval_{\smallthreshold},
    \end{align}
with probability at least $1 - \delta$. Moreover, 
  there exists a universal constant $\unicontwo$ such that for any fixed
  $\delta \in (0, 1)$, $\smallthreshold \in (0, \frac{1}{4}]$, and $r\in (0, \frac{1}{8})$ we have
\begin{align}
\label{eq:define_e_r_multi}
  \Prob \brackets{ \sup \limits_{\theta \in \ball(0,r)}  
  \enorm{M_{n, d}(\theta) - \PseudoND( \theta)} \leq 
  \unicontwo r\sqrt{\frac{d\log(1/\delta)}{n}}} 
  \geq  1 - \delta - e^{- (n d)^{ 4 \smallthreshold}/ 8}.
\end{align}
\end{lemma}
\end{subequations}
\noindent See Appendix~\ref{sec:proof_of_lemma_lemma:d_dim_operator} for
the proof.\\

Lemma~\ref{lemma:d_dim_operator} shows that the operator $\PseudoND$ has
a faster contraction (order $1-\enorm{\theta}^2$) towards zero, when compared
to its univariate-version (order $1-\theta^6$ cf.~\eqref{eq:lemma_pop_event}).
This difference between the univariate and the multivariate case had already
been highlighted in Section~\ref{sub:some_illustrative_examples} in 
Figure~\ref{fig:symm_fit_all_results}. Indeed substituting $\dims=1$ in the
bound~\eqref{eq:pop_em_multivariate_contraction_bound} gives us a vacuous
bound for the univariate case, providing further evidence for the benefit of 
sharing variance among different dimensions in multivariate setting of symmetric fit~\eqref{eq:symmetric_fit}. With Lemma~\ref{lemma:d_dim_operator} at hand, the proof of Theorem~\ref{thm:em_multivariate}
follows by using the localization argument from the paper~\cite{Raaz_Ho_Koulik_2018}.
Mimicking the arguments similar to equation~\eqref{eq:local_simple}, we
obtain the following statistical rate:\footnote{Moreover,
similar to the arguments made in the paper~\cite{Raaz_Ho_Koulik_2018}, 
localization argument is necessary to derive a sharp rate. Indeed, a direct
application of the framework introduced by Balakrishnan~\etal~\cite{Siva_2017}
for our setting implies a sub-optimal rate of order $(d/n)^{1/6}$ for the
Euclidean error $\Vert{\theta_{n}^{ t} - \thetastar}\Vert$ (cf.~\eqref{eq:siva_simple}
and \eqref{eq:local_simple}).}
\begin{align}
\label{eq:local_simple_d}
  \frac{\epsilon \cdot r/\sqrt{\obs}}{1 - \gamma(\epsilon)} =
  \epsilon
  \quad \Longrightarrow 
  \frac{ \epsilon r/ \sqrt{ n}}{ \epsilon^2} = \epsilon \quad
  \Longrightarrow \quad \epsilon \sim \obs^{-\myfrac{1}{4}}.
\end{align}

Much of the work in the proof of Theorem~\ref{thm:em_multivariate}
is to establish Lemma~\ref{lemma:d_dim_operator}.
With the bounds~\eqref{eq:pop_em_multivariate_contraction_bound} and \eqref{eq:define_e_r_multi}
at hand, using the localization
argument (in a manner similar to the proof of Theorem~\ref{thm:em_univariate}),
easily leads to the statistical rate of order $(d/n)^{1/4}$ as claimed in
Theorem~\ref{thm:em_multivariate}. The detailed proof is thereby omitted.
% subsection proof_of_theorem_thm:em_multivariate (end)

%%%%%%%%%%%%%%%%%%%%%%%%%%%%%%%%%%%%%%%%%%%%%%%%%%%%%%%%%%%%%%%%%%%%%%%

\subsection{Proof of Lemma~\ref{lemma:correct_pop_contraction_and_radem_three}}
\label{sub:proof_of_lemma_lemma:correct_pop_contraction_and_radem_three}

We now prove Lemma~\ref{lemma:correct_pop_contraction_and_radem_three} which
provides the basis for the two-staged proof of Theorem~\ref{thm:em_univariate}.

The proof for the contraction property~\eqref{eq:lemma_cor_pop} of the corrected
population operator $\corop$ is similar to that of the property~\eqref{eq:lemma_pop_event}
pseudo-population operator $\Mtil_{\obs,1}$ (albeit with a few high probability
arguments replaced by deterministic arguments). Hence, while we provide
a complete proof of the bound~\eqref{eq:lemma_pop_event} (in Section~\ref{ssub:proof_of_contraction_pseudo_EM_uni}),
we only provide a proof sketch for the bound~\eqref{eq:lemma_cor_pop}
at its end.
Moreover the proofs of bounds~\eqref{eq:define_e_r_old} and \eqref{eq:define_e_r_new}
are provided in Sections~\ref{ssub:proof_of_perturbation_bound_for_Mtil}
and \ref{ssub:proof_of_perturbation_bound_for_corop} respectively.

%%%%%%%%%%%%%%%%%%%%%%%%%%%%%%%%%%%%%%%%%%%%%%%%%%%%%%%%%%%%%%%%%%%%%%%

\subsubsection{Contraction bound for population operator $\Mtil_{\obs,1}$}
\label{ssub:proof_of_contraction_pseudo_EM_uni}

We begin by defining some notation. For $\smallthreshold \in (0, 1/
12]$ and \mbox{$\alpha \geq 1/ 2 - 6 \smallthreshold$}, we define the
  event $\event_\alpha$ and the interval $\interval_{\alpha,
    \smallthreshold}$ as follows\:
\begin{align}
\label{eq:event_x_square}
  \event_\alpha & = \braces{ \bigg \vert \sum_{j = 1}^\obs X_j^2/ \obs
    - 1 \bigg \vert \leq {\obs^{- \alpha}}}, \quad \text{ and}, \\
  \label{eq:interval_definition}
  \interval_{\alpha, \smallthreshold} &= [3 n^{- 1/ 12 +
      \smallthreshold}, \sqrt{9/ 400 - n^{- \alpha}}],
\end{align}
where in the above notations we have omitted the dependence on $\obs$,
as it is clear from the context.  We also use the scalars $a$ and $b$
to denote the following:
\begin{align*}
a \defn 1 - n^{-\alpha} \quad \text{ and } \quad b \defn 1 +
n^{-\alpha}.
\end{align*}
With the above notation in place, observe that standard chi-squared
tail bounds yield that $ \Prob[ \event_\alpha] \geq 1 - e^{- n^{1 - 2
    \alpha}/ 8} \geq 1 - \delta$.  Moreover, invoking the lower bound
on $n$ in Theorem~\ref{thm:em_univariate}, we have that $[3
  n^{- 1/ 12 + \smallthreshold}, 1/ 10] \subseteq \interval_{\alpha,
  \smallthreshold}$.  Now conditional on the high probability event
$\event_\alpha$, the population EM update $\Mtil_{\obs,1}(\theta)$, in absolute value,
can be upper and lower bounded as follows:

\begin{align*}
\abss{\Mtil_{\obs,1}( \theta)} & \leq \Exs_Y \brackets{ Y \tanh \parenth{
    \frac{ Y \abss{ \theta}}{a - \theta^2}}} = \abss{ \theta}
\underbrace{ \Exs_Y \brackets{ \frac{Y} { \abss{ \theta}} \tanh
    \parenth{ \frac{ \abss{ \theta} X} {a - \theta^2}}}}_{ \rdefn
  \unicontraup( \theta)}, \quad \text{and}, \\
\abss{\Mtil_{\obs,1}( \theta)} & \geq \Exs_Y \brackets{ Y \tanh \parenth{
    \frac{ X \abss{ \theta}}{b - \theta^2}}} = \abss{ \theta}
\underbrace{ \Exs_Y \brackets{ \frac{Y} {\abss{ \theta}} \tanh
    \parenth{ \frac{ \abss{ \theta} Y}{b - \theta^2}}}}_{ \rdefn
  \unicontralo( \theta)},
\end{align*} 
where the last two inequalities follows directly from the definition
of $\Mtil_{\obs,1}(\theta)$ in equation~\eqref{eq:pop_one_d}, and from the
fact that for any fixed $y, \theta \in \real$, the function $w \mapsto
y\tanh(y\abss{\theta}/(w-\theta^2))$ is non-increasing in $w$ for
$w>\theta^2$.  Consequently, in order to complete the proof, it
suffices to establish the following bounds:
\begin{align}
\label{eq:gamma_bound}
   1 -3 \theta^6/ 2 \leq \unicontralo( \theta), \quad \text{and} \quad
   \unicontraup( \theta) \leq (1 - \theta^6/ 5).
\end{align}
The following properties of the hyperbolic function $x \mapsto x
\tanh(x)$ are useful for our proofs:
\begin{lemma} 
\label{lemma:inequality_tanh_sech_function}
For any $x \in \Rspace$, the following holds
\begin{align}
\text{(Lower bound):} \quad & x \tanh(x) \geq x^2 - \frac{x^4}{3} +
\frac{2 x^6}{15} - \frac{17 x^8}{ 315}, \nonumber \\ \text{(Upper
  bound):} \quad & x \tanh(x) \leq x^2 - \frac{ x^4}{ 3} + \frac{2
  x^6}{ 15} - \frac{17 x^8}{315} + \frac{62 x^{10}}{ 2835}. \nonumber
  \end{align}
\end{lemma}
\noindent See Appendix~\ref{sub:proof_of_lemma_univariate_balanced}
for its proof.

Given the bounds in Lemma~\ref{lemma:inequality_tanh_sech_function},
we derive the upper and lower bounds in the
inequality~\eqref{eq:gamma_bound} separately.

%%%%%%%%%%%%%%%%%%%%%%%%%%%%%%%%%%%%%%%%%%%%%%%%%%%%%%%%%%%%%%%%%%%%%%%%

\paragraph{Upper bound for $\unicontraup(\theta)$:} 

Invoking the upper bound on $x\tanh(x)$ from
Lemma~\ref{lemma:inequality_tanh_sech_function}, we find that
\begin{align*}
\unicontraup( \theta) & \leq \frac{a - \theta^{2}}{ \theta^2} \biggr(
\frac{ \theta^2}{(a - \theta^2)^2} \Exs \brackets{ Y^2} - \frac{
  \theta^4}{3 (a - \theta^2)^4} \Exs \brackets{ Y^4} + \frac{ 2
  \theta^6}{15 (a - \theta^2)^6} \Exs \brackets{ Y^6} \nonumber
\\ & \hspace{ 12 em} - \frac{17 \theta^8}{315 (a - \theta^2)^8} \Exs
\brackets{ Y^8} + \frac{62 \theta^{10}}{2835 (a - \theta^2)^{10}} \Exs
\brackets{ Y^{10}} \biggr).
\end{align*}
Recall that, for $Y \sim \Ncal(0,1)$, we have $\Exs \brackets{ Y^{2
    k}} = (2 k - 1)!!$ for all $k \geq 1$. Therefore, the last
inequality can be simplified to
\begin{align}
\label{eq:gamma_simple_bound}  
\unicontraup( \theta) \leq \frac{1}{a - \theta^2} - \frac{
  \theta^2}{(a - \theta^2)^3} + \frac{2 \theta^4}{(a - \theta^2)^5} -
\frac{17 \theta^6}{3 (a - \theta^2)^7} + \frac{62 \theta^8}{3 (a -
  \theta^2)^9}.
\end{align}
When $n^{- \alpha} + \theta^2 \leq 9/ 400$, we can verify that the
following inequalities hold:
\begin{align*}
\frac{1}{1 - n^{- \alpha} - \theta^2} & \leq 1 + (n^{- \alpha} +
\theta^2) + (n^{- \alpha} + \theta^2)^2 + (n^{- \alpha} + \theta^2)^3
+ 2( n^{- \alpha} + \theta^2)^4, \\
- \frac{ \theta^2}{(1 - n^{- \alpha}- \theta^2)^3} & \leq - \theta^2
\parenth{1 + 3 (n^{- \alpha} + \theta^2) + 6 (n^{- \alpha} +
  \theta^2)^2 + 10 (n^{- \alpha} + \theta^2)^3}, \\
\frac{ \theta^4}{(1 - n^{- \alpha} - \theta^2)^5} & \leq \theta^4
\parenth{1 + 5 (n^{- \alpha} + \theta^2) + 16 (n^{- \alpha} +
  \theta^2)^2}, \\
- \frac{ \theta^6}{(1 - n^{- \alpha} - \theta^2)^7} & \leq - \theta^6
\parenth{1 + 7 (n^{- \alpha} + \theta^2)}, \\
\frac{ \theta^8}{(1 - n^{- \alpha} - \theta^2)^9} & \leq 5 \theta^8/
4. \nonumber
\end{align*} 
Substituting $a = 1-n^{-\alpha}$ into the
bound~\eqref{eq:gamma_simple_bound} and doing some algebra with the
above inequalities and using the fact that $\max \braces {\theta, n^{-
    \alpha}} \leq 1$ we have that
\begin{align*}
\unicontraup( \theta) \leq 1 - \frac{2}{3} \theta^6 + \frac{61}{6}
\theta^8 + 100n^{-\alpha} \leq 1 - \frac{2}{5} \theta^6 + 100
n^{-\alpha} \leq 1 - \frac{1}{5} \theta^6.
\end{align*}
The second last inequality above follows since $\theta \leq 3/20$, and
the last inequality above utilizes the fact that if $\alpha \geq 1/ 2
- 6 \smallthreshold$, then $\theta^6/ 5 \geq 100 \obs^{- \alpha}$ for
all \mbox{$\theta \geq 3 \obs^{- 1/ 12 + \smallthreshold}$.} This
completes the proof of the upper bound of $\unicontraup(\theta)$.

%%%%%%%%%%%%%%%%%%%%%%%%%%%%%%%%%%%%%%%%%%%%%%%%%%%%%%%%%%%%%%%%%%%%%%%

\paragraph{Lower bound for $\unicontralo(\theta)$:}  

We start by utilizing the lower bound of $x\tanh(x)$ in the expression
for $\unicontralo(\theta)$, which yields:
\begin{align}
\label{eq:gamma_low_simple_bound}
\unicontralo(\theta) \geq \frac{1}{b - \theta^2} - \frac{ \theta^2}{(b
  - \theta^2)^3} + \frac{2 \theta^4}{(b - \theta^2)^5} - \frac{17
  \theta^6}{3 (b - \theta^2)^7}.
\end{align}
Since $\abss{ \theta} \in [3 n^{- 1/ 12 + \smallthreshold}, \sqrt{9/
    400 - n^{- \alpha}}]$ by assumption, we have the following lower
bounds:
\begin{align*}
\frac{1}{1 + n^{- \alpha} - \theta^2} & \geq 1 + ( \theta^2 - n^{-
  \alpha}) + (\theta^2 - n^{- \alpha})^2 + (\theta^2 - n^{- \alpha})^3
+ (\theta^2 - n^{- \alpha})^4, \\
- \frac{ \theta^2}{(1 + n^{- \alpha} \theta^2)^3} & \geq - \theta^2
  - \parenth{1 + 3 (\theta^2 - n^{- \alpha}) + 6 (\theta^2 - n^{-
  - \alpha})^2 + 11( \theta^2 - n^{- \alpha})^3}, \\
\frac{\theta^4}{(1 + n^{- \alpha} - \theta^2)^5} & \geq \theta^4
\parenth{1 + 5 (\theta^2 - n^{- \alpha}) + 15 (\theta^2 - n^{-
    \alpha}}, \\
- \frac{\theta^6}{(1 + n^{- \alpha} - \theta^2)^7} & \geq - \theta^6
\parenth{1 + 8 (\theta^2 - n^{- \alpha})}.
\end{align*}
Substituting $b = 1 + n^{-\alpha}$ into the
bound~\eqref{eq:gamma_low_simple_bound} and doing some algebra with
the above inequalities and using the fact that $\max \braces {\theta,
  n^{- \alpha}} \leq 1$ we have that
\begin{align*}
\unicontralo( \theta) \geq 1 - \frac{2}{3} \theta^6 - \frac{76}{3}
\theta^8 - 100n^{- \alpha} \geq 1 -\frac{5}{4} \theta^6 -100 n^{-
  \alpha} \geq 1- \frac{3}{2} \theta^6,
\end{align*}
The second last inequality above follows since $\theta \leq 3/20$, and
the last inequality above utilizes the fact that if $\alpha \geq 1/ 2
- 6 \smallthreshold$, then $\theta^6/ 4 \geq 100 \obs^{- \alpha}$ for
all \mbox{$\theta \geq 3 \obs^{- 1/ 12 + \smallthreshold}$.} This
completes the proof of the lower bound of $\unicontralo(\theta)$.

%%%%%%%%%%%%%%%%%%%%%%%%%%%%%%%%%%%%%%%%%%%%%%%%%%%%%%%%%%%%%%%%%%%%%%%
\paragraph{Proof of contraction bound for $\corop$:} 

Note that it suffices to repeat the arguments with $a=1$ and $b=1$
in the RHS of the inequalities~\eqref{eq:gamma_simple_bound} and \eqref{eq:gamma_low_simple_bound}
respectively. Given the other computations, the remaining steps are 
straightforward algebra and are thereby omitted.
%%%%%%%%%%%%%%%%%%%%%%%%%%%%%%%%%%%%%%%%%%%%%%%%%%%%%%%%%%%%%%%%%%%%%%%

\subsubsection{Proof of perturbation bound for $\Mtil_{n, 1}$} % (fold)
\label{ssub:proof_of_perturbation_bound_for_Mtil}
We now prove the bound~\eqref{eq:define_e_r_old} which is based on standard 
arguments to derive Rademacher complexity bounds. 
We first symmetrize with Rademacher variables, and apply the Ledoux-Talagrand
contraction inequality.  We then invoke results on sub-Gaussian and
sub-exponential random variables, and finally perform the associated
Chernoff-bound computations to obtain the desired result.

To ease the presentation, we denote $\alpha \mydefn 1/ 2 - 2 \beta$
and $\mathcal{I} \mydefn [1 - \obs^{- \alpha} - 1/ 64, 1 - \obs^{-
    \alpha}]$. Next we fix $r \in [ 0, 1/ 8]$ and define
$\widetilde{r} \defn \frac{r}{1 - \obs^{- \alpha} - 1/ 64}$.  For
sufficiently large $n$, we have $\newr \leq 2 r$. 
Recall the definition~\eqref{eq:event_x_square} of the event:
$\event_\alpha = \{  \vert \sum_{j = 1}^\obs
X_j^2/ \obs - 1  \vert \leq {\obs^{- \alpha}}\}$.
Conditional on the event $\event_\alpha$, the following inequalities hold
\begin{align*}
\abss{\samm( \theta)  -  \Mtil_{\obs,1}( \theta)}
      & \leq \sup_{\theta \in \ball(0,  r), \sigma^2 \in \mathcal{I}} 
      \abss{ \frac{1}{n} \sum_{i = 1}^n X_i \tanh \parenth{ \frac{X_i
      \theta}{ \sigma^2}}
      - \Exs \left[ Y \tanh \parenth{ \frac{ Y \theta}{ \sigma^2}}
      \right]} \\
      & \leq \sup_{\newtheta \in \ball(0, \newr)}
      \abss{ \newm_n( \newtheta) - \newm( \newtheta)},
\end{align*}
with all them valid for any $\theta \in \ball(0, r)$.  Here $Y$
denotes a standard normal variate $\Ncal( 0, 1)$ whereas the operators
$\newm$ and $\newm_n$ are defined as
\begin{align*}
\newm( \newtheta) 
      \mydefn \Exs[ Y \tanh( Y \newtheta)]  
      \quad \text{and} \quad
\newm_n( \newtheta) 
      \mydefn \frac{1}{ n} \sum_{i = 1}^n X_i \tanh( X_i 
      \newtheta).
\end{align*}
To facilitate the discussion later, we define the unconditional 
random 
variable
\begin{align*}
Z \mydefn \sup_{ \newtheta \in \ball (0, \newr)}
      \abss{ \newm_n( \newtheta) - \newm( \newtheta)}.
\end{align*}
Employing standard symmetrization argument from empirical 
process theory~\cite{Vaart_Wellner_2000}, we find that
\begin{align*}
\Exs[ \exp( \lambda Z)] 
      \leq \Exs \brackets{ \exp \parenth{ \sup_{ \newtheta \in
      \ball (0, \newr)} \frac{2 \lambda}{ n} 
      \sum_{i = 1}^n \radem_i \tanh( X_i \newtheta) X_i}},
\end{align*}
where $\radem_i, i \in [ n]$ are i.i.d. Rademacher random 
variables 
independent of $\braces{X_i, i\in[n]}$. Noting that, the 
following 
inequality with hyperbolic function $\tanh(x)$ holds
\begin{align*}
\abss{ \tanh( x \newtheta) - \tanh( x \newtheta')} 
      \leq  \abss{( \newtheta - \newtheta') x} 
      \quad \text{ for all } x.
\end{align*}
Consequently for any given $x$, the function $\newtheta \mapsto \tanh(
x \newtheta)$ is Lipschitz. Invoking the Ledoux-Talagrand contraction
result for Lipschitz functions of Rademacher
processes~\cite{Ledoux_Talagrand_1991} and following the proof
argument from Lemma 1 in the paper~\cite{Raaz_Ho_Koulik_2018}, we
obtain that
\begin{align*}
Z \leq c \newr \sqrt{ \frac{ \log( 1/ \delta)}{ n}}, 
      \quad \text{ with probability } \geq 1 - \delta,
\end{align*}
for some universal constant $c$.
Finally, using $\newr \leq 2r$ for large $n$, we obtain that
\begin{align*}
\abss{ \samm( \theta) - \Mtil_{\obs,1}( \theta)} 
      \leq 2 c r \sqrt{ \frac{ \log(1/ \delta)}{ n}}, 
      \quad \text{ with probability }
      \geq 1 - \delta - e^{- n^{1 - 2 \alpha}/ 8},
\end{align*}
where we have also used the fact that $\Prob[\event_\alpha] \geq 1 -
e^{- n^{1 - 2 \alpha}/ 8}$ from standard chi-squared tail bounds. The 
bound~\eqref{eq:define_e_r_old} follows and we are done.

%%%%%%%%%%%%%%%%%%%%%%%%%%%%%%%%%%%%%%%%%%%%%%%%%%%%%%%%%%%%%%%%%%%%%%%

\subsubsection{Proof of perturbation bound for $\corop$} % (fold)
\label{ssub:proof_of_perturbation_bound_for_corop}

We now prove the bound~\eqref{eq:define_e_r_new}.
Note that it suffices to establish
the following point-wise result:
\begin{align*}
\abss{ \corop( \theta)- \samm( \theta)} 
    \precsim \frac{\abss{ \theta}^3 \log^{ 10} (5 n / \delta)} 
    {\sqrt{\obs}}
    \quad \mbox{for all $\quad \abss{\theta} \precsim 
    \obs^{- 1/ 16}$},
\end{align*}
with probability at least $1 - \delta$ for any given $\delta > 0$. For
the reader's convenience, let us recall the definition of these
operators
\begin{subequations}
\begin{align}
\label{eq:recall_cor_pop}
\corop(\theta) & = \Exs \Big[ X \tanh( X \theta/ (1 - \theta^2))
  \Big], \\
\label{eq:recall_sample_op}
\samm(\theta) & = \frac{1}{ \obs} \sum_{i = 1}^n X_i \tanh \Big( X_i
\theta/ (a_n - \theta^2) \Big),
\end{align}  
\end{subequations}
where $a_n \defn \sum_{i=1}^\obs X_i^2/n$. We further denote $\popmean{\idx} \defn
\Exs_{X\sim\NORMAL(0, 1)}[ X^ {\idx}]$, and $\sammean{\idx} \defn
\frac{1}{\obs}\sum_{i = 1}^n X_i^{ \idx}$.  From known results on
Gaussian moments, we have \mbox{$\popmean{2\idx} = (2k-1)!!$} for each
integer $k = 1, 2, \ldots$.

For any given $x$ and scalar $\scalartwo$, consider the map $\theta
\mapsto x \tanh( x \theta/ (\scalartwo - \theta^2))$.  The $9$-th
order Taylor series for this function around $\theta=0$ is given by
\begin{align}
x \tanh( x \theta/ (\scalartwo - \theta^2)) & = \frac{ \theta
  x^2}{\scalartwo} - \frac{ \theta^3 (x^4 - 3 \scalartwo x^2) }{3
  \scalartwo^3} + \theta^5 \parenth{\frac{2 x^6}{15 \scalartwo^5} -
  \frac{x^4} {\scalartwo^4} + \frac{x^2}{\scalartwo^3}} \notag\\ & +
\theta^7 \parenth{ -\frac{17 x^8}{315 \scalartwo^7} + \frac{2 x^6}{3
    \scalartwo^6} - \frac{2 x^4}{\scalartwo^5} +
  \frac{x^2}{\scalartwo^4} } \notag\\ & + \theta^9 \parenth{\frac{62
    x^{10}}{2835 \scalartwo^9} - \frac{17 x^8}{45 \scalartwo^8} +
  \frac{2 x^6}{\scalartwo^7} - \frac{10 x^4}{3 \scalartwo^6} +
  \frac{x^2}{\scalartwo^5} } + \remainder,
      \label{eq:ninth_order_expansion}
\end{align}
where the remainder $\remainder$ satisfies $\remainder \leq \order{\theta^{11}}$.
Plugging in this expansion with $\scalartwo=1$ on RHS of equation~\eqref{eq:recall_cor_pop}
and taking expectation over $X \sim \NORMAL(0, 1)$, we obtain
\begin{subequations}
  \begin{align}
  \label{eq:cor_pop_expansion}
    \corop( \theta) 
    & = \theta + \theta^3 \big(
      \sum_{\idx = 1}^2 \constant_{3, \idx} \popmean{2 \idx}\big)
      + \theta^5 \big(
      \sum_{\idx = 1}^3 \constant_{5, \idx} \popmean{2 \idx}\big)
      + \theta^7 \big(
      \sum_{\idx = 1}^4 \constant_{7, \idx} \popmean{2 \idx}\big)
      + \theta^9 \big(
      \sum_{\idx = 1}^5 \constant_{9, \idx} \popmean{2 \idx}\big)
      + \remainder,
  \end{align}
  where we have used the notation $\popmean{\idx} := \Exs_{X\sim\NORMAL
  (0, 1)}[ X^{\idx}]$ and $\constant_{j, k}$ denote universal
  constants.
  Furthermore, plugging in the same expansion~\eqref{eq:ninth_order_expansion} 
  with $\scalartwo=a_n$ on RHS of equation~\eqref{eq:recall_sample_op},
  we obtain the following expansion for the sample EM operator
  \begin{align}
  \label{eq:sample_em_expansion}
    \samm( \theta) 
    & = \theta + \theta^3 \big(
      \sum_{\idx = 1}^2 \constant_{3, \idx} 
      \frac{\sammean{2 \idx}}{a_n^{1 + k}} \big) + 
      \theta^5 \big( \sum_{\idx = 1}^3 \constant_{5, \idx} 
      \frac{ \sammean{2 \idx}}{a_n^{2 + k}} \big)
      + \theta^7 \big( \sum_{\idx = 1}^4 \constant_{7, \idx} 
      \frac{ \sammean{2 \idx}}{a_n^{3 + k}} \big) 
      + \theta^9 \big( \sum_{\idx = 1}^5 \constant_{9, \idx} 
      \frac{\sammean{2 \idx}}{a_n^{4 + k}} \big)
      % & \hspace{ 33 em}
       + \remainder_n,
  \end{align}
  where $\sammean{\idx}$ denotes the sample mean of $X^k$,
  i.e., $\sammean{\idx} \defn \frac{1} {\obs}\sum_{i = 1}^n X_i^{\idx}$.
\end{subequations}
In order to lighten the notation, we introduce the following
convenient shorthand:
\begin{align}
\label{eq:defn_betas}
\simconst_j = \sum_{\idx = 1}^{\frac{j + 1}{2}} \constant_{j, \idx}
\popmean{2 \idx} \quad \text{and} \quad \simconsthat_j = \sum_{\idx =
  1}^{\frac{j + 1}{2}} \constant_{j, \idx} \frac{ \sammean{2
    \idx}}{a_n^{\frac{j - 1}{2} + \idx}} \quad \text{for } j \in
\braces{3, 5, 7, 9} \rdefn \jset.
\end{align}
A careful inspection reveals that $\simconst_{3} = \simconst_5 =
0$. With the above notations in place, we find that
\begin{align*}
\abss{ \corop( \theta) - \samm( \theta)} & = \big \vert{\sum_{j \in
    \jset} \theta^j (\simconst_j - \simconsthat_j)} \big \vert +
\remainder \\ & \rdefn \term_1 + \term_2.
\end{align*}
Therefore, it remains to establish that
\begin{align}
\label{eq:simple_bounds_u1_u2}
\term_1 \precsim \frac{ \abss{ \theta}^3 \log^{ 5} (5 n / \delta)} {
  \sqrt{\obs}} \quad \text{and} \quad \term_2 \precsim \frac{\abss{
    \theta}^3 \log^{ 5} (5 n / \delta)}{ \sqrt{\obs}},
\end{align}
with probability at least $1 - \delta$ for any given $\delta > 0$.
Since the remainder term is of order $\theta^{11}$, the assumption
$\abss{\theta}\precsim n^{-1/16}$ ensures that the remainder term is
bounded by a term of order $\theta^3/\sqrt{\obs}$ and thus the
bound~\eqref{eq:simple_bounds_u1_u2} on the second term $\term_2$
follows.

We now use concentration properties of Gaussian moments in order to
prove the bound~\eqref{eq:simple_bounds_u1_u2} on the first term
$\term_1$. Since $\abss{\theta}\leq 1$, it suffices to show that
\begin{align}
\label{eq:bound_betas}
\sup_{j \in \jset} \abss{ \simconst_j - \simconsthat_j } 
    \precsim \frac{\log^{ 5} (5 n / \delta)} {\sqrt{n}}
\end{align}
with probability at least $1 - \delta$.
Using the relation~\eqref{eq:defn_betas}, we
find that
\begin{align}
\abss{\simconst_j - \simconsthat_j }
      = \big \vert \sum_{\idx = 1}^{\frac{j + 1}{2}}          
      \parenth{ \constant_{j, \idx} \popmean{2 \idx}
        - \constant_{j, \idx} 
        \frac{ \sammean{2 \idx}}{a_n^{\frac{j - 1}{2} + \idx}}} 
        \big \vert
      & \leq \sum_{\idx = 1}^{\frac{j + 1}{2}} 
      \frac{\constant_{j, 
      \idx}}{a_n^{\frac{j - 1}{2} + \idx}} 
      \big \vert{ \popmean{2 \idx} - \sammean{2 \idx} } 
      \big \vert + \constant_{j, \idx} (1 - a_n^{- \frac{j - 1} 
      {2} - \idx}) \popmean{2 \idx} \nonumber \\
      & \leq C \sum_{\idx = 1}^{\frac{j + 1}{2}} \parenth{\big \vert{  
      \popmean{2 \idx} - \sammean{2 \idx}
        } \big \vert + \frac{ \popmean{2 \idx}}{\sqrt{\obs}}},
        \label{eq:beta_diff_bounds}
\end{align}
for any $j \in \jset$.
Here in the last step we have used the following bounds:
\begin{align*}
\max_{j \in \jset, \idx \leq \frac{j + 1}{2}} \constant_{j, \idx} 
    \leq C \quad \text{and} \quad
\max_{j \in \jset, \idx \leq \frac{j + 1}{2}}(1 -a_n^{- \frac{j - 1}{2} -
    \idx}) \leq \frac{C}{\sqrt{ \obs}}
\end{align*}
for some universal constant $C$.  Thus a lemma for the
$1/\sqrt{n}$-concentration\footnote{The bound from
  Lemma~\ref{lemma:concentration_high_order} is sub-optimal for $k=1$
  but is sharper than the standard tail bounds for Gaussian
  polynomials of degree $2k$ for $k\geq 2$. The $1/\sqrt{n}$
  concentration of higher moments is necessary to derive the sharp
  rates stated in our results.} of higher moments of Gaussian random
variable is now useful:
\begin{lemma} \label{lemma:concentration_high_order}
Let $X_{1}, \ldots, X_{n}$ are i.i.d. samples from $
\Ncal(0,1)$ and let $\popmean{2\idx} := \Exs_{X\sim\NORMAL(0, 1)}[ X^{2\idx}]$
and $\sammean{2\idx} : = \frac{1}{\obs}\sum_{i = 1}^n X_i^{ 2\idx}$.
Then, we have
\begin{align*}
\Prob \parenth{ \abss{ \sammean{2 \idx} - \mu_{2 \idx}} 
    \leq  \frac{C_{ \idx} \log^{ \idx}(n/ \delta)}{\sqrt{ n}}} 
    \geq 1 - \delta 
  \quad\text{for any $\quad k \geq 1$},
\end{align*}
where $C_{\idx}$ denotes a universal constant depending only on $\idx$.
\end{lemma}
\noindent See the Appendix~\ref{sub:proof_of_lemma_concentration_high_order}
for the proof.

For any $\delta > 0$, consider the event
\begin{align}
\label{eq:gassian_event}
\event \mydefn 
    \left\{ \big \vert{ \popmean{2 \idx} - \sammean{2 \idx}} \big \vert 
    \leq  \frac{C_{ \idx} \log^{ \idx}(5 n/ \delta)}{\sqrt{ \obs}}
        \quad \text{ for all } \idx \in
        \braces{2,
        4, \ldots, 10} \right\}.
\end{align}
Straightforward application of union bound with 
Lemma~\ref{lemma:concentration_high_order} yields that 
$\Prob \brackets{ \event} \geq 1 - \delta$.
conditional on the event $\event$ inequality~\eqref{eq:bound_betas}
implies that
\begin{align}
  \sup_{j \in \jset} \abss{ \simconst_j - \simconsthat_j }
  &\leq C \sup_{j \in \jset}
  \sum_{\idx = 1}^{\frac{j + 1}{2}} \parenth{\big \vert{  
      \popmean{2 \idx} - \sammean{2 \idx}
        } \big \vert + \frac{ \popmean{2 \idx}}{\sqrt{\obs}}}\notag\\
  &\leq  C \sup_{j \in \braces{3, 5, 7, 9}} \frac{j + 1}{2}
   \parenth{\big \vert{  
      \popmean{j+1} - \sammean{j+1}
        } \big \vert + \frac{ (j+1)!!}{\sqrt{\obs}}}\notag\\
  &\stackrel{(i)}{\leq}  C \sup_{j \in \braces{3, 5, 7, 9}} (j-1)
   \parenth{\big \vert{ 
        C_{\frac{j + 1}{2}} \frac{\log^{\frac{j + 1}{2}}(5n/\delta)}{\sqrt{n}}
        } \big \vert + \frac{ (j+1)!! }{\sqrt{\obs}}}\notag\\
  &\stackrel{(ii)}{\leq} C \frac{\log^{ 5} (5 n / \delta)}{ 
      \sqrt{ \obs}},
\end{align}
where step~(i) follows from the definition of the event~\eqref{eq:gassian_event}
and in step~(ii) using the fact that $j \leq 9$ is bounded we absorbed
all the constants into a single constant.
Since the event $\event$ has probability at least $1-\delta$, the 
claim~\eqref{eq:bound_betas} now follows.

\subsubsection{Sharpness of bounds of Lemma~\ref{lemma:correct_pop_contraction_and_radem_three}}
\begin{figure}[t]
  \begin{center}
    \begin{tabular}{c}
    \widgraph{0.45\textwidth}{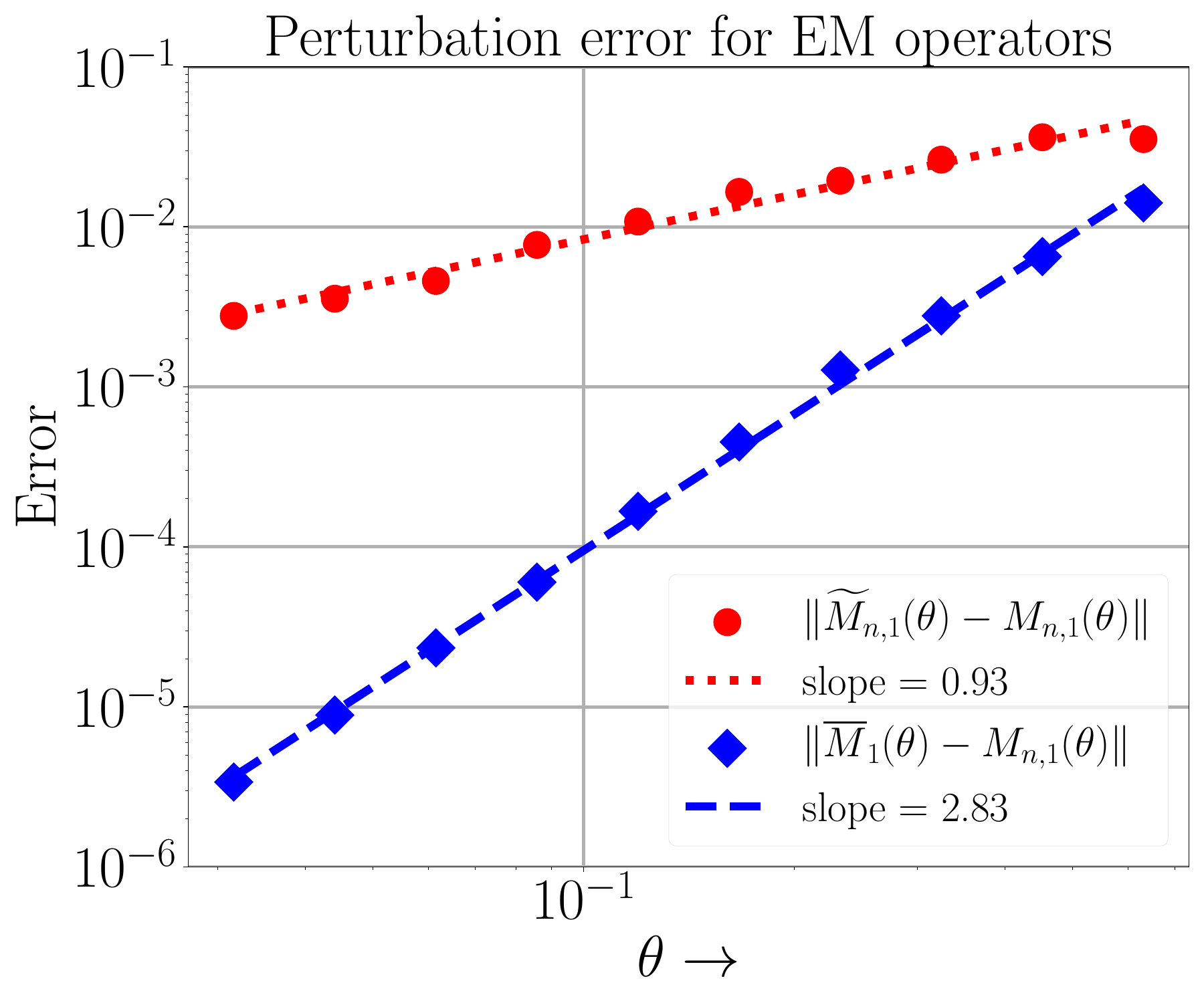}
    \end{tabular}
    \caption{Plots of the perturbation errors for the
      pseudo-population operator $\Mtil_{n, 1}$~\eqref{eq:pop_one_d} and
      the corrected population operator
      $\corop$~\eqref{eq:stage_2_pop_em_operator} with respect to the
      sample EM operator $\samm$~\eqref{eq:sample_em_simple_operator},
      as a function of $\theta$. From the least-squares fit on the
      log-log scale, we see that the error $\Vert{\Mtil_{n, 1}(
        \theta) - \samm}(\theta) \Vert$ scales linearly with $\theta$,
      the error $\Vert{\corop( \theta) - \samm( \theta)} \Vert$ has a
      cubic dependence on $\theta$, in accordance
      with Lemma~\ref{lemma:correct_pop_contraction_and_radem_three}.}
    \label{fig:radem}
  \end{center}
\end{figure}

In Figure~\ref{fig:radem}, we numerically verify the linear and cubic scaling
of the bounds stated in Lemma~\ref{lemma:correct_pop_contraction_and_radem_three}.

%%%%%%%%%%%%%%%%%%%%%%%%%%%%%%%%%%%%%%%%%%%%%%%%%%%%%%%%%%%%%%%%%%%%%%%%%%%%%%%%%%%%%%%%%%%%
% section proofs (end)

%%%%%%%%%%%%%%%%%%%%%%%%%%%%%%%%%%%%%%%%%%%%%%%%%%%%%%%%%%%%%%%%%

% \appendix

\section{Minimax bound} % (fold)
\label{sec:minimax_bound}

We now show that the error of order $\obs^{-\frac{1}{8}}$ (up to logarithmic
factors) is, in fact, tight in the standard minimax sense. Given a compact
set $\Omega \subset \Rspace \times (0, \infty)$, and a set of true parameters
$(\theta^*, \sd^*) \in \Omega$, suppose that we draw $n$
i.i.d. samples $\{X_i\}_{i=1}^n$ from a two-Gaussian mixture of the form
$\frac{1}{2}\mathcal{N}(\theta^*, (\sd^*)^2) + \frac{1}{2}\mathcal{N}(
-\theta^*,(\sd^*)^2)$.
Let $(\estloca, \estsca) \in \Omega$ denote any estimates---for
the respective parameters---measurable with respect to the observed samples
$X_{1}, \ldots, X_{n} \stackrel{\text{i.i.d.}}
{\sim}f_{\theta^*, \sd^*}$ and let $\Exs_{ (\theta^*, \sigma^*)}$ denote
the corresponding expectation.
\begin{proposition}
  \label{prop:minimax_lower_bound}
There exists a universal constant $c_\Omega > 0$ (depending only on $\Omega$),
such that 
\begin{align*}
\inf_{(\estloca, \estsca)} \sup_{(\theta^*, \sigma^*)}
\Exs_{(\theta^*, \sigma^*)} \brackets{ \big({\vert{\estloca}\vert -
    \abss{\theta^*}}\big)^2 + \abss{(\estsca)^2 - (\sigma^*)^2}} \geq
c_\Omega n^{- \frac{1}{4} - \delta} \qquad \mbox{for any $\delta >
  0$.}
\end{align*}
\end{proposition}
\noindent See Appendix~\ref{sub:proof_of_prop_EM_univariate_balanced} for
the proof.

Based on the connection between location parameter $\theta_{n}^{t}$ and scale parameter $\sd_{n}^{t}$ in the EM updates (cf. Equation~\eqref{eq:sample_em_operators}), the minimax lower bound in Proposition~\ref{prop:minimax_lower_bound}
shows that the (non-squared) error of EM location updates $|\vert{\theta^{t}_n}\vert -
\abss{\theta^*}|$ is lower bounded by a term (arbitrarily close to)
$\obs^{-\myfrac{1}{8}}$. 
% section minimax_bound (end)
\subsection{Proof of Proposition~\ref{prop:minimax_lower_bound}} % (fold)
\label{sec:proof_of_proposition_prop:minimax_lower_bound}
We now present the proof of the minimax bound.
We introduce the shorthand $v \mydefn
\sigma^{2}$ and $\eta \mydefn (\theta, v)$.  First of all, we claim
the following key upper bound of Hellinger distance between mixture
densities $f_{\eta_{1}}$, $f_{\eta_{2}}$ in terms of the distances
among their corresponding parameters $\eta_{1}$ and $\eta_{2}$:
\begin{align}
    \label{eqn:minimax_bound_key}  
\inf_{\eta_{1}, \eta_{2} \in \Omega} \frac{h \parenth{f_{\eta_{1}},
    f_{\eta_{2}}}}{\parenth{ \parenth{\abss{\theta_{1}} -
      \abss{\theta_{2}}}^{2} + \abss{ v_{1} - v_{2}}}^r} = 0 \qquad
\mbox{for any $r \in (1,4)$.}
\end{align}
Moreover, for any two densities $p$ and $q$, we denote the total variation
distance between $p$ and $q$ by $V(p ,q) := (1/ 2) \int \abss{p(x) - q(x)}
dx$. Similarly, the squared Hellinger
distance between $p$ and $q$ is given as $h^2(p, q) = (1/ 2) \int
\parenth{\sqrt{p(x)} - \sqrt{q(x)}}^2 dx$.

Taking the claim~\eqref{eqn:minimax_bound_key} as given for the
moment, let us complete the proof of
Proposition~\ref{prop:minimax_lower_bound}. Our proof relies on Le
Cam's lemma for establishing minimax lower bounds. In
particular, for any $r \in (1,4)$ and for any $\epsilon > 0$
sufficiently small, according to the result in
equation~\eqref{eqn:minimax_bound_key}, there exist $\eta_{1} =
(\theta_{1}, v_{1})$ and $\eta_{2} = (\theta_{2}, v_{2})$ such that
$\parenth{\abss{\theta_{1}} - \abss{\theta_{2}}}^{2} + \abss{ v_{1} -
  v_{2}} = 2 \epsilon$ and $h \parenth{f_{\eta_{1}}, f_{\eta_{2}}}
\leq \unicon \epsilon^{r}$ for some universal constant $\unicon$.
From Lemma 1 from Yu~\cite{Yu-1997}, we obtain that
\begin{align*}
\sup_{ \eta \in \{\eta_{1}, \eta_{2}\}} \Exs_{\eta} \brackets{
  \parenth{\abss{\estloca} - \abss{\theta}}^2 + \abss{(\estsca)^2 -
    (\sigma)^2}} \gtrsim \epsilon \parenth{1 - V(f_{\eta_{1}}^{n},
  f_{\eta_{2}}^{n})},
\end{align*}
where $f_{\eta}^{n}$ denotes the product of mixture densities
$f_{\eta}$ of the data $X_{1}, \ldots, X_{n}$. A standard relation
between total variation distance and Hellinger distance leads to
\begin{align*}
V(f_{\eta_{1}}^{n}, f_{\eta_{2}}^{n}) \leq h(f_{\eta_{1}}^{n},
f_{\eta_{2}}^{n}) = \sqrt{1 - \brackets{1 - h^2(f_{\eta_{1}},
    f_{\eta_{2}})}^n} \leq \sqrt{1 - \brackets{1 - \unicon
    \epsilon^{r}}^{n}}.
\end{align*}
By choosing $c \epsilon^{r} = 1/ n$, we can verify that
\begin{align*}
\sup_{ \eta \in \{\eta_{1}, \eta_{2} \}} \Exs_{\eta} \brackets{
  \parenth{\abss{\estloca} - \abss{\theta}}^2 + \abss{(\estsca)^2 -
    (\sigma)^2}} \gtrsim \epsilon \asymp n^{- 1/ r},
\end{align*}
which establishes the claim of
Proposition~\ref{prop:minimax_lower_bound}.

%%%%%%%%%%%%%%%%%%%%%%%%%%%%%%%%%%%%%%%%%%%%%%%%%%%%%%%%%%%%%%%%%%%%%%%

\subsubsection{Proof of claim~\eqref{eqn:minimax_bound_key}} % (fold)
\label{ssub:proof_of_claim_eqn:minimax_bound_key}

In order to prove claim~\eqref{eqn:minimax_bound_key}, it is
sufficient to construct sequences $\eta_{1, n} = (\theta_{1, n}, v_{1,
  n})$ and $\eta_{2, n} = (\theta_{2, n}, v_{2, n})$ such that
\begin{align*}
h \parenth{f_{\eta_{1, n}}, f_{\eta_{2, n}}} \big/
\parenth{\parenth{\abss{ \theta_{1, n}} - \abss{\theta_{2, n}}}^{2} +
  \abss{ v_{1, n} - v_{2, n}}}^r \to 0
\end{align*}
as $n \to \infty$. Indeed, we construct these sequences as follows:
$\theta_{2, n} = 2 \theta_{1, n}$ and $v_{1, n} - v_{2, n} = 3
\parenth{ \theta_{1, n}}^{2}$ for all $n \geq 1$ while $\theta_{1, n}
\to 0$ as $n \to \infty$.  Direct computation leads to
\begin{align*}
f_{\eta_{1, n}}(x) - f_{\eta_{2, n}}(x) = \frac{1}{2}
\underbrace{\parenth{ \normDensity(x; - \theta_{1, n}, v_{1, n}) -
    \normDensity(x; - \theta_{2, n}, v_{2, n})} }_{T_{1, n}} +
\frac{1}{2} \underbrace{ \parenth{ \normDensity(x; \theta_{1, n},
    v_{1, n}) - \normDensity(x; \theta_{2, n}, v_{2, n})}}_{T_{2, n}}.
\end{align*}
Invoking Taylor expansion up to the third order, we obtain that
\begin{align*}
T_{1, n} & = \sum_{\abss{ \alpha} \leq 3} \frac{(\theta_{2, n} -
  \theta_{1, n})^{ \alpha_{1}} (v_{1, n} - v_{2, n})^{ \alpha_{2}}}
{\alpha_{1}! \alpha_{2}!} \frac{ \partial^{ \abss{ \alpha}}
  {\normDensity}}{\partial{ \theta}^{ \alpha_{1}} \partial{v}^{
    \alpha_{2}}} (x; - \theta_{2, n}, v_{2, n}) + R_{1}(x), \\
T_{2, n} & = \sum_{\abss{ \alpha} \leq 3} \frac{(\theta_{1, n} -
  \theta_{2, n})^{ \alpha_{1}} (v_{1, n} - v_{2, n})^{ \alpha_{2}}}
{\alpha_{1}! \alpha_{2}!} \frac{ \partial^{ \abss{ \alpha}}
  {\normDensity}}{\partial{ \theta}^{ \alpha_{1}} \partial{v}^{
    \alpha_{2}}} (x; \theta_{2, n}, v_{2, n}) + R_{2}(x)
\end{align*}
where $\abss{ \alpha} = \alpha_{1} + \alpha_{2}$ for $\alpha =
(\alpha_{1}, \alpha_{2})$. Here, $R_{1}(x)$ and $R_{2}(x)$ are Taylor
remainders that have the following explicit representations
\begin{align*}
R_{1}(x) & \mydefn 4 \sum \limits_{\abss{ \beta} = 4}
\frac{(\theta_{2, n} - \theta_{1, n})^{\beta_{1}}(v_{1, n} - v_{2,
    n})^{\beta_{2}}} {\beta_{1}! \beta_{2}!} \\ & \times \int
\limits_{0}^{1} (1 - t)^{3} \dfrac{\partial^{4}{\normDensity}}
       {\partial{ \theta^{ \beta_{1}}} \partial{ v^{\beta_{2}}}}
       \parenth{x; -\theta_{2, n} + t (\theta_{2, n} - \theta_{1, n}),
         v_{2, n} + t( v_{1, n} - v_{2, n})} dt, \nonumber \\ R_{2}(x)
       & \mydefn 4 \sum \limits_{\abss{ \beta} = 4} \frac{(\theta_{1,
           n} - \theta_{2, n})^{ \beta_{1}}(v_{1, n} - v_{2, n})^{
           \beta_{2}}} {\beta_{1} ! \beta_{2} !} \\ & \times \int
       \limits_{0}^{1}(1 - t)^{3} \dfrac{ \partial^{4}{\normDensity}}{
         \partial{ \theta^{\beta_{1}}} \partial{ v^{ \beta_{2}}}}
       \parenth{x; \theta_{2, n} + t ( \theta_{1, n} - \theta_{2, n}),
         v_{2, n} + t (v_{1, n} - v_{2, n})}dt.
\end{align*}
Recall from equation~\eqref{EqnAlgebra} that the univariate location-scale
Gaussian distribution
has the PDE structure of the following form
\begin{align} 
\dfrac{ \partial^{2}{ \normDensity}}{ \partial{ \theta^{2}}} (x;
\mean, \sd^2) = 2 \dfrac{ \partial { \normDensity}} { \partial{
    \sigma^{2}}} (x; \mean, \sd^2). \notag
\end{align}
Therefore, we can write the formulations of $T_{1, n}$ and 
$T_{2, n}$ as follows:
\begin{align*}
T_{1, n} & = \sum_{\abss{ \alpha} \leq 3} \frac{(\theta_{2, n} -
  \theta_{1, n})^{ \alpha_{1}} (v_{1, n} - v_{2, n})^{ \alpha_{2}}}
{2^{\alpha_{2}} \alpha_{1}! \alpha_{2}!} \frac{ \partial^{ \alpha_{1}
    + 2 \alpha_{2}} {\normDensity}}{\partial{ \theta}^{ \alpha_{1} + 2
    \alpha_{2}}} (x; - \theta_{2, n}, v_{2, n}) + R_{1}(x), \\
T_{2, n} & = \sum_{\abss{ \alpha} \leq 3} \frac{(\theta_{1, n} -
  \theta_{2, n})^{ \alpha_{1}} (v_{1, n} - v_{2, n})^{ \alpha_{2}}}
{2^{ \alpha_{2}} \alpha_{1}! \alpha_{2}!}  \frac{ \partial^{
    \alpha_{1} + 2 \alpha_{2}} {\normDensity}}{\partial{ \theta}^{
    \alpha_{1} + 2 \alpha_{2}}} (x; \theta_{2, n}, v_{2, n}) +
R_{2}(x).
\end{align*}
Via a Taylor series expansion, we find that
\begin{align*}
\frac{ \partial^{ \alpha_{1} + 2 \alpha_{2}} {\normDensity}}{\partial{
    \theta}^{ \alpha_{1} + 2 \alpha_{2}}} (x; \theta_{2, n}, v_{2, n})
= \sum_{\tau = 0}^{3 - \abss{ \alpha}} \frac{(2 \theta_{2, n})^{
    \tau}}{\tau !}  \frac{ \partial^{ \alpha_{1} + 2 \alpha_{2} +
    \tau} {\normDensity}}{\partial{ \theta}^{ \alpha_{1} + 2
    \alpha_{2} + \tau}} (x; - \theta_{2, n}, v_{2, n}) + R_{2,
  \alpha}(x)
\end{align*}
for any $\alpha = (\alpha_{1}, \alpha_{2})$ such that $1 \leq \abss{
  \alpha} \leq 3$. Here, $R_{2, \alpha}$ is Taylor remainder admitting
the following representation
\begin{align*}
R_{2, \alpha}(x) & = \sum \limits_{\tau = 4 - \abss{\alpha}}
\frac{\tau \parenth{2 \theta_{2, n}}^{\tau}}{ \tau !}  \int
\limits_{0}^{1}(1 - t)^{\tau - 1} \dfrac{ \partial^{4}{ \normDensity}}
       {\partial{ \theta^{\alpha_{1} + \tau}} \partial{ v^{
             \alpha_{2}}}} \parenth{x; -\theta_{2, n} + 2 t \theta_{2,
           n}, v_{2, n}} dt.
\end{align*}
Governed by the above results, we can rewrite $f_{\eta_{1, n}}(x) -
f_{\eta_{2, n}}(x)$ as
\begin{align*}
f_{\eta_{1, n}}(x) - f_{\eta_{2, n}}(x) = \sum_{l = 1}^{6} A_{l, n}
\dfrac{ \partial^{l} {\normDensity}}{ \partial{\theta^{ l}}} (x;
-\theta_{2, n}, v_{2, n})+ R(x)
\end{align*}
where the explicit formulations of $A_{l, n}$ and $R(x)$ are given by
\begin{align*}
A_{l, n} & \mydefn \frac{1}{2} \sum \limits_{\alpha_{1},\alpha_{2}}
\dfrac{1}{2^{\alpha_{2}}} \dfrac{(\theta_{2, n} - \theta_{1,
    n})^{\alpha_{1}}(v_{1, n} - v_{2, n})^{ \alpha_{2}}}{\alpha_{1} !
  \alpha_{2} !}  \nonumber \\
& \hspace{ 3 em} + \frac{1}{2} \sum \limits_{\alpha_{1}, \alpha_{2},
  \tau} \dfrac{1}{2^{ \alpha_{2}}} \dfrac{2^{ \tau}( \theta_{2,
    n})^{\tau} (\theta_{1, n} - \theta_{2, n})^{ \alpha_{1}} (v_{1, n}
  - v_{2, n})^{ \alpha_{2}}}{\tau !  \alpha_{1} !  \alpha_{2} !},
\nonumber \\
R(x) & \mydefn \frac{1}{2} R_{1}(x) + \frac{1}{2} R_{2}(x) +
\sum_{\abss{ \alpha} \leq 2} \dfrac{1}{2^{ \alpha_{2}}}
\dfrac{(\theta_{1, n} - \theta_{2, n})^{\alpha_{1}} (v_{1, n} - v_{2,
    n})^{\alpha_{2}}} {\alpha_{1} !  \alpha_{2} !} R_{2, \alpha}(x)
\nonumber
\end{align*}
for any $l \in [6]$ and $x \in \mathbb{R}$. Here the ranges of
$\alpha_{1}, \alpha_{2}$ in the first sum of $A_{l, n}$ satisfy
$\alpha_{1} + 2 \alpha_{2} = l$ and $1 \leq \abss{ \alpha} \leq 3$
while the ranges of $\alpha_{1}, \alpha_{2}, \tau$ in the second sum
of $A_{l, n}$ satisfy $\alpha_{1} + 2 \alpha_{2} + \tau = l$, $0 \leq
\tau \leq 3 - \abss{\alpha}$, and $1 \leq \abss{ \alpha} \leq 3$.

From the conditions that $\theta_{2, n} = 2 \theta_{1, n}$ and $v_{1,
  n} - v_{2, n} = 3 \parenth{ \theta_{1, n}}^{2}$, we can check that
$A_{l, n} = 0$ for all $1 \leq l \leq 3$.  Additionally, we also have
\begin{align*}
\max \{\abss{A_{4, n}}, \abss{A_{5, n}}, \abss{A_{6, n}}\} \precsim
\abss{\theta_{1, n}}^{4}.
\end{align*}
Given the above results, we claim that
\begin{align}
\label{eqn:upper_hellinger}
h \parenth{f_{\eta_{1, n}}, f_{\eta_{2, n}}} \precsim \abss{
  \theta_{1, n}}^{8}.
\end{align}
Assume that the claim~\eqref{eqn:upper_hellinger} is given. From the
formulations of sequences $\eta_{1, n}$ and $\eta_{2, n}$, we can
verify that
\begin{align*}
\parenth{ \parenth{\abss{ \theta_{1, n}} - \abss{\theta_{2, n}}}^{2} +
  \abss{ v_{1, n} - v_{2, n}}}^r \asymp \abss{ \theta_{1, n}}^{2 r}.
\end{align*}
Since $1 \leq r < 4$ and $\theta_{1, n} \to 0$ as $n \to \infty$, the
above results lead to
\begin{align*}
h \parenth{f_{\eta_{1, n}}, f_{\eta_{2, n}}} \big/ \parenth{
  \parenth{\abss{ \theta_{1, n}} - \abss{\theta_{2, n}}}^{2} + \abss{
    v_{1, n} - v_{2, n}}}^r \precsim \abss{ \theta_{1, n}}^{8 - 2 r}
\to 0.
\end{align*}
As a consequence, we achieve the conclusion of the
claim~\eqref{eqn:minimax_bound_key}.

%%%%%%%%%%%%%%%%%%%%%%%%%%%%%%%%%%%%%%%%%%%%%%%%%%%%%%%%%%%%%%%%%%%%%

\subsubsection{Proof of claim~\eqref{eqn:upper_hellinger}} % (fold)
\label{ssub:proof_of_claim_eqn:upper_hellinger}

The definition of Hellinger distance leads to the following equations
\begin{align}
2 h^2 \parenth{f_{\eta_{1, n}}, f_{\eta_{2, n}}} & = \int
\dfrac{\parenth{ f_{\eta_{1, n}}(x) - f_{\eta_{2, n}}(x)}^2}
      {\parenth{ \sqrt{ f_{\eta_{1, n}}(x)} + \sqrt{ f_{\eta_{2,
                n}}(x)}}^2} dx \nonumber \\
  & = \int \dfrac{(\sum \limits_{l = 4}^{6} A_{l, n} \dfrac{
            \partial^{l}{ \normDensity}} {\partial{ \theta^{l}}}(x; -
          \theta_{2, n},v_{2, n}) + R(x))^2}{\parenth{ \sqrt{
            f_{\eta_{1, n}}(x)} + \sqrt{ f_{\eta_{2, n}}(x)}}^2} dx
      \nonumber \\
  \label{eqn:minimax_first}      
& \precsim \int \dfrac{ \sum_{l = 4}^{6} \parenth{ A_{l, n}}^2
          \parenth{ \dfrac{ \partial^{l}{ \normDensity}} {\partial{
                \theta^{l}}}(x; - \theta_{2, n},v_{2, n})}^2 +
          R^2(x)}{\parenth{ \sqrt{ f_{\eta_{1, n}}(x)} + \sqrt{
              f_{\eta_{2, n}}(x)}}^2} dx,
\end{align}
where the last inequality is due to Cauchy-Schwarz's
inequality. According to the structure of location-scale Gaussian
density, the following inequalities hold
\begin{align}
\int \dfrac{ \parenth{ \dfrac{ \partial^{l}{ \normDensity}} {\partial{
        \theta^{ l}}}(x; - \theta_{2, n}, v_{2, n})}^{2}}{\parenth{
    \sqrt{ f_{\eta_{1, n}}(x)} + \sqrt{ f_{\eta_{2, n}}(x)}}^2} dx
\precsim \int \dfrac{\parenth{ \dfrac{ \partial^{l}{ \normDensity}}
    {\partial{ \theta^{ l}}}(x; - \theta_{2, n}, v_{2, n})}^{2}} {
  \normDensity(x; - \theta_{2, n}, v_{2, n})} dx <
\infty \label{eqn:minimax_second}
\end{align}
for $4 \leq l \leq 6$. Note that, for any $\beta = (\beta_{1},
\beta_{2})$ such that $\abss{ \beta} = 4$, we have
\begin{align*}
\abss{\theta_{2, n} - \theta_{1, n}}^{\beta_{1}} \abss{ v_{1, n} -
  v_{2, n}}^{\beta_{2}} \asymp \abss{ \theta_{1, n}}^{4 + \beta_{2}}
\precsim \abss{ \theta_{1, n}}^4.
\end{align*}
With the above bounds, an application of Cauchy-Schwarz's inequality
leads to
\begin{align*}
& \int \frac{ R_{1}^2( x)}{\parenth{ \sqrt{ f_{\eta_{1, n}}(x)} +
      \sqrt{ f_{\eta_{2, n}}(x)}}^2} dx \\ & \precsim \abss{
    \theta_{1, n}}^8 \sum_{\abss{ \beta} = 4} \int \frac{\sup
    \limits_{t \in [0, 1]} \parenth{ \dfrac{\partial^{4}
        {\normDensity}} {\partial{ \theta^{ \beta_{1}}} \partial{
          v^{\beta_{2}}}} \parenth{x; -\theta_{2, n} + t (\theta_{2,
          n} - \theta_{1, n}), v_{2, n} + t( v_{1, n} - v_{2,
          n})}}^2}{ \normDensity(x; - \theta_{2, n}, v_{2, n})} dx
  \precsim \abss{ \theta_{1, n}}^8.
\end{align*}
With a similar argument, we also obtain that
\begin{align*}
\int \frac{ R_{2}^2(x)}{\parenth{ \sqrt{ f_{\eta_{1, n}}(x)} + \sqrt{
      f_{\eta_{2, n}}(x)}}^2} dx & \precsim \abss{ \theta_{1, n}}^8,
\ \
\max_{1 \leq \abss{ \alpha} \leq 4} \int \frac{ R_{2, \alpha}^2(x)}
    {\parenth{ \sqrt{ f_{\eta_{1, n}}(x)} + \sqrt{ f_{\eta_{2,
              n}}(x)}}^2} dx & \precsim \abss{ \theta_{1, n}}^8.
\end{align*}
Governed by the above bounds, another application of Cauchy-Schwarz's
inequality implies that
\begin{align}
\int \frac{ R^2( x)}{\parenth{ \sqrt{ f_{\eta_{1, n}}(x)} + \sqrt{
      f_{\eta_{2, n}}(x)}}^2} dx \precsim \abss{ \theta_{1,
    n}}^{8}. \label{eqn:minimax_third}
\end{align}
Combining the results from equations~\eqref{eqn:minimax_first},
~\eqref{eqn:minimax_second}, and~\eqref{eqn:minimax_third}, we achieve
the conclusion of the claim~\eqref{eqn:upper_hellinger}.

%%%%%%%%%%%%%%%%%%%%%%%%%%%%%%%%%%%%%%%%%%%%%%%%%%%%%%%%%%%%%%%%%

%%%%%%%%%%%%%%%%%%%%%%%%%%%%%%%%%%%%%%%%%%%%%%%%%%%%%%%%%%%%%%%%%
\section{Proofs of auxiliary results} % (fold)
\label{sec:proofs_of_other_auxiliary_results}
In this appendix, we collect the proofs of several auxiliary results
stated throughout the paper.

\subsection{Proof of Corollary~\ref{cor:convergence_rate_sample_EM_univariate_balanced}} 
\label{sub:proof_of_prop_EM_univariate_balanced}

In order to ease the presentation, we only provide the proof sketch
for the localization argument with this corollary. The detail proof
argument for the corollary can be argued in similar fashion as that of
Theorem~\ref{thm:em_univariate}.  In particular, we
consider the iterations $t$ such that $\theta_n^t \in [n^{- a_ \ell},
  n^{- a_r} ]$ where $a_\ell > a_r$. For all such iterations with
$\theta_n^t$, invoking Lemma~\ref{lemma:correct_pop_contraction_and_radem_three}, we find that
\begin{align*}
\abss{ \Mtil_{\obs,1}( \theta_n^t)} \lesssim \underbrace{ (1 - n^{- 6{
      a_\ell}})}_{ \rdefn \gamma_{ a_ \ell}} \abss{ \theta_n^t} \quad
\text{ and} \quad \abss{ \samm( \theta_n^t) - \Mtil_{\obs,1}( \theta_n^t)}
\lesssim n^{- a_r}/ \sqrt{n}.
\end{align*}
Therefore, we obtain that
\begin{align*}
  \abss{ \theta_{n}^{t + T}} & \leq \abss{ \Mtil_{\obs,1}( \theta_{n}^{
      t + T - 1})} + \abss{ \Mtil_{\obs,1}( \theta_{n}^{t + T - 1}) -
    \samm( \theta_{n}^{ t + T - 1})} \leq \gamma_{ a_\ell}
  \theta_{n}^{ t + T - 1} + n^{- a_r}/ \sqrt{n}.
\end{align*}
Unfolding the above inequality $T$ times, we find that
\begin{align*}
\abss{ \theta_{n}^{t + T}} \leq \gamma_{ a_\ell}^2 ( \theta_{n}^{t + T
  - 2}) + n^{- a_r}/ \sqrt{n} (1 + \gamma_m) & \leq \gamma_{ a_\ell}^T
\theta_{n}^{t} + (1 + \gamma_{ a_\ell} + \ldots + \gamma_{ a_\ell}^{T
  - 1}) n^{- a_r}/ \sqrt{n} \\ & \leq e^{-T n^{- 6 {a_\ell}} } n^{-
  a_r} + \frac{1}{1 - \gamma_{ a_\ell}} \cdot n^{- a_r}/ \sqrt{n}.
\end{align*}
As $T$ is sufficiently large such that the second term is the dominant
term, we find that that
\begin{align*}
\abss{ \theta_{n}^{t + T} } \lesssim \frac{ 1}{1 - \gamma_{ a_\ell}}
\cdot n^{- a_r}/ \sqrt{n} = n^{6 a_\ell - a_r - 1/2}.
\end{align*}
Setting the RHS equal to $n^{- a_\ell}$, we obtain the recursion that
\begin{align}
\label{eq:a_l_r_univariate}
  a_\ell = \frac{a_r}{ 7} + \frac{ 1}{ 14}.
\end{align}
Solving for the limit $a_\ell = a_r = a_\star$ yields that 
$a_\star = 1/ 12$. It suggests that we eventually have 
$\theta_n^t \to \ball(0, n^{-\myfrac{1}{12}})$. As a consequence, 
we achieve the conclusion of the corollary.

%%%%%%%%%%%%%%%%%%%%%%%%%%%%%%%%%%%%%%%%%%%%%%%%%%%%%%%%%%%%%%%%%%%%%%%%%
\subsection{Proof of Lemma~\ref{lemma:non_expansive}} 
\label{sub:proof_of_lemma_lemma:non_expansive}
Without loss of generality, we can assume that $\abss{\theta} \in 
[\staterr^{- a_{\ind + 1}}, \staterr^{- a_\ind}]$.
Conditional on the event $\event$,
we have that
\begin{align*}
  \abss{ \corop( \theta)} \leq (1 - \staterr^{ -6 a_{\ind + 1}}/ 5) 
  \abss{ \theta}
  \quad \text{ and} \quad
  \abss{ \samm( \theta) - \corop( \theta) }
  \leq c_2 \staterr^{- 3 a_\ind} \staterr^{-\myfrac{1}{2}}.
\end{align*}
As a result, we have
\begin{align*}
  \abss{ \samm( \theta)}
   \leq \abss{ \samm( \theta) - \corop( \theta)} 
   + \abss{ \corop( \theta)}
  & \leq ( 1 - \staterr^{ - 6 a_{\ind + 1}}/ 5 ) \abss{ \theta} 
  + c_2 \staterr^{-\myfrac{1}{2}} \staterr^{- 3 a_{ \ind}} \\
  & \leq ( 1 - \staterr^{ - 6 a_{\ind + 1}}/ 5 + c_2 \staterr^{-\myfrac{1}{2}} \staterr^{- 2 a_\ind}) 
  \staterr^{- a_\ind} \\
  & \leq \staterr^{- a_\ind}.
\end{align*}
Here, to establish the last inequality, we have used the following observation:
for \mbox{$\omega = n/c_{n, \delta}$} and that 
$n \geq (c')^{1 / \smallthreshold} c_{n, \delta}$, we have
\begin{align*}
  5 c_{2} \staterr^{ 6 a_{\ind + 1} - 2 a_{\ind} - 1/ 2} 
  \leq 5 c_{2} \staterr^{ 4 a_{\ind} - 1/ 2}
  \leq c' \staterr^{4 a_{\lstar} - 1/ 2} 
  \leq c' \staterr^{- 4 \smallthreshold} 
  \leq 1/ (c')^3 \leq 1,
\end{align*}
which leads to 
$ - \staterr^{ - 6 a_{\ind + 1}}/ 5 + c_2 \staterr^{-\myfrac{1}{2}} 
\staterr^{- 2 a_\ind} \leq 0$.
As a consequence, we achieve the conclusion of the lemma.%%%%%%%%%%%%%%%%%%%%%%%%%%%%%%%%%%%%%%%%%%%%%%%%%%%%%%%%%%

\subsection{Proof of Lemma~\ref{lemma:d_dim_operator}} % (fold)
\label{sec:proof_of_lemma_lemma:d_dim_operator}
The proof of the perturbation bound~\eqref{eq:define_e_r_multi}
is a standard extension of $\dims = 1$ case presented above
in Section~\ref{ssub:proof_of_perturbation_bound_for_Mtil}, and thereby
is omitted.

We now present the proof of the contraction bound~\eqref{eq:pop_em_multivariate_contraction_bound},
which has several similarities with the proofs of 
bounds~\eqref{eq:lemma_pop_event} and \eqref{eq:lemma_cor_pop} from Lemma~\ref{lemma:correct_pop_contraction_and_radem_three}. 
In order to simplify notation, we use the shorthand $Z_{n, d} \defn
\frac{1}{\obs d} \sum_{j = 1}^\obs \enorm{ X_j}^2$.  Recalling the
definition~\eqref{eq:pop_like_em_operator} of operator
$\PseudoND( \theta)$, we have
\begin{align}
    \label{eq:pseudo_pop_thm_norm}
\enorm{ \PseudoND( \theta)} & = \left \| \Exs_{Y \sim \Ncal(0,
  1)} \brackets{Y \tanh \parenth{\frac{Y^{\top} \theta}{Z_{\obs, d} -
      \enorm{ \theta}^2/ d}}} \right \|_2.
\end{align}
We can find an orthonormal matrix $R$ such that \mbox{$R \theta =
  \enorm{ \theta} \basis_{1}$,} where $\basis_{1}$ is the first
canonical basis in $\Rspace^{d}$. Define the random vector $V = R Y$.
Since $Y \sim \Ncal(0, I_{d})$, we have that $V \sim \Ncal(0, I_{d})$.
On performing the change of variables $Y = R^{\top} V$, we find that
\begin{align*}
\enorm{ \Exs_{Y} \brackets{ Y \tanh \parenth{\frac{ Y^{ \top} \theta}{
        Z_{n, d}- \enorm{ \theta}^2/ d}}}} & = \enorm{\Exs_{V}
  \brackets{ R^{ \top} V \tanh\parenth{ \frac{ \enorm{ \theta} V_{1}}{
        Z_{n, d} - \enorm{ \theta}^2/ d}}}} \\
& = \abss{ \Exs_{V_{1}} \brackets{ V_{1} \tanh \parenth{ \frac{
        \enorm{ \theta} V_{1}}{Z_{n, d} - \enorm{ \theta}^2/ d}}}}
\end{align*}
where the final equality follows from the fact that
\begin{align*}
\Exs[ R^\top V f(V_1)] = R^\top \Exs [V f(V_1)] = R^\top (\Exs[ V_1
  f(V_1)], 0, \ldots, 0)^\top.
\end{align*}
Furthermore, the orthogonality of the matrix $R$ implies that $\enorm{
  \Exs[ R^\top V f(V_1)]}^2 = \abss{ \Exs[ V_1 f(V_1)]}^2$.
 
In order to simplify the notation, we define the scalars $a, b$ and
the event $\event_{\alpha, d}$ as follows:
\begin{subequations}
\begin{align}
\label{eq:defn_a_b_event}
  a \defn 1 - (n d)^{- \alpha}, \quad b \defn 1 + (n d)^{- \alpha},
  \quad\text{and}\quad \event_{\alpha, d} = \braces{ \vert Z_{n, d} -
    1 \vert \leq {(\obs d)^{- \alpha}}},
\end{align}
where $\alpha$ is a suitable scalar to be specified later.  Note that
standard chi-squared tail bounds guarantee that
\begin{align}
\label{eq:event_prob_bound}
 \Prob[ \event_{ \alpha, d}] \geq 1 - 2 e^{- {\dims^{2 \alpha} \obs^{1
       - 2\alpha}}/ {8}}.
\end{align}
\end{subequations}
Now conditional on the event $\event_{\alpha, d}$, we have
\begin{align*}
\enorm{ \PseudoND( \theta)} & \leq \abss{ \Exs_{ V_{1}}
  \brackets{ V_{1} \tanh \parenth{ \frac{ \enorm{\theta} V_{1}}{a -
        \enorm{ \theta}^2/ d}}}} = \enorm{ \theta} \underbrace{ \Exs_{
    V_{1}} \brackets{ \frac{ V_{1}}{ \enorm{ \theta}} \tanh \parenth{
      \frac{ \enorm{ \theta} V_1}{ a - \enorm{ \theta} ^2/ d}}}}_{
  \rdefn \locacontramutiup( \theta)}, \quad \text{and},\nonumber
\\ \enorm{ \PseudoND( \theta)} & \geq \abss{ \Exs_{ V_{1}}
  \brackets{ V_{1} \tanh \parenth{ \frac{ \enorm{ \theta} V_{1}}{b -
        \enorm{ \theta}^2/ d}}}} = \enorm{ \theta} \underbrace{ \Exs_{
    V_{1}} \brackets{ \frac{ V_{1}}{ \enorm{ \theta}} \tanh \parenth{
      \frac{ \enorm{ \theta} V_1}{b - \enorm{ \theta}^2/ d}}}}_{
  \rdefn \locacontramutilo ( \theta)},
\end{align*}
where the above inequalities follow from the fact that for any fixed
$y, \theta \in \realdim$, the function $w \mapsto y \tanh( y \enorm{
  \theta}/ (w - \enorm{\theta}^2/d))$ is non-increasing in $w$ for $w>
\enorm{ \theta}^2/ d$.

Substituting $\alpha = 1/2 - 2 \smallthreshold$ in the
bound~\eqref{eq:event_prob_bound} and invoking the large sample size
assumption in the theorem statement, we obtain that
$\Prob[\event_{\alpha, d}] \geq 1 - \delta$.  Putting these
observations together, it remains to prove that
\begin{align}
\label{eq:rho_bound_multi}
 \locacontramutilo( \theta) \geq \parenth{1 - \dfrac{ 3 \enorm{ \theta}^
 {2}}{ 4}}\enorm{ \theta}^2, \quad \text{ and} \quad
 \locacontramutiup ( \theta) \leq \parenth{1 - \parenth{1 -
     \frac{1}{d}} \frac{ \enorm{ \theta}^2}{ 4}} \enorm{ \theta}^2,
\end{align}
for all $5 (d/ n)^{- 1/ 4 + \smallthreshold} \leq \enorm{ \theta}^2
\leq (d - 1)/(6 d - 1)$ conditional on the event $\event_{\alpha, d}$
for $\alpha = 1/2-6\smallthreshold$ to obtain the conclusion of the
theorem.

The proof of the claims in equation~\eqref{eq:rho_bound_multi} relies
on the following bounds on the hyperbolic function $\tanh(x)$.
For any $x \in \Rspace$, the following bounds hold:
\begin{align}
\label{eq:lemma_six_bounds}
\text{(Upper bound)} \quad  x^2 - \frac{x^4} {3} +
\frac{2 x^6}{ 15}\  \geq \ \  x \tanh( x) \ \   \geq\  x^2 - \frac{x^4}
{3} \quad\text{(Upper bound)}.
\end{align}
We omit the proof of
these bounds, as it is very
similar to that of similar results stated and proven later in
Lemma~\ref{lemma:inequality_tanh_sech_function}.  We now turn to
proving the bounds stated in equation~\eqref{eq:rho_bound_multi}
one-by-one.

%%%%%%%%%%%%%%%%%%%%%%%%%%%%%%%%%%%%%%%%%%%%%%%%%%%%%%%%%%%%%%%%%%%%%%

\paragraph{Bounding $\locacontramutiup( \theta)$:}

Applying the upper bound~\eqref{eq:lemma_six_bounds} for $x\tanh(x)$, we obtain that
\begin{align*} 
\locacontramutiup( \theta) & \leq \frac{a - \enorm{ \theta}^2/ d}{
  \enorm{ \theta}^2} \biggr( \frac{ \enorm{ \theta}^2}{(a - \enorm{
    \theta}^2/ d)^2} \Exs \brackets{ V_{1}^2} - \frac{ \enorm{
    \theta}^4}{3 (a - \enorm{ \theta}^2/ d)^4} \Exs \brackets{
  V_{1}^4} + \frac{ 2 \enorm{ \theta}^6}{15 (a - \enorm{ \theta}^2/
  d)^6} \Exs \brackets{ V_{1}^6}\biggr).
\end{align*}
Substituting $\Exs \brackets{V_{1}^{2k}} = (2k-1)!!$ for $k = 1, 2, 3$
in the RHS above, we find that
\begin{align} 
\locacontramutiup( \theta) \leq \frac{1}{a - \enorm{ \theta}^2/ d} -
\frac{ \enorm{ \theta}^2}{(a - \enorm{ \theta}^2/ d)^3} + \frac{ 2
  \enorm{ \theta}^4}{(a - \enorm{ \theta}^2/ d)^5}.
  \label{eq:upper_rho_final}
\end{align}
The condition $\enorm{ \theta}^2 + (n d)^{ - \alpha} \leq \frac{d -
  1}{6 d - 4} < 1/ 6$ implies the following bounds:
\begin{align*}
& \frac{1}{{1 - (n d)^{ - \alpha} - \enorm{ \theta}^2/ d}} \leq 1 +
  \parenth{(n d)^{ - \alpha} + \enorm{ \theta} ^2/ d } + 3/ 2 \cdot
  \parenth{(n d)^{- \alpha} + \enorm{\theta}^2/d }^{2}, \\
& \frac{1}{({1 - (n d)^{- \alpha} - \enorm{ \theta}^2/ d})^3} \geq 1 +
  3 \parenth{(n d)^{- \alpha} + \enorm{ \theta} ^2/ d }, \\ &
  \frac{1}{({1 - (n d)^{- \alpha} - \enorm{ \theta}^2/ d})^5} \leq
  3/2.
\end{align*}
Substituting the definitions~\eqref{eq:defn_a_b_event} of $a$ and $b$
and plugging the previous three bounds on the RHS of the
inequality~\eqref{eq:upper_rho_final} yields that
\begin{align*}
\locacontramutiup( \theta) & \leq 1 + \frac{ \enorm{ \theta}^2}{ d} +
\frac{ 3 \enorm{ \theta}^4}{2 d^2} - \enorm{ \theta} ^2 \parenth{1 +
  \frac{ 3 \enorm{ \theta}^2}{ d}} + 3 \enorm{ \theta}^4 + \frac{ 11}{
  2}(n d)^{- \alpha} \\
& \leq 1 - \parenth{1 - \frac{ 1}{ d}} \enorm{ \theta}^2 + \parenth{3
  - \frac{ 2}{ d}} \enorm{ \theta}^4 + \frac{ 11}{ 2} (n d)^{- \alpha}
\\ & \leq 1 - \parenth{1 - \frac{ 1}{ d}} \frac{ \enorm{ \theta}^2}{
  4}
\end{align*}
where the last step follows from the following observations that
\begin{align}
  \label{eq:theta_lower_bound_one}  
(3 - 2/ d) \enorm{ \theta}^4 & \leq (1 - 1/ d) \enorm{ \theta}^2/ 2,
  \quad \text{for all } \enorm{ \theta} \leq (d - 1)/ (6 d - 4),
  \\
  \label{eq:theta_lower_bound_two}
  11(n d)^{- \alpha}/2 & \leq (1 - 1/ d) \enorm{ \theta}^{2}/ 4, \quad
  \text{for all } \enorm{ \theta} \geq 5 ( d/ n)^{- 1/ 4 +
    \smallthreshold} \ \text{when } \alpha = 1/ 2 - 2 \smallthreshold.
\end{align}
Therefore, the claim with an upper bound of $\locacontramutiup(
\theta)$ now follows.

%%%%%%%%%%%%%%%%%%%%%%%%%%%%%%%%%%%%%%%%%%%%%%%%%%%%%%%%%%%%%%%%%%%%%%%%%%%

\paragraph{Bounding $\locacontramutilo( \theta)$:}

Using the lower bound~\eqref{eq:lemma_six_bounds} for $x\tanh(x)$, 
we find that
\begin{align}
\locacontramutilo( \theta) & \geq \frac{b - \enorm{ \theta}^2/ d}{
  \enorm{ \theta} ^2} \biggparenth{ \frac{ \enorm{ \theta}^2}{({b -
      \enorm{ \theta}^2/ d})^2} \Exs \brackets{ V_{1}^2} - \frac{
    \enorm{ \theta}^4} {3( b - \enorm{ \theta}^2/ d)^4} \Exs
  \brackets{ V_{1} ^4}} \\
\label{eq:small_rho_bound}
& = \frac {1} {b - \enorm{ \theta}^2/ d} - \frac{ \enorm{ \theta}^2}
      {(b - \enorm{ \theta}^2/ d)^3}.
\end{align}
The condition $\enorm{ \theta} - (n d)^{- \alpha} \geq 0$ leads to
\begin{align*}
\frac{1} {{1 + (n d)^{- \alpha} - \enorm{ \theta}^2/ d}} & \geq 1 +
\parenth{ \enorm{ \theta}^2/ d -(n d)^{- \alpha}} + \parenth{ \enorm{
    \theta}^2/ d - (n d)^{- \alpha}}^{2}, \\
\frac{1}{({1 + (n d)^{- \alpha} - \enorm{ \theta}^2/ d})^3} & \leq 1 +
4 \parenth{ \enorm{ \theta}^2/ d -(n d)^{- \alpha}}.
\end{align*}
Applying these inequalities to the bound~\eqref{eq:small_rho_bound},
we obtain that
\begin{align*} 
\locacontramutilo(\theta) & \geq 1 + \frac{ \enorm{ \theta}^2} { d} +
\frac{ \enorm{ \theta}^4}{ d^2} - \enorm{ \theta}^2 \parenth{1 +
  \frac{ 4 \enorm{ \theta}^2} { d}} - 2(n d)^{- \alpha} \\
& \stackrel{ (i)}{ \geq} 1 - \enorm{ \theta}^2 \parenth{1 - \frac{ 1}{
    d}} - \frac{ \enorm{ \theta}^2} {6} \parenth{ \frac{ 4}{ d} -
  \frac{ 1}{ d^2}} - \frac{ \enorm{ \theta}^{ 2}(1 - 1/ d)}{ 11} \\
& \geq 1 - \dfrac{3 \enorm{ \theta}^{ 2}}{ 4}
\end{align*}
where step~(i) in the above inequalities follows from the
observations~\eqref{eq:theta_lower_bound_one}-\eqref{eq:theta_lower_bound_two}
above.  The lower bound~\eqref{eq:rho_bound_multi} for
$\locacontramutilo(\theta)$ now follows.
% subsection proof_of_lemma_lemma:d_dim_operator (end)
%%%%%%%%%%%%%%%%%%%%%%%%%%%%%%%%%%%%%%%%%%%%%%%%%%%%%%%%%%%%%%%%%

\subsection{Proof of Lemma~\ref{lemma:inequality_tanh_sech_function}} 
\label{sub:proof_of_lemma_univariate_balanced}

The proof of this lemma relies on an evaluation of 
coefficients with $x^{2 k}$ as $k \geq 1$. In particular, we divide 
the proof of the lemma into two key parts:
\paragraph{Upper bound:} From the definition of hyperbolic function $
\tanh( x)$, it is 
sufficient to demonstrate that
\begin{align*}
x \parenth{ \exp( x) - \exp(- x)} 
    \leq \parenth{ x^2 - \frac{ x^4}{ 3} + \frac{ 2 x^6}{ 15} 
    - \frac{17 x^8}{ 315} + \frac{ 62 x^{ 10}}{ 2835}}
    \parenth{ \exp( x) + \exp(- x)}.
\end{align*}
Invoking the Taylor series of $\exp( x)$ and $\exp (- x)$, the above
inequality is equivalent to
\begin{align*}
\sum \limits_{k = 0}^{ \infty} \frac{2 x^{2 k + 2}}{(2 k + 1) !} 
    \leq \parenth{ x^2 - \frac{ x^4}{ 3} + 
    \frac{2 x^6}{ 15} - \frac{ 17 x^8}{ 315} 
    + \frac{ 62 x^{10}}{ 2835}} \parenth{ \sum \limits_{k = 0}
    ^{ \infty} \dfrac{ 2 x^{ 2 k}}{(2 k) !}}. 
\end{align*}
Our approach to solve the above inequality is to show that the 
coefficients of $x^{2 k}$ in the LHS is smaller than that of $x^{2 k}$ in 
the RHS for all $k \geq 1$. 
In fact, when $1 \leq k \leq 3$, we can quickly check that the previous 
observation holds. 
For $k \geq 4$, it suffices to validate that
\begin{align*}
    \frac{ 2}{ ( 2 k) !} - \frac{ 2}{ 3 (2 k - 2) !} 
    + \frac{ 4}{ 15 ( 2 k - 4) !} - \frac{ 34}{ 315 (2 k - 6) !}
    + \frac{ 124}{ 2835 ( 2 k - 8) !} - \frac{ 2}{( 2 k + 1) !} 
    \geq 0. \nonumber
\end{align*}
Direct computation with the above inequality leads to
\begin{align*}
    (k - 1) (k - 2) (k - 3) (k - 4)
    (496 k^4 - 1736 k^3 + 1430 k^2 + 446 k - 381) 
    \geq 0
\end{align*}
for all $k \geq 4$, which is always true. As a consequence, we achieve 
the conclusion with the upper bound of the lemma.
\paragraph{Lower bound:} For the lower bound of the lemma, it is 
equivalent to prove that
\begin{align*}
\sum \limits_{k = 0}^{ \infty} \frac{2 x^{ 2 k + 2}}{( 2 k + 1) !} 
    \geq \parenth{ x^2 - \frac{ x^4}{ 3} + \frac{ 2 x^6}{ 15} 
    - \frac{ 17 x^8}{ 315}} \parenth{ \sum \limits_{ k = 0}
    ^{ \infty} \dfrac{ 2 x^{ 2 k}}{ (2 k) !}}. 
\end{align*}
Similar to the proof technique with the upper bound, we only need to 
verify that
\begin{align*}
    \frac{ 2}{ ( 2 k) !} - \frac{ 2}{ 3 (2 k - 2) !} 
    + \frac{ 4}{ 15 ( 2 k - 4) !} - \frac{ 34}{ 315 ( 2 k - 6) !}
     - \frac{2}{(2 k + 1)!} \leq 0 \nonumber
\end{align*}
for any $k \geq 3$. The above inequality is identical to
\begin{align}
    (k - 1) (k - 2) (k - 3) (4352 k^3 - 4352 k^2 - 512 k + 1472)
     \geq 0 \nonumber
\end{align}
for all $k \geq 3$, which always holds. Therefore, we obtain the 
conclusion with the lower bound of the lemma.

%%%%%%%%%%%%%%%%%%%%%%%%%%%%%%%%%%%%%%%%%%%%%%%%%%%%%%%%%%%%%%%%%%%%%%%%%%%%%%%%%%%%%%%%%%%

\subsection{Proof of Lemma \ref{lemma:concentration_high_order}} 
\label{sub:proof_of_lemma_concentration_high_order}

The proof of this lemma is based on appropriate truncation argument.
More concretely, given any positive scalar 
$\truncate$, and the random variable $X \sim \Ncal(0,1)$, consider the pair of
truncated random variables $(Y,  Z)$ defined by: 
\begin{align}
\label{eqn:truncation}
  Y \defn  X^{2 k} \mathbb{I}
    _{\abss{ X} \leq \truncate} \quad  \text{and} \quad 
    Z \defn  X^{2 k} \mathbb{I}_{ \abss{X} \geq \truncate}.
\end{align}
With the above notation in place, for $\obs$ i.i.d. samples
$X_1, \ldots , X_n$ from $\Ncal(0,1)$, we have
\begin{align*}
\frac{ 1}{ n} \sum_{i = 1}^{n} X_{i}^{2 k} = 
\frac{1}{n} \sum_{i = 1}^{n} Y_{i} + \frac{1}{n} \sum_{i = 1}
    ^{n} Z_{i}   \mydefn S_{Y, n} + S_{Z, n}.
\end{align*}
where $S_{Y, n}$ and $S_{Z, n}$, denote the averages of the random
variables $Y_i's$ and $Z_i's$ respectively.  Observe that $\abss{
  Y_{i}} \leq \truncate^{2 k}$ for all $i \in [n]$; consequently, by
standard sub-Gaussian concentration of bounded random variables, we
have
\begin{align}
\Prob \parenth{ |S_{Y, n} - \Exs \brackets{Y} | \geq   t_{1} } 
    \leq 2 \exp \parenth{ - \frac{n t_{1}^2}{2 \truncate^{4 k}}}.
    \label{eqn:subgaussian_X2k}
\end{align}
Next, applying Markov's inequality with the non-negative random
variable $S_{Z, n}$, we find that
\begin{align}
\Prob \parenth{ S_{Z, n} \geq t_{2} } \leq \frac{\Exs \brackets{ S_{Z,
      n} }}{ t_{2}} = \frac{ \Exs \brackets{Z} }{ t_{2}}.
    \label{eqn:Markov_Zk}
\end{align}  
By definition of the truncated random variable $Y$, we have
$\Exs [Y]  \leq \Exs [X^{2k}]$; moreover, an application of Holder's inequality 
to $\Exs \brackets{Z}$ yields 
\begin{align*}
\Exs \brackets{Z} 
    = \Exs \parenth{ X^{2 k} \mathbb{I}_{\abss{X} \geq \truncate}} 
    \leq \sqrt{\Exs \brackets{X^{4 k}} } \sqrt{\Prob \parenth{ \abss{X} 
    \geq \truncate}} 
    \leq \sqrt{2 \Exs \brackets{ X^{4 k}} } \exp(- \truncate^2/ 4).
\end{align*}
Combining the bounds on $\Exs [Y]$ and $\Exs[Z]$ with the
inequalities~\eqref{eqn:subgaussian_X2k} and~\eqref{eqn:Markov_Zk} 
we deduce that
\begin{subequations} 
\begin{align}
& \frac{ \sum_{i = 1}^{n} X_{i}^{2 k}}{ n} 
    \leq \Exs \brackets{ Y} + t_{1} + t_{2} 
    \leq \Exs \brackets{ X^{2 k}} + t_{1} + t_{2}, \ \ 
    \text{and}, 
    \label{eq:xi_2_upper_bound}
     \\
& \frac{\sum_{i = 1}^{n} X_{i}^{2 k}}{ n} 
    \geq \Exs \brackets{ X^{2 k}} - t_{1} - t_{2} \sqrt{2 \Exs 
    \brackets{ X^{4 k}} }
    \exp(- \truncate^2/ 4) 
    \label{eq:xi_2_lower_bound} 
\end{align}
\end{subequations}
with probability at least $1 - \exp \parenth{ - \frac{n t_{1}^2}{2 
\truncate^{4 k}}} -  \sqrt{2 \Exs \brackets{ X^{4 k}} }\exp(- 
\truncate^2/ 4)$. Finally, given any $\delta > 0$, choose the scalars
$\truncate, t_1, t_2$ as follows:
\begin{align*}
  \truncate = 2 \sqrt{ \log \parenth{ \frac{2 
    \sqrt{2 \obs \Exs \brackets{  X^{4k}} }}{\delta}}},
    \quad  t_1 = \truncate^2 \sqrt{ \frac{1}{n} 
    \log \left(\frac{2}{\delta}\right)}    \quad \text{and}
    \quad t_2 = \frac{1}{\sqrt{\obs}}.
\end{align*}
Substituting the choice of $t_{1}, t_{2}$ and $\truncate$, in 
bounds~\eqref{eq:xi_2_upper_bound} and~\eqref{eq:xi_2_lower_bound}
we conclude that with probability at least $1 - \delta$ 
\begin{align*}
\abss{ \frac{\sum_{i = 1}^{n} X_{i}^{2 k}}{ n} 
    - \Exs \brackets{ X^{2 k}} }
    \leq \frac{ C_{k} \log^{ k}(n/ \delta)}{ \sqrt{ n}},
\end{align*}
where $C_{k}$ is a universal constant that depends only on $k$. This 
completes the proof of Lemma~\ref{lemma:concentration_high_order}.
%%%%%%%%%%%%%%%%%%%%%%%%%%%%%%%%%%%%%%%%%%%%%%%%%%%%%%%%%%%%%%%%%

\subsection{Proof of one step bound for population EM} 
\label{sub:proof_of_lemma_lemma:one_step_update}
We now describe a special one-step contraction property of the population
operator.
\begin{lemma}
  \label{lemma:one_step_update}
  For any vector $\theta^0$ such that $\Vert \theta^0\Vert \leq
  \sqrt{d}$, we have $\Vert\PseudoND(\theta^0)\Vert \leq \sqrt{2/\pi}$
  with probability at least $1-\delta$.
\end{lemma}
The proof of this lemma is a straightforward application of the proof
argument in Lemma~\ref{lemma:d_dim_operator} in
Appendix~\ref{sec:proof_of_lemma_lemma:d_dim_operator}. In order to
simplify notations, we use the shorthand $Z_{n, d} = \sum_{j = 1}^\obs
\enorm{ X_j}^2/( \obs d)$.  Recalling the
definition~\eqref{eq:pop_like_em_operator} of operator
$\PseudoND$, we have
\begin{align*}
\enorm{ \PseudoND( \theta)} & = \Biggr \| \Exs_{Y \sim
  \Ncal(0, 1)} \brackets{Y \tanh \parenth{\frac{Y^{\top}
      \theta}{Z_{\obs, d} - \enorm{ \theta}^2/ d}}} \Biggr \|_2.
\end{align*}
As demonstrated in the proof of Theorem~\ref{thm:em_multivariate}, we
have the equivalence
\begin{align*}
\enorm{ \PseudoND(\theta)} = \Exs \brackets{ V_{1} \tanh
  \parenth{ \frac{ \enorm{ \theta} V_{1}} {Z_{\obs, d} - \enorm{
        \theta}^2/ d}}}
\end{align*}
where $V_{1} \sim \Ncal(0, 1)$. Since the function $x \tanh
\parenth{ \frac{ \enorm{ \theta} x}{a - \enorm{ \theta}^2/ d}}$ is an 
even function in terms of $x$ for any given $a$, we find that
\begin{align*}
\Exs \brackets{ V_{1} \tanh \parenth{ \frac{
    \enorm{ \theta} V_{1}}{Z_{\obs, d} - \enorm{ \theta}^2/ d}}} 
    & = \Exs \brackets{ \abss{ V_{1}} \tanh
    \parenth{ \frac{ \enorm{ \theta} \abss{ V_{1}}}
    {Z_{\obs, d} - \enorm{ \theta}^2/ d}}} \\
    & \leq \Exs \brackets{ \abss{ V_{1}}} 
    = \sqrt{ \frac{ 2}{ \pi}} \nonumber
\end{align*}
where the second inequality is due to the basic inequality $\tanh( x)
\leq 1$ for all $x \in \Rspace$. The inequality in the above display
implies that regardless of the initialization $\theta^0$, we always
have $\enorm{ \PseudoND(\theta)} \leq \sqrt{ 2/ \pi}$, as
claimed.

\section{Wasserstein Distance}
\label{sec:wass_dist}
In Figures~\ref{fig:gaussian_vs_mixture} and \ref{fig:mixture_vs_mixture},
we use EM to estimate all the parameters of the fitted Gaussian mixture
(e.g., the parameters $\braces{w_i, \mu_i, \Sigma_i, i\in[k]}$ if the fitted
mixture
were $\mathcal{G}=\sum_{i=1}^k w_i\Ncal(\mu_i,\Sigma_i)$) and use first-order
Wasserstein distance between the fitted model and the true model to measure
the quality of the estimate.
Here we briefly summarize the definition of the first-order
Wasserstein distance and refer the readers to the book~\cite{Villani-09}
and the paper~\cite{Ho-Nguyen-Ann-16} for more details.
Given two Gaussian mixture distributions of the form
\begin{align*}
  \mathcal{G} = \sum_{i=1}^k w_i\Ncal(\mu_i,\Sigma_i)
  \qtext{and}
  \mathcal{G}' = \sum_{j=1}^{k'} w_j\Ncal(\mu_j',\Sigma_j),
\end{align*}
the first-order Wasserstein distance between the two is given by
\begin{align}
  W_1(\mathcal{G}, \mathcal{G}') = \inf_{q \in \mathcal{Q}} 
  \sum_{i=1}^k\sum_{j=1}^{k'}q_{ij} \parenth{\enorm{\theta_i-\theta'_j}+ \fronorm{\Sigma_i-\Sigma'_j}
  },
  \label{eq:wass_dist}
\end{align}
where $\fronorm{A}$ denotes the Frobenius norm of the matrix $A$ (which
in turn is defined as $\sqrt{\sum_{ij}A_{ij}^2}$).
Moreover, $\mathcal{Q}$ denotes the set of all couplings on $[k] \times
[k']$ such that
\begin{align*}
  q_{ij} \in [0, 1], \qquad
  \sum_{i=1}^k q_{ij} = w'_j\qtext{and}
  \sum_{j=1}^{k'} q_{ij} = w_i\qtext{for all} i\in[k], j\in[k'].
\end{align*}
We note that the optimization problem~\eqref{eq:wass_dist} is a linear program
in the $k \times k'$ dimensional variable $q$ and standard
linear program solvers can be used for solving it.
Also, we remark that here we have abused the notation slightly since the
the definition of the Wasserstein distance above is typically used for the
mixing measures which only depends on the parameters of the Gaussian mixture
(and not the Gaussian density).
Finally, applying definition~\eqref{eq:wass_dist}, we can directly conclude
that for the symmetric fit~\eqref{eq:symmetric_fit}, we have
\begin{align}
  W_1\parenth{\frac{1}{2}\Ncal(\theta, \sigma^2 I_d)+
  \frac{1}{2} \Ncal(-\theta, \sigma^2 I_d), \Ncal(\theta_\star,
  \sigma_{\star}^2I_d) 
  }
  = \enorm{\theta-\theta_\star}  + \sqrt{d} \sqrt{\abss{\sigma^2 - \sigma_\star^2}},
\end{align}
where we have assumed that $\min\braces{\enorm{\theta-\theta_\star},
\enorm{-\theta-\theta_\star}} =
\enorm{\theta-\theta_\star}$.
\end{document}